\documentclass[12pt]{article}
\usepackage{ct-arxiv}
\usepackage[utf8]{inputenc}

\specs{n}{vol (issue)}{year}{}



\title{On quasi-polynomials counting planar tight maps}

\author[1,2]{Jérémie Bouttier\thanks{%
    Supported by the Agence Nationale de la Recherche via the grants
    ANR-18-CE40-0033 ``Dimers'' and ANR-19-CE48-0011 ``Combiné''.}}
\author[2]{Emmanuel Guitter}
\author[3]{Grégory Miermont\thanks{Institut Universitaire de France (IUF)}}

\affil[1]{%
Sorbonne Université and Université Paris Cité, CNRS, IMJ-PRG, F-75005 Paris, France

\email{jeremie.bouttier@imj-prg.fr}%
}

\affil[2]{%
Université Paris-Saclay, CNRS, CEA, Institut de physique théorique,
91191, Gif-sur-Yvette, France

\email{emmanuel.guitter@ipht.fr}
}

\affil[3]{%
ENS de Lyon, UMPA, CNRS UMR 5669, 46 allée d’Italie, 69364 Lyon Cedex 07, France
  
\email{gregory.miermont@ens-lyon.fr}%
}

\newcommand{\fig}[2]{\includegraphics[width=#1\textwidth]{#2}}
\newcommand{\Z}{\mathbb{Z}}
\newcommand{\iv}[1]{\left[#1\right]}

\begin{document}

\maketitle


\begin{abstract}
  A \emph{tight map} is a map with some of its vertices marked, such
  that every vertex of degree $1$ is marked. We give an explicit
  formula for the number $N_{0,n}(d_1,\ldots,d_n)$ of planar tight
  maps with $n$ labeled faces of prescribed degrees $d_1,\ldots,d_n$,
  where a marked vertex is seen as a face of degree $0$. It is a
  quasi-polynomial in $(d_1,\ldots,d_n)$, as shown previously by
  Norbury.  Our derivation is bijective and based on the {\em slice
    decomposition} of planar maps. In the non-bipartite case, we also
  rely on enumeration results for two-type forests.  We discuss the
  connection with the enumeration of non necessarily tight maps. In
  particular, we provide a generalization of Tutte's classical
  slicings formula to all non-bipartite maps.
\end{abstract}



\section{Introduction}\label{sec:introduction-1}

\subsection{Tight maps}\label{sec:tight-maps}

The main purpose of this paper is to study the enumeration problem for a class
of maps, called {\em tight maps}.

\begin{definition}\label{sec:introduction}
  A {\em tight map} is a connected map with some of its vertices
  marked, such that every vertex of degree $1$ is marked. In a tight
  map, the faces as well as the marked vertices are called
  \emph{boundaries}.
\end{definition}

Even though this definition makes sense for maps on arbitrary
surfaces, we will restrict in this paper to the planar case. We refer
to \cite{Schaeffer15} and \cite{triskell} for the standard definitions
and terminology about maps.

\begin{figure}[t]
  \centering
  \includegraphics[width=.9\textwidth]{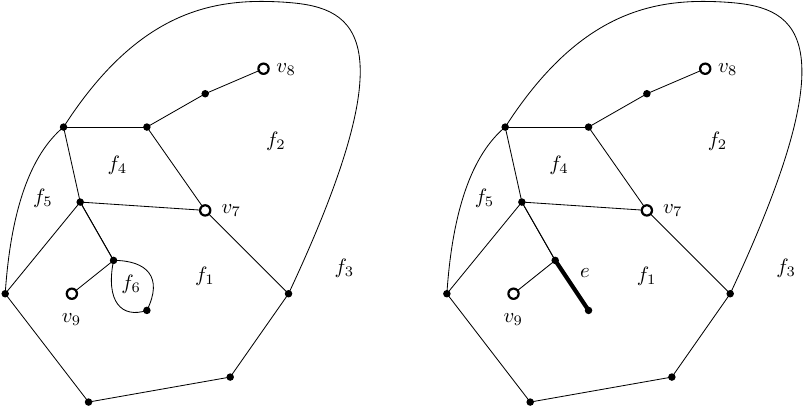}
  \caption{On the left, a tight map with $9$ boundaries: six faces
    $f_1,f_2,f_3,f_4,f_5,f_6$ with respective degrees
    $12, 8, 5, 4, 3, 2$, and three marked vertices $v_7,v_8,v_9$,
    shown in white. On the right: a rooted tight map, whose opening
    operation gives the tight map on the left. As this figure
    demonstrates, a vertex incident to the root in a rooted tight map
    can be of degree $1$ without being necessarily marked. These maps
    can be also seen as pointed (respectively pointed rooted) tight
    maps, for example by distinguishing the marked vertex $v_9$. }
  \label{fig:sampletightmap}
\end{figure}

Usually, we will endow a tight map with some extra structure, in
particular by labeling its faces and some or all of
its marked vertices, or distinguishing one marked vertex. 
For instance, we will call {\em pointed tight map} a tight map
with one distinguished marked vertex. We will adopt a slightly unusual
notion of {\em rooted tight map} compared to the well-established notion of
rooting of a map: in particular, the root will be an unoriented edge.  
If $e$ is a distinguished edge in a map $\mathbf{m}$, there is a natural
opening operation $O(\mathbf{m},e)$ consisting in cutting along the edge, thereby
creating a face of degree $2$. If $\mathbf{m}$ has marked or labeled
elements (vertices or faces), then $O(\mathbf{m}, e)$ naturally
inherits these elements.  
We say that a map $\mathbf{m}$ with some of its vertices marked and a distinguished
edge $e$ is a {\em rooted tight map} if $O(\mathbf{m},e)$ is a tight map. Note that
$\mathbf{m}$ may not be a tight map itself, as the distinguished edge may be
incident to an unmarked vertex of degree $1$.
Finally, a {\em pointed rooted tight map} is a rooted tight map
with one distinguished marked vertex. See Figure
\ref{fig:sampletightmap} for examples of tight maps.

We define the \emph{length} of a boundary in a tight map as being
equal to its degree for a face, and to zero for a marked vertex. In
other words, we interpret the marked vertices as boundaries of length
$0$.

\begin{figure}[!]
  \centering
  \includegraphics[width=.7\textwidth]{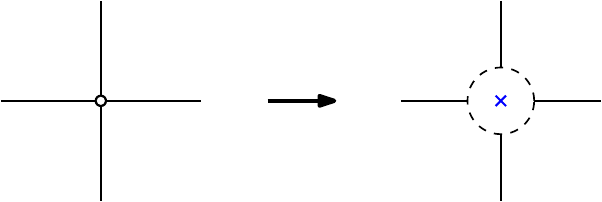}
  \caption{Blowing a boundary-vertex into a face of degree $0$ with
    one puncture. }
  \label{fig:boundaryvertex}
\end{figure}

The terminology of tight maps comes from \cite{triskell}. Let us
discuss it in some detail. In a general map $\mathbf{m}$, drawn on a
surface $S$, and with some faces and vertices marked, let us call the
marked elements the {\em boundaries}, and the unmarked elements the
{\em internal} faces and vertices. We let $S'$ be the space obtained
from $S$ by removing one point (i.e.\ creating a puncture) inside each
of the boundaries of $\mathbf{m}$.  A boundary of $\mathbf{m}$ is
called {\em tight} if its contour path has minimal possible length
among all paths in the map that are freely homotopic to it in
$S'$. When the map has marked vertices, then all such vertices are
automatically tight boundaries, and the previous minimality condition
on the contours of the boundary-faces must be understood in a slightly
modified map obtained by blowing every marked vertex of degree $k$
into a $k$-cycle with edges of ``length $0$'', hence creating a new
face in the map, which we view as having degree $0$, with one
puncture. See Figure \ref{fig:boundaryvertex} for an example.

This being said, it is straightforward to see that a map is tight
according to Definition \ref{sec:introduction}, if and only if it is a
map with no internal faces (i.e.\ all its faces are marked as
boundaries), whose boundaries are all tight. Indeed, starting from a
tight map, we see that all its boundaries are necessarily tight: in
the modified map with every marked vertex blown into a degree-$0$
face, the contour of a given face is in fact the unique
non-backtracking path of edges in its free homotopy class in the
punctured surface $S'$. Conversely, a map which is not tight contains
an unmarked vertex $v$ of degree $1$. The contour of the face $f$
incident to $v$ can be deformed into a strictly shorter path by
shortcutting the edge incident to $v$, meaning that $f$ is not a tight
boundary.

In fact, we view the results of the present paper as a first step
towards a better (in particular, bijective) understanding of the
counting problem for general maps with tight boundaries 
and possibly with internal faces. This problem was addressed in the
case of planar maps with three boundaries in \cite{triskell}, and
remains a challenge in more complex topologies.  

\subsection{Lattice count polynomials}\label{sec:counting-tight-maps}

For any choice of nonnegative integers $d_1,d_2,\ldots,d_n$ not all
equal to $0$ and for any nonnegative integer $g$, we let
$N_{g,n}(d_1,d_2,\ldots,d_n)$ be the number of tight maps of genus
$g$, with $n$ labeled boundaries of respective lengths
$d_1,d_2,\ldots,d_n$, and where each map is weighted by its inverse
number of automorphisms. In fact, the latter number is always $1$ as
soon as $n\geq 3$, while the only tight maps with two boundaries
having a non-trivial automorphism group are (genus $0$) cycles of
length $p\geq 2$, with automorphism group $\Z/p\Z$ (such maps are thus
weighted by $1/p$). For $n=1$, the only planar ($g=0$) example of a tight map
is the (marked) vertex-map, which is excluded from the discussion since its
only boundary has length $0$, but in higher genera, there are many
tight maps with one face, and these can have a non-trivial automorphism group. 

The numbers $N_{g,n}(d_1,\ldots,d_n)$ have been extensively studied in
particular by Norbury and Do
\cite{Norbury2010,Norbury2013,DoNo11}\footnote{Here, we should warn
  the reader that these references use the notion of fatgraphs, which
  is different but equivalent to the language of maps used in this
  paper.}, in the broader context of the study of invariants of
spectral curves appearing in Eynard and Orantin's topological
recursion \cite{Eynard2007}.  Norbury proved that $N_{g,n}(d_1,\ldots,\allowbreak d_n)$ is a
quasi-polynomial in the variables $d_1^2,\ldots,d_n^2$, depending on
their parities. This means that for every $k\in \{0,1,\ldots,n\}$,
there exists a polynomial $\mathsf{N}^{(k)}_{g,n}(x_1,\ldots,x_n)$ in
the variables $x_1^2,\ldots,x_n^2$, symmetric under permutations of
the first $k$ variables and of the last $n-k$ variables, such that, if
the numbers $d_1,\ldots,d_k$ are odd and the numbers
$d_{k+1},\ldots,d_n$ are even, then
\begin{equation}
  \label{eq:1}
  N_{g,n}(d_1,\ldots,d_n)=\mathsf{N}^{(k)}_{g,n}(d_1,\ldots,d_n)\, .
\end{equation}
For $k$ odd, these polynomials are equal to $0$. In fact, the first
two papers mentioned above assume that $d_1,\ldots,d_n$ are all
non-zero, while the third paper considers the general case where some,
but not all $d_i$ may vanish. Definition 2.7 in \cite{DoNo11} is
indeed equivalent to our definition of tight maps, while Proposition
2.8 therein proves that the extension of the quasipolynomials to some
zero values do solve the enumeration problem of tight maps with marked
vertices.  Norbury \cite{Norbury2010} also proves that evaluating the
polynomials at $(0,0,\ldots,0)$ gives interesting geometric
information, although the combinatorial meaning of this evaluation is
not clear. Note that the theory of enumeration of integer points in
polytopes implies that $N_{g,n}(d_1,\ldots,d_n)$ is a piecewise
quasi-polynomial in $d_1,\ldots,d_n$, see for instance the discussion
around~\cite[Proposition~4]{Budd2020a} (in the case $b=0$). Therefore,
it is surprising that $N_{g,n}(d_1,\ldots,d_n)$ is actually a genuine
quasi-polynomial, furthermore in the squared variables.

The approach taken in \cite{Norbury2010,Norbury2013,DoNo11} is to
prove the wanted properties using recursions for these polynomials,
called {\em lattice count polynomials}, that also allows one to
effectively compute them. These recursions are in turn consequences of
combinatorial recursion relations with a geometric flavor, similar to
Tutte's equations used for instance in~\cite{Tutte1962}, and to the
topological recursion originating in Eynard and Orantin's
work~\cite{Eynard2007}.

In this paper, focusing on the planar case
$g=0$, our main goal in to show how one can obtain the above quasipolynomiality results by \emph{bijective}
techniques, which in passing yield new explicit formulas for the lattice count
polynomials. We will use two different strategies: the first one, discussed in Section \ref{sec:slicings},
is based on a \emph{substitution} approach using as an input Tutte's classical
{\em slicings formula}~\cite{Tutte1962}. This formula holds however only for planar maps which 
are \emph{bipartite} or \emph{quasi-bipartite}, namely with a number of faces
of odd degree equal to $0$ or $2$ respectively. As a consequence, the substitution approach is limited  
to the enumeration of planar tight bipartite and quasi-bipartite maps, corresponding to respectively $k=0$ and $k=2$ in \eqref{eq:1}. 
The second, purely bijective, strategy is based on the so-called {\em slice decomposition} of maps introduced in \cite{hankel}, 
and its extensions developed in \cite{irredmaps,BouttierHDR}. 
We will first discuss it in Section~\ref{sec:bijective} in the easier case of planar tight bipartite 
and quasi-bipartite maps, and then extend it in Section~\ref{sec:extens-non-bipart}
to the general case of planar tight maps with an arbitrary number $k$ of faces of odd degree.
Using then the substitution approach backwards, our general expression for
tight maps allows us to extend Tutte's slicings formula to non necessarily tight maps with an arbitrary
number of faces of odd degree.

\medskip

The paper is organized as follows. Section \ref{sec:main-results-2}
provides a self-contained presentation of our main results.
Section~\ref{sec:main-results}
deals with the simpler case of bipartite and quasi-bipartite
tight maps: after introducing and studying in
Section~\ref{sec:polynomials} the required basic univariate and
multivariable polynomials, we state our main theorems which connect
these polynomials to the numbers of planar tight bipartite maps in
Section~\ref{sec:mainresult} (Theorem~\ref{thm:maintheorem}) and of
planar tight quasi-bipartite maps in Section~\ref{sec:quasibipartite}
(Theorem~\ref{thm:quasibiptheorem}). We then state our enumeration
result for general tight maps in Section~\ref{sec:quasipolmain}
(Theorem~\ref{thm:quasimain}) after introducing the appropriate
univariate and multivariate quasi-polynomials.
Section~\ref{sec:slicings} discusses the connection with the
enumeration of non necessarily tight maps by the substitution
approach: Section \ref{sec:equiv-with-tutt} is devoted to the derivation of Theorems
\ref{thm:maintheorem} and \ref{thm:quasibiptheorem} from Tutte's
slicings formula, and Section \ref{sec:non-bipart-slic} uses this
approach backwards to obtain from Theorem~\ref{thm:quasimain} an
extension of the slicings formula to maps with an arbitrary number of
faces of odd degree, see Theorem~\ref{thm:extendedslicings}. We then
discuss in Sections~\ref{sec:bijective} and
\ref{sec:extens-non-bipart} the bijective
approach based on the slice decomposition of maps. Section~\ref{sec:bijective}
concentrates again on the simpler bipartite and quasi-bipartite cases,
discussing first tight maps with a single face
(Section~\ref{sec:oneface}), tight maps with two faces
(Section~\ref{sec:twofaces}), pointed rooted tight maps in connection
with tight slices (Section~\ref{sec:pointedrooted}) and finally tight
maps which are neither pointed nor rooted
(Section~\ref{sec:generalbip}) using the slice decomposition
of annular maps. All these bijective results are then extended to the
non-bipartite or quasi-bipartite case in
Section~\ref{sec:extens-non-bipart}, which requires the preliminary
enumeration of so-called petal trees (Section~\ref{sec:petaltrees}),
petal necklaces (Section~\ref{sec:petalnecklaces}) and non-bipartite
slices (Section~\ref{sec:generalslices}). Our most general enumeration
result for planar tight maps with arbitrary prescribed boundary lengths is
given in Section~\ref{sec:generalnonbip} by
Theorem~\ref{thm:themaintheorem}, which presents a single formula
encompassing Theorems~\ref{thm:maintheorem}, \ref{thm:quasibiptheorem}
and \ref{thm:quasimain}.  We gather our concluding remarks in
Section~\ref{sec:conclusion}, while a few appendices detail the
derivation of some technical results.

\section{Main results}\label{sec:main-results-2}

\subsection{Polynomials counting planar tight bipartite or quasi-bipartite maps}\label{sec:main-results}

In this section, we provide explicit expressions for the lattice count
polynomials $\mathsf{N}^{(0)}_{0,n}$ and $\mathsf{N}^{(2)}_{0,n}$,
which correspond to planar tight bipartite and quasi-bipartite maps,
respectively.

\subsubsection{Definition and properties of the polynomials}
\label{sec:polynomials}

Let us start by introducing families of polynomials which appear in
the explicit expression of $\mathsf{N}^{(0)}_{0,n}$. Here, we
concentrate on the very definitions of these polynomials and on their
resulting algebraic properties. The connection with tight map
enumeration will be discussed in Section~\ref{sec:mainresult}.

\paragraph{Basic univariate polynomials.}

Our first basic polynomials are functions of a single variable $m$ and are defined as follows:
for any integer $k\geq 0$, we set
\begin{equation}
  \label{eq:defpkqk}
  \begin{split}
   & p_k(m):=\frac{1}{(k!)^2}\, \prod_{i=1}^k\left(m^2-i^2\right) = \binom{m-1}{k} \binom{m+k}{k} \\
   & q_k(m):=\frac{1}{(k!)^2}\, \prod_{i=0}^{k-1}\left(m^2-i^2\right) = \binom{m}{k} \binom{m+k-1}{k} \\
 \end{split}
\end{equation}
with the usual convention $p_0(m)=q_0(m)=1$ for the empty product,
and with $\binom{x}{k}=x(x-1)\cdots(x-k+1)/k!$ viewed as a polynomial in $x$.
 
Clearly, $p_k$ and $q_k$ are polynomials of degree $k$ in the variable
$m^2$, and $p_k(m)$ and $q_k(m)$ are integers if $m$ is an
integer. The two families of polynomials are linked by the relation
\begin{equation}
\label{eq:relpkqk}
q_k(m)=p_k(m)+p_{k-1}(m),\qquad k\geq 0 
\end{equation}
with the convention $p_{-1}:=0$.  A combinatorial interpretation of
this relation based on the enumeration of tight maps with a single
face will be given in Section~\ref{sec:mainresult}. We also record the
identities
\begin{align}
  \label{eq:pkunivstr}
  (k+1) p_{k+1}(m) &= (m-k-1) p_k(m) + \sum_{j=1}^{m-1} (2j) p_k(j) \\
  \label{eq:qkunivstr}
  (k+1) q_{k+1}(m) &= (m-k) q_k(m) + \sum_{j=1}^{m-1} (2j) q_k(j)
\end{align}
which are valid for $m$ a positive integer, and which may be checked by induction.

\paragraph{Multivariate polynomials.}

The above univariate polynomials may be extended to multivariate polynomials, functions of $n$ variables
$m_1,m_2,\ldots,m_n$ as follows:
for any integer $k\geq 0$ and any integer $n\geq 1$, we set
 \begin{equation}
 \label{eq:defpkqkmultivariate}
  \begin{split}
    & p_k(m_1,m_2,\ldots,m_n):=\sum_{k_1,k_2,\ldots,k_n\geq 0\atop k_1+k_2+\cdots+k_n=k}p_{k_1}(m_1)q_{k_2}(m_2)\cdots q_{k_n}(m_n) ,  \\
    & q_k(m_1,m_2,\ldots,m_n):=\sum_{k_1,k_2,\ldots,k_n\geq 0\atop k_1+k_2+\cdots+k_n=k}q_{k_1}(m_1)q_{k_2}(m_2)\cdots q_{k_n}(m_n) .\\
  \end{split}
\end{equation}
(In the right-hand side of the first line, all factors except the first one are $q_{k_i}$'s.)

Clearly, $p_k$ and $q_k$ are polynomials of degree $k$ in the variables $m_1^2,m_2^2,\ldots,m_n^2$.
Note that the notation is consistent for $n=1$ with
\eqref{eq:defpkqk}. 
From the identity $p_k(1)=\delta_{k,0}$, we get the identification:
\begin{equation}
\label{eq:addingroot}
p_k(1,m_1,m_2,\ldots,m_n)=q_k(m_1,m_2,\ldots,m_n)
\end{equation}
and, by~\eqref{eq:relpkqk}, we find that $q_k=p_k+p_{k-1}$ for any
number of variables.

For bookkeeping purposes, let us record the explicit expressions of
$p_k$ and $q_k$ for $k=1,2,3$: 
\begin{equation}
  \label{eq:firstfewpk}
  \begin{split}
    p_1(m_1,\ldots,m_n) &= \left(\sum_{i=1}^n m_i^2 \right) - 1, \\
    p_2(m_1,\ldots,m_n) &= \frac14 \left(\sum_{i=1}^n m_i^4 \right) +  \sum_{i<j} m_i^2 m_j^2  - \frac54 \left(\sum_{i=1}^n m_i^2 \right) + 1, \\
     p_3(m_1,\ldots,m_n) &= \frac1{36} \left(\sum_{i=1}^n m_i^6 \right) + \frac14 \left( \sum_{i \neq j} m_i^4 m_j^2 \right) + \sum_{i<j<h} m_i^2 m_j^2 m_h^2 \\
    & \qquad 
    - \frac7{18} \left(\sum_{i=1}^n m_i^4 \right) - \frac32 \left(
      \sum_{i<j} m_i^2 m_j^2 \right) + \frac{49}{36}
    \left(\sum_{i=1}^n m_i^2 \right) - 1\, ,\\
     q_1(m_1,\ldots,m_n) &= \sum_{i=1}^n m_i^2 , \\
    q_2(m_1,\ldots,m_n) &= \frac14 \left(\sum_{i=1}^n m_i^4 \right) +  \sum_{i<j} m_i^2 m_j^2  - \frac14 \left(\sum_{i=1}^n m_i^2 \right) , \\
     q_3(m_1,\ldots,m_n) &= \frac1{36} \left(\sum_{i=1}^n m_i^6 \right) + \frac14 \left( \sum_{i \neq j} m_i^4 m_j^2 \right) + \sum_{i<j<h} m_i^2 m_j^2 m_h^2 \\
    & \qquad 
    - \frac5{36} \left(\sum_{i=1}^n m_i^4 \right) - \frac12 \left( \sum_{i<j} m_i^2 m_j^2 \right) + \frac19 \left(\sum_{i=1}^n m_i^2 \right) .
  \end{split}
\end{equation}

\begin{proposition}
  \label{prop:pqsym}
  For any integer $k \geq 0$, $p_k$ and $q_k$ are symmetric
  functions. In other words, for any integer $n \geq 1$,
  $p_k(m_1,m_2,\ldots,m_n)$ and $q_k(m_1,m_2,\ldots,m_n)$ are
  symmetric polynomials in $m_1,m_2,\ldots,m_n$ which satisfy the
  consistency relation
  \begin{equation}
    \label{eq:pkqkwithzeroes}
    \begin{split}
      &p_k(m_1,m_2,\ldots,m_n,0)= p_k(m_1,m_2,\ldots,m_n),\\
      & q_k(m_1,m_2,\ldots,m_n,0)=q_k(m_1,m_2,\ldots,m_n).\\
\end{split}
\end{equation}
\end{proposition}
\begin{proof}
The symmetry of $q_k(m_1,m_2,\ldots,m_n)$ is apparent from its very definition in \eqref{eq:defpkqkmultivariate}.
As for $p_k(m_1,m_2,\ldots,m_n)$, its symmetry is also made apparent from the following alternative and manifestly symmetric expression:
\begin{equation}
\label{eq:alternativepk}
p_k(m_1,m_2,\ldots,m_n)=\sum_{k_0,k_1,k_2,\ldots,k_n\geq 0\atop k_0+k_1+k_2+\cdots+k_n=k}{n-1\choose k_0} p_{k_1}(m_1)p_{k_2}(m_2)\cdots p_{k_n}(m_n).
\end{equation}
To get this latter expression, we use again the relation \eqref{eq:relpkqk} to write, in the expression \eqref{eq:defpkqkmultivariate}
for $p_k(m_1,m_2,\ldots,m_n)$, each $q_{k_i}$ for $i=2$ to $n$ as the
sum of $p_{k_i}$ and $p_{k_i-1}$ and distribute the two terms in the product so as to get a sum of terms of the form
$p_{k'_1}(m_1)p_{k'_2}(m_2)\cdots p_{k'_n}(m_n)$ with summation variables $k'_1=k_1$ and $k'_i=k_i$ or $k_i-1$ for $i\geq 2$. The number of  
terms in the sum having exactly $k_0$ indices $i$ for which $k'_i=k_i-1$ is ${n-1\choose k_0}$ and, for such terms, 
the sum rule $k_1+k_2+\cdots+k_n=n$ in \eqref{eq:defpkqkmultivariate}
becomes  $k_0+k'_1+k'_2+\cdots+k'_n=n$. This leads to \eqref{eq:alternativepk} upon renaming the summation variable $k'_i$ as $k_i$.

As for the consistency relation~\eqref{eq:pkqkwithzeroes}, it is
a direct consequence of the identity $q_k(0)=\delta_{k,0}$ (here and in the following, we will always implicitly assume that $k$ is a non-negative integer).
\end{proof}

Finally, let us state some recursion relations obeyed by the $p_k$'s,
which we call the dilaton and string equations as we shall see later
that they correspond to the recursions obtained in~\cite{Norbury2013}
in the bipartite case.

\begin{proposition}[Dilaton and string equations]
  \label{prop:pqdilstr}
 We have the dilaton equation
  \begin{equation}
    \label{eq:pdil}
    p_k(m_1,\ldots,m_n,1) - p_k(m_1,\ldots,m_n,0) = p_{k-1}(m_1,\ldots,m_n)
  \end{equation}
  and the string equation, valid for non-negative integer
  $m_1,\ldots,m_n$: 
  \begin{multline}
    \label{eq:pstr}
    (k+1) p_{k+1}(m_1,\ldots,m_n) = (m_1 + \cdots + m_n - k - 1) p_k(m_1,\ldots,m_n) + \\ \sum_{i=1}^n \sum_{j=1}^{m_i-1} (2j) p_k(m_1,\ldots,m_{i-1},j,m_{i+1},\ldots,m_n)
  \end{multline}
\end{proposition}

\begin{proof}
  From~\eqref{eq:addingroot} and Proposition~\ref{prop:pqsym}, the
  dilaton equation boils down to the relation $q_k=p_k+p_{k-1}$ noted
  above. The string equation is nothing but the multivariate extension
  of~\eqref{eq:pkunivstr}, and is obtained by a linear combination of
  it and~\eqref{eq:qkunivstr} (precisely, we take~\eqref{eq:pkunivstr}
  at $k=k_1$ and $m=m_1$ times $q_{k_2}(m_2) \cdots q_{k_n}(m_n)$ and
  add, for $i=2,\ldots,n$, the relation~\eqref{eq:qkunivstr} at
  $k=k_i$ and $m=m_i$ times $p_{k_1}(m_1) \cdots q_{k_n}(m_n)$ with
  the factor $q_{k_i}(m_i)$ omitted).
\end{proof}

\subsubsection{Enumeration results in the bipartite case}
\label{sec:mainresult}

We are now ready to state our first enumerative result:
\begin{theorem}
\label{thm:maintheorem}
For $n\geq 3$ and for non-negative integers $m_1,m_2,\ldots,m_n$ not
all equal to zero, the 
number $N_{0,n}(2m_1,2m_2,\ldots,2m_n)$ of planar tight bipartite maps with $n$ boundaries labeled from $1$ to $n$ 
with respective lengths $2m_1,2m_2,\ldots,2m_n$ is given by the polynomial
\begin{equation}
\label{eq:maintheorem}
\mathsf{N}^{(0)}_{0,n}(2m_1,2m_2,\ldots,2m_n)=(n-3)!\, p_{n-3}(m_1,m_2,\ldots,m_n).
\end{equation}
\end{theorem}

For $n=4,5,6$, the expression for $\mathsf{N}^{(0)}_{0,n}$ that we
obtain from~\eqref{eq:firstfewpk} is in agreement with the polynomials
given in \cite[Table~1]{Budd2020a} for $g=0$ and $b=0$, as expected.
Note that the constant term of $\mathsf{N}^{(0)}_{0,n}$, obtained by
setting all the $m_i$'s to zero, is equal to
$(n-3)!p_{n-3}(0,\ldots,0)=(-1)^{n-3}(n-3)!$, and this quantity was
interpreted in~\cite{Norbury2010} as the orbifold Euler characteristic
of the moduli space $\mathcal{M}_{0,n}$.

By combining Proposition \ref{prop:pqdilstr} with Theorem
\ref{thm:maintheorem}, we recover the string and dilaton equations
found by Norbury in \cite{Norbury2010} in the planar bipartite
case. Note that, in this reference, the string equation corresponds to
the addition of a face of degree $2$, while \eqref{eq:pstr} for
$k=n-3$ may be interpreted as the addition of a vertex. Norbury's
original equation
can however be recovered by combining it with the dilaton
equation. Note finally that Proposition \ref{prop:pqdilstr} holds more
generally for any $k$. When $k\geq n-3$ we can naturally interpret it
in terms of adding new marked vertices, but the combinatorial meaning of the
polynomial $p_k(m_1,\ldots,m_n)$ for $k<n-3$ is more elusive. 

A first derivation of Equation~\eqref{eq:maintheorem} will be
presented in Section \ref{sec:slicings} below by showing that, up to
some appropriate transformation accounting for the tight nature of the
maps, it is actually equivalent to Tutte's celebrated slicings
enumeration formula~\cite{Tutte1962}. We shall then present in Section
\ref{sec:bijective} a direct bijective proof of
Theorem~\ref{thm:maintheorem} upon using some canonical \emph{slice
  decomposition} of the maps at hand
\cite{hankel,irredmaps,BouttierHDR}. As it appears, it will be
convenient for that purpose to proceed gradually and first derive
\eqref{eq:maintheorem} for a number of specialized cases
(Propositions~\ref{prop:pkinterpret}, \ref{prop:qkinterpret},
\ref{prop:pk2interpret} and \ref{prop:pointedrooted} below) before
addressing the result in all generality.

\paragraph{Maps with one face.}
Taking $n=k+3$ in \eqref{eq:maintheorem} with $m_1=m\neq 0$ and $m_2=\cdots=m_{k+3}=0$, we get
\begin{equation}
\label{eq:oneface}
N_{0,k+3}(2m,\underbrace{ 0,\ldots,0}_{k+2})=k!\,
p_k(m,\underbrace{0,\ldots,0}_{k+2})=k!\, p_k(m),\qquad k\geq 0\, ,
\end{equation}
which already appeared in \cite[Corollary
5.6]{Norbury2013}\footnote{Incidentally, we note that our main results (Theorems
\ref{thm:maintheorem}, \ref{thm:quasibiptheorem} and \ref{thm:quasimain}) answer the question
raised in the paragraph just before this corollary
 about finding a general formula for $N_{0,n}$.}. 
Upon dividing by $k!$, which amounts to considering that all but two
of the marked vertices are unlabeled, we obtain the following
combinatorial interpretation of $p_k(m)$:
\begin{proposition}
\label{prop:pkinterpret}
For $k\geq 0$ and $m\geq 1$, $p_k(m)$ is the number of planar tight bipartite maps with one face of degree $2m$ and 
$k+2$ distinct marked vertices, two of them distinguished and labeled, say as vertex $1$ and vertex $2$, and the remaining $k$ unlabeled. 
\end{proposition}
\noindent Note that a planar map with a single face of degree $2m$ is nothing but a plane tree with $m$ edges. It is tight if and only if
all its leaves are marked.

\medskip
Similarly, taking $n=k+3$ in \eqref{eq:maintheorem} with $m_1=m$ and $m_2=1$ and $m_3=\cdots=m_{k+3}=0$, we now 
get
\begin{equation}
\label{eq:onefacebis}
N_{0,k+3}(2m,2,\underbrace{ 0,\ldots,0}_{k+1})=k!\, p_k(m,1,\underbrace{0,\ldots,0}_{k+1})=k!\, q_k(m,\underbrace{0,\ldots,0}_{k+1})=k!\, q_k(m),
\end{equation}
where we used \eqref{eq:addingroot} (and the symmetry in exchanging the variables) to switch from $p_k$ to $q_k$.
Dividing by $k!$ and viewing the face of degree $2$ as a split root edge, as discussed in the introduction, we obtain a combinatorial interpretation of $q_k(m)$:
\begin{proposition}
\label{prop:qkinterpret}
For $k\geq 0$ and $m\geq 1$, $q_k(m)$ is the number of pointed rooted
planar tight bipartite maps with one face of degree $2m$ and $k$
additional unlabeled marked vertices (distinct from each other and
from the pointed vertex).
\end{proposition}

Proofs of Propositions~\ref{prop:pkinterpret} and
\ref{prop:qkinterpret} will be presented in Section~\ref{sec:oneface}
by a direct enumeration of the trees at hand.  From the above
interpretations of $p_k(m)$ and $q_k(m)$, we may now understand the
identity \eqref{eq:relpkqk}, for integer values of $m$, in a
combinatorial way using the map language. Indeed, for each tree
enumerated by $q_k(m)$, we may transfer the marking of its root edge
into a marking of that of its endpoints further away from the pointed
vertex. Let us for clarity label the newly marked vertex as vertex $2$
and the pointed vertex as vertex $1$. Note that vertices $1$ and $2$
are necessarily distinct by construction, but that the vertex $2$ may
very well coincide with one of the $k$ additional marked vertices in
the map enumerated by $q_k(m)$.  The marking transformation is clearly
reversible, the root edge being recovered as the only edge incident to
vertex $2$ that belongs to the branch (that is, the unique simple path) from vertex $2$ to vertex
$1$. We may thus interpret $q_k(m)$ as counting plane trees with $m$
edges, and with two distinct marked vertices $1$ and $2$ and $k$ other
marked vertices distinct from each other and from the vertex $1$. This
yields a map enumerated by $p_k(m)$ when none of the $k$ marked
vertices coincide with the vertex $2$---note that this may happen even
if the vertex $2$ is a leaf since, as an endpoint of the root edge, it
needs not being marked in the map enumerated by $q_k(m)$. Otherwise,
it yields a map enumerated by $p_{k-1}(m)$ by ignoring the
``redundant'' additional marking of vertex $2$. This yields the
desired relation \eqref{eq:relpkqk}.

\paragraph{Maps with two faces.}
Taking $n=k+3$ in \eqref{eq:maintheorem} with $m_1,m_2 \geq 1$ and $m_3=\cdots=m_{k+3}=0$, we get
\begin{equation}
\label{eq:twofaces}
N_{0,k+3}(2m_1,2m_2,\underbrace{ 0,\ldots,0}_{k+1})=k!\, p_k(m_1,m_2,\underbrace{0,\ldots,0}_{k+1})=k!\, p_k(m_1,m_2).
\end{equation}
Upon dividing by $k!$ we get a combinatorial interpretation of $p_k(m_1,m_2)$:
\begin{proposition}
\label{prop:pk2interpret}
For $k\geq 0$ and $m_1,m_2 \geq 1$, $p_k(m_1,m_2)$ is the number of planar tight bipartite maps with two faces of respective degrees $2m_1,2m_2$ and $k+1$ distinct marked vertices, one of them distinguished and labeled, say as vertex $1$, and the remaining $k$ unlabeled. 
\end{proposition}
A direct bijective proof of this proposition will be presented in
Section~\ref{sec:twofaces}. More generally, and although we will not
use it in our bijective proof in Section \ref{sec:bijective}, we note the relation,
valid for integers $m_1,m_2,m_3$ not all equal to $0$
\begin{equation}
  \label{eq:3}
  N_{0,k+3}(2m_1,2m_2,2m_3, \underbrace{0,\ldots,0}_{k})=k!\, p_k(m_1,m_2,m_3)
\end{equation}
so that $p_{k}(m_1,m_2,m_3)$ counts planar tight bipartite maps with
three labeled boundaries with lengths $2m_1,2m_2,2m_3$ and $k$ unlabeled
marked vertices. It is relatively straightforward to adapt the
bijective proof of Proposition~\ref{prop:pk2interpret} to
prove~\eqref{eq:3} directly, we leave it as an exercise to the reader.

\paragraph{Pointed rooted maps.}
Taking \eqref{eq:maintheorem} with $n\to n+2$ and specializing to $m_{n+1}=2$ and $m_{n+2}=0$, we get:
\begin{proposition}
\label{prop:pointedrooted}
  For $n \geq 1$, the number of pointed rooted planar tight bipartite
  maps with $n$ labeled boundaries of respective lengths
  $2m_1,\ldots,2m_n$ (in addition to the marked vertex and to the root
  edge) is given by
\begin{equation}
\label{eq:pointedrooted}
\mathsf{N}_{0,n+2}^{(0)},(2m_1,2m_2,\dots,2m_n,2,0)=(n-1)!\, q_{n-1}(m_1,m_2,\ldots,m_n).
\end{equation}
\end{proposition}
\begin{proof}
This is a direct application of formula \eqref{eq:maintheorem}, together with the identity $p_{n-1}(m_1,m_2,\allowbreak \ldots,m_n,1,0)=q_{n-1}(m_1,m_2,\ldots,m_n)$
from \eqref{eq:addingroot}.
\end{proof}
The expression \eqref{eq:pointedrooted} will be proved bijectively in Section~\ref{sec:pointedrooted} by a direct decomposition of the maps into slices.

\subsubsection{Enumeration results in the quasi-bipartite case}
\label{sec:quasibipartite}

Recall that a planar quasi-bipartite map is a planar map whose all
faces but two have even degree. To give the explicit expression of the
corresponding lattice count polynomial $\mathsf{N}^{(2)}_{0,n}$, we
need to introduce the following new family of univariate polynomials:
for any integer $k\geq 0$, we set
 \begin{equation}
 \label{eq:defptildek}
 \tilde{p}_k(m):=\frac{1}{(k!)^2}\, \prod_{i=1}^k\left(m^2-\left(i-\frac{1}{2}\right)^2\right) = \binom{m-\frac12}{k}\binom{m+k-\frac12}{k}\\
    \end{equation}
with the convention $\tilde{p}_0(m)=1$.
Note that $\tilde{p}_k$ is again a polynomial of degree $k$ in $m^2$
and that $\tilde{p}_k(m)$ is an integer if $m$ is a half-integer. It
satisfies the following counterpart of~\eqref{eq:pkunivstr}
\begin{equation}
  \label{eq:tildepkunivstr}
  (k+1) \tilde{p}_{k+1}(m) = \left(m-k-\tfrac12\right) \tilde{p}_k(m) + \sum_{0<j<m} (2j) \tilde{p}_k(j)
\end{equation}
where it is understood that $m$ and $j$ are now half-integers. 

\medskip
The multivariate extension of $\tilde{p}_k$ is then defined
for any integer $k\geq 0$ and any integer $n\geq 2$ as
 \begin{equation}
 \label{eq:defptildekmultivariate}
 \tilde{p}_k(m_1,m_2;m_3,\ldots,m_n):=\sum_{k_1,k_2,\ldots,k_n\geq 0\atop k_1+k_2+\cdots+k_n=k}\tilde{p}_{k_1}(m_1)
 \tilde{p}_{k_2}(m_2)q_{k_3}(m_3)\cdots q_{k_n}(m_n).
 \end{equation}
\noindent Again, we may append an arbitrary number of $0$'s to the arguments of $\tilde{p}_k$ without changing its 
value. Note that $\tilde{p}_k$ is in general \emph{not symmetric} in all its variables, but only in $m_1$ and $m_2$ on the one hand, and in $m_3,\ldots,m_n$ on the other hand. Note also that $\tilde{p}_k(m,1/2)=\tilde{p}_k(m)$.
The quasi-bipartite analog of Theorem \ref{thm:maintheorem} is then:
\begin{theorem}
\label{thm:quasibiptheorem}
For $n\geq 3$, for $m_1,m_2 \in \Z_{\geq 0} + \frac12$ and
$m_3,\ldots,m_n \in \Z_{\geq 0}$, the number
$N_{0,n}(2m_1,2m_2,2m_3,\ldots,2m_n)$ of planar tight quasi-bipartite
maps with $n$ boundaries labeled from $1$ to $n$ with respective
lengths $2m_1,2m_2,2m_3,\ldots,2m_n$ is given by
\begin{equation}
\label{eq:quasibiptheorem}
\mathsf{N}_{0,n}^{(2)}(2m_1,2m_2,2m_3,\ldots,2m_n)=(n-3)!\, \tilde{p}_{n-3}(m_1,m_2;m_3,\ldots,m_n).
\end{equation}
\end{theorem}
\noindent In particular, for $n=k+3$ and $m_3=\cdots=m_{k+3}=0$, we get, for $m_1, m_2$ half-integers:
\begin{equation}
\label{eq:twofacesquasi}
N_{0,k+3}(2m_1,2m_2,\underbrace{ 0,\ldots,0}_{k+1})=k!\, \tilde{p}_k(m_1,m_2;\underbrace{0,\ldots,0}_{k+1})=k!\, \tilde{p}_k(m_1,m_2).
\end{equation}
Upon dividing by $k!$, we obtain a combinatorial interpretation of $\tilde{p}_k(m_1,m_2)$:
\begin{proposition}
\label{prop:tildepk2interpret}
For $k\geq 0$ and $m_1,m_2 \in \Z_{\geq 0}+\frac12$,
$\tilde{p}_k(m_1,m_2)$ is the number of planar tight quasi-bipartite
maps with two faces of odd degrees $2m_1,2m_2$ and $k+1$ distinct
marked vertices, one of them distinguished and labeled, say as vertex
$1$, and the remaining $k$ unlabeled.
\end{proposition}
\noindent Setting $m_2=1/2$, we obtain a combinatorial interpretation
of the univariate polynomial $\tilde{p}_k(m)=\tilde{p}_k(m,1/2)$,
which will be given a direct bijective derivation in Section
\ref{sec:oneface}: 

\begin{proposition}
  \label{sec:quasi-bipartite-case}
  For $k\geq 0$ and $m \in \Z_{\geq 0} + \frac12$, $\tilde{p}_k(m)$ is
  the number of planar tight quasi-bipartite maps with one face of odd
  degree $2m$, one face of degree one, and $k+1$ distinct marked
  vertices, one of them distinguished and labeled, say as vertex $1$,
  and the remaining $k$ unlabeled.
\end{proposition}

The polynomials $\tilde{p}_k$ obey the dilaton equation
 \begin{equation}
    \label{eq:tildepdil}
\tilde{p}_k(m_1,m_2;\ldots,m_n,1) - \tilde{p}_k(m_1,m_2;\ldots,m_n,0) = \tilde{p}_{k-1}(m_1,m_2;\ldots,m_n)
  \end{equation}
 and the string equation (for $m_i$'s as in Theorem~\ref{thm:quasibiptheorem})
  \begin{multline}
    \label{eq:tildepstr}
    (k+1) \tilde{p}_{k+1}(m_1,m_2;\ldots,m_n) = (m_1 + \cdots + m_n - k - 1) \tilde{p}_k(m_1,m_2;\ldots,m_n) + \\ 
    \sum_{i=1}^n \sum_{0<j<m}(2j) \tilde{p}_k(m_1,\ldots,m_{i-1},j,m_{i+1},\ldots,m_n)
  \end{multline}
  where we sum over half-integer values of $j$ for $i=1$ and $2$ and
  over integer values of $j$ for $i\geq 3$.  The proof is similar to
  that of Proposition~\ref{prop:pqdilstr} and uses
  now~\eqref{eq:tildepkunivstr}. Again, this corresponds to Norbury's
  dilaton and string equations in the planar quasi-bipartite case.

\begin{remark}
  We have the relation
  $\tilde{p}_k(1/2,1/2;m_1,\ldots,m_n)=q_k(m_1,\ldots,m_n)$ which
  implies that
  $N_{0,n+2}(2m_1,2m_2,2m_3,\ldots,2m_n,1,1)=N_{0,n+2}(2m_1,2m_2,2m_3,\ldots,2m_n,2,0)$
  for $m_1,\ldots,m_n$ integers. The latter equality can be explained
  via a ``slit-slide-sew'' bijection in the spirit
  of~\cite{Bettinelli2020}.
\end{remark}

\subsection{Quasi-polynomials counting planar tight maps with more
  odd faces}
\label{sec:quasipolmain}

In this section, we provide explicit expressions for the lattice count polynomials
$\mathsf{N}^{(k)}_{0,n}$, enumerating planar tight maps with $k$ boundaries
of odd lengths, and $n-k$ boundaries of even lengths, for an arbitrary
value of 
$k\geq 3$. 

To this end, similarly to the bipartite and quasi-bipartite cases
discussed in the preceding section, we first need to introduce a
two-parameter family of univariate polynomials which generalizes those
introduced above: for $k$ a non-negative integer and $e\in \Z$, we
define
\begin{equation}
\label{eq:pkedef1}
  p_{k,e}(m):=\frac{1}{(k!)^2}\, \prod_{i=1}^k\left(m^2-\left(i-\frac{e}2\right)^2\right) = \binom{m+\frac{e}2-1}{k} \binom{m-\frac{e}2+k}{k} .
\end{equation}
We recover the polynomials $p_k$, $\tilde{p}_k$, and $q_k$ of
Section~\ref{sec:main-results} for $e=0,1,2$, respectively. We will
provide combinatorial interpretations of these polynomials in Section
\ref{sec:extens-non-bipart}.  

Next, let $r,s$ be non-negative integers and $\epsilon\in \Z$
be fixed. For $m\in \Z/2$, we let 
 \begin{equation}
    \label{eq:pieps}
    \pi_{r,s}^{(\epsilon)}(m) :=
    \begin{cases}
      \binom{r+s}s p_{r+s,s+1+\epsilon}(m) & \text{if $m-\frac{s+1+\epsilon}2 \in \Z$,} \\
      0 & \text{otherwise.}
    \end{cases}
  \end{equation}
  For every choice of $r,s,\epsilon$, this defines a quasi-polynomial
  in the variable $2m$. For the purposes of stating the main theorem
  of this section (Theorem \ref{thm:quasimain}), only the cases $\epsilon \in \{0,1\}$ will be of
  interest. Note that for $m=0$ we have
\begin{equation}
  \label{eq:piepsm0}
    \pi_{r,s}^{(\epsilon)}(0) =
    \begin{cases}
     \delta_{r,0}\delta_{s,0} & \text{if $\epsilon=1$,} \\
      0 & \text{if $\epsilon=0$.}
    \end{cases}
\end{equation} 

We may now state the main theorem of this section. In \eqref{eq:7} below and
later, for $\epsilon \in
\{0,1\}$, we will write $\bar\epsilon:=1-\epsilon$ to lighten the
notation. 

\begin{theorem}
  \label{thm:quasimain}
  For $n \geq 3$, for $m_1,m_2,\ldots,m_n \in \Z_{\geq 0}/2$, with at
  least three of the $m_i$ being half-integers, the number
  $N_{0,n}(2m_1,2m_2,\ldots,2m_n)$ of planar tight maps with $n$
  boundaries labeled from $1$ to $n$ with respective lengths
  $2m_1,2m_2,\ldots,2m_n$ is given by the symmetric quasi-polynomial
  \begin{multline}
    \label{eq:quasimain}
    N_{0,n}(2m_1,2m_2,\ldots,2m_n)= \\
    \sum_{\left(\substack{\epsilon_1,\ldots,\epsilon_n\\r_1,\ldots,r_n\\s_1,\ldots
       ,s_n}\right)\in I_n}
    \left( \sum_{i=1}^n r_i \right)! \left( \sum_{i=1}^n \epsilon_i s_i \right) \left( \sum_{i=1}^n s_i -1 \right)!
     \prod_{i=1}^n \pi_{r_i,s_i}^{(\epsilon_i)}(m_i),
  \end{multline}
where     $I_n$ is the (finite) subset of
$\{0,1\}^{n} \times \Z_{\geq 0}^n \times \Z_{\geq 0}^n$ defined by
\begin{equation}\label{eq:7}
  I_n := \left\{ \left(
      \begin{array}{c}
        \epsilon_1,\ldots,\epsilon_n\\r_1,\ldots,r_n\\s_1,\ldots
       ,s_n
     \end{array}
   \right):
    \begin{array}{l}
      \sum\limits_{i=1}^n \epsilon_i=\sum\limits_{i=1}^n r_i +1  \\
      \sum\limits_{i=1}^n \bar\epsilon_i=\sum\limits_{i=1}^n s_i+2
    \end{array}
\, , \sum_{i=1}^ns_i\geq 1\right\}\, .
\end{equation}

\end{theorem}

\begin{remark}
  \label{rem:unification}
  We note that the right-hand side of \eqref{eq:quasimain} is equal to
  $0$ when the number $k$ of faces of odd degree is equal to $0$ or
  $2$, so that the formula does {\em not} hold in these cases, which
  have been respectively dealt with above in Theorems
  \ref{thm:maintheorem} and
  \ref{thm:quasibiptheorem}. Theorem~\ref{thm:quasimain} yields a
  non-trivial result only when $k\geq 4$ is an even number, since a
  map necessarily has an even number of faces of odd degrees. As a
  sanity check, it is not difficult to see that the right-hand side of
  \eqref{eq:quasimain} vanishes when $k$ is odd. Indeed, assume
  without loss of generality that $2m_1,\ldots,2m_k$ are odd numbers,
  and that $2m_{k+1},\ldots,2m_n$ are even.  By \eqref{eq:pieps} and
  \eqref{eq:piepsm0}, the product term in the sum \eqref{eq:quasimain}
  is non-zero only if $s_i+\epsilon_i$ is even for $1\leq i\leq k$,
  and odd for $k+1\leq i\leq n$. On the other hand, the constraints in
  the definition of the summation index implies that
  $\sum_{i=1}^n(s_i+\epsilon_i)+2=n$, which after reduction modulo
  $2$, shows that $n-k$ and $n$ have the same parity, so that $k$ is
  necessarily even. A fully unified formula, encompassing
  Theorems~\ref{thm:maintheorem}, \ref{thm:quasibiptheorem}
  and~\ref{thm:quasimain} is given in Theorem~\ref{thm:themaintheorem}
  below.
\end{remark}

The proof of Theorem~\ref{thm:quasimain} will follow an architecture
similar to the bijective proof of Theorems \ref{thm:maintheorem} and
\ref{thm:quasibiptheorem}, building from elementary examples of maps
with explicit enumeration formulas, to construct general ones.  In
particular, note that \eqref{eq:quasimain} reduces to
$\pi_{r,s}^{(1)}(m)$ if we specialize it to $n=r+s$ with $m_1=m$,
$m_2=\cdots=m_{s+3}=1/2$ and $m_{s+4}=\cdots=m_{r+s}=0$.  However, it
is not a priori obvious to obtain an interpretation of
$\pi^{(0)}_{r,s}(m)$ by specializing formula \eqref{eq:quasimain}. For
this reason, we will need to investigate in some depth these
quasi-polynomials and relate them to the combinatorial notion of {\em
  petal trees}. This will be the object of Section
\ref{sec:petaltrees}, but let us record right away the definition of
this notion so as to state one important result,
Proposition~\ref{prop:pieps}, which will be used at the end of Section
\ref{sec:slicings}.

\begin{figure}
\centering
\includegraphics{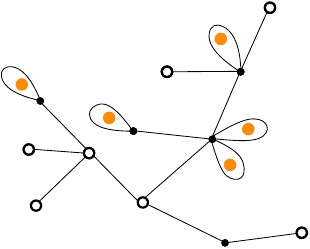}
\caption{A tight petal tree with 5 petals (displayed with yellow marks) and 7 marked 
vertices (shown in white).}
\label{fig:petaltree}
\end{figure}

We call {\em petal} a face of degree $1$. 
A \emph{petal tree} is a planar map having an
\emph{exterior} face of arbitrary degree, and such that every other
face is a petal. A tight
petal tree is just a petal tree with marked vertices, which is tight as
a map. See Figure~\ref{fig:petaltree} for an illustration.

\begin{proposition}
  \label{prop:pieps}
  For $r,s$ nonnegative integers, $m \in \Z_{> 0}/2$, and
  $\epsilon\in \{-1,0,1\}$, the number of tight petal trees with an exterior
  face of degree $2m$, $s+1+\epsilon$ petals, $1+\epsilon$ of which
  distinguished, and $r+1-\epsilon$ marked vertices, $1-\epsilon$ of
  which distinguished, is equal to $\pi_{r,s}^{(\epsilon)}(m)$. 
\end{proposition}

\begin{remark} 
Proposition \ref{prop:pieps}, which will be proved in Section \ref{sec:petaltrees}, holds for $m>0$. It can be extended to $m=0$, provided we restrict the value of $\epsilon$ 
to the set $\{0,1\}$. Indeed, in this case, \eqref{eq:piepsm0}
is consistent with Proposition \ref{prop:pieps} upon understanding the exterior face of degree $0$ as a distinguished marked vertex. 
For $\epsilon=0$ or $1$, the only possible map with such a marked vertex, $s+1+\epsilon$ petals (and no other face) and 
$r+1-\epsilon$ other marked vertices is made of
a single loop connecting the distinguished marked vertex and separating two petals. It has $s+1+\epsilon=2$ and  $r+1-\epsilon=0$,
hence $\epsilon=1$ and $r=s=0$ (note that the distinction of the two petals does not create any degeneracy).
\end{remark}
\begin{remark} 
The situation above is quite similar to the case $m=1/2$, for which we obtain
\begin{equation}
  \label{eq:piepsmonehalf}
    \pi_{r,s}^{(\epsilon)}\left(\frac{1}2\right) =
    \begin{cases}
    0 & \text{if $\epsilon=-1$,} \\
       \delta_{r,0}\delta_{s,0} & \text{if $\epsilon=0$,}\\
    0 & \text{if $\epsilon=1$.}
    \end{cases}
\end{equation}
This agrees with Proposition \ref{prop:pieps} since the only possible map with an exterior face of degree $1$, 
$s+1+\epsilon$ petals and $r+1-\epsilon$ marked vertices is made of a single loop connecting a unique vertex and separating
the exterior face from a unique petal. It has $s+1+\epsilon=1$ and $r+1-\epsilon\leq 1$, hence $\epsilon=0$, $s=0$ and $r=0$ (note that the unique vertex is therefore marked).
\end{remark}

\section{Connection with the enumeration of non necessarily tight maps}
\label{sec:slicings}

\subsection{Equivalence with Tutte's slicings formula in the
  (quasi-)bipartite case}\label{sec:equiv-with-tutt}

One of the earliest results in map enumeration is Tutte's
\emph{slicings formula}~\cite{Tutte1962} which, in our current
terminology, asserts that the number $M(\ell_1,\ldots,\ell_n)$ of
planar (non necessarily tight) \emph{bipartite} maps with $n \geq 3$
labeled faces of prescribed even degrees $2\ell_1,\ldots,2\ell_n$ is
given by
\begin{equation}
  \label{eq:slicings}
  M(\ell_1,\ldots,\ell_n) = (\ell_1 + \cdots + \ell_n-1)_{n-3} \prod_{i=1}^n
  \binom{2\ell_i-1}{\ell_i},
\end{equation}
where $(\ell)_k:=\ell(\ell-1)\cdots(\ell-k+1)$ denotes the falling
factorial. Note that the faces are assumed unrooted and that the
formula extends to the case where some, but not all, of the $\ell_i$
vanish, with the convention $\binom{-1}{0}=1$, upon again understanding that a face of degree $0$ is a marked
vertex (planar maps with three or more labeled faces or vertices have
no symmetries).

In this section, we explain how Tutte's slicings formula is related to
Theorem~\ref{thm:maintheorem} giving a formula for the number of
planar \emph{tight} bipartite maps with $n$ boundaries of prescribed
lengths. As we shall see, the two formulas can be deduced from one
another.

\begin{figure}[t]
  \centering
  \includegraphics[width=.6\textwidth]{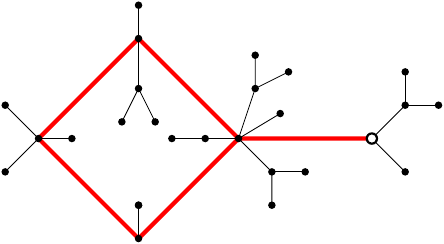}
  \caption{The tight core (thick red edges) of a planar bipartite map
    with two faces and one marked vertex (shown in white). The rest of
    the map consists of rooted plane trees (possibly empty) attached
    in the corners of the tight core.}
  \label{fig:tightcore}
\end{figure}

The key observation, already made in~\cite[Section~4]{Budd2020a}, is
that an arbitrary (non necessarily tight) map, possibly with marked
vertices, can be bijectively decomposed into a tight map (which we
call the \emph{tight core}) and a collection of rooted plane trees
(without marked vertices) attached to the corners of the tight
core. See Figure~\ref{fig:tightcore} for an illustration. Precisely,
the tight core has the same number of faces and marked vertices as the
arbitrary map, and a face of degree $2\ell$ in the arbitrary map
yields a face of degree $2m$ in the tight core, for some
$m \leq \ell$, together with a plane forest made of $2m$ trees having
$\ell-m$ edges in total. The number of such plane forests is equal to
\begin{equation}
  \label{eq:Adef}
  A_{\ell,m} := \frac{2m}{2\ell} \binom{2\ell}{\ell-m},
\end{equation}
see for instance \cite[$\triangleright$ I.38]{FlSe09}, with conventionally $A_{0,0}=1$ and $A_{\ell,m}=0$ for $m>\ell$. As a
consequence, we have
\begin{equation}
  \label{eq:MNrel}
  M(\ell_1,\ldots,\ell_n) = \sum_{m_1=0}^{\ell_1} \cdots \sum_{m_n=0}^{\ell_n}
  A_{\ell_1,m_1} \cdots A_{\ell_n,m_n} N_{0,n}(2m_1,\ldots,2m_n)
\end{equation}
where $N_{0,n}$ is the number of planar tight bipartite maps with $n$
boundaries of lengths $2m_1,\ldots,$ $2m_n$, as defined in
Section~\ref{sec:counting-tight-maps}. Note that the matrix
$(A_{\ell,m})_{\ell,m \geq 0}$ is unitriangular, hence the
formula~\eqref{eq:MNrel} can be inverted as
\begin{equation}
  \label{eq:NMrel}
  N_{0,n}(2m_1,\ldots,2m_n) = \sum_{\ell_1=0}^{m_1} \cdots \sum_{\ell_n=0}^{m_n}
  B_{m_1,\ell_1} \cdots B_{m_n,\ell_n} M(\ell_1,\ldots,\ell_n)
\end{equation}
where $B$ is the inverse of $A$. This inverse is given explicitly by
$B_{m,\ell}=(-1)^{m-\ell} \binom{m+\ell-1}{m-\ell}$ but we will not
use its expression in the following.

Now, let us substitute Tutte's slicings formula~\eqref{eq:slicings}
into~\eqref{eq:NMrel}. By the Chu-Vandermonde identity, we may expand
the falling factorial as
\begin{equation}
  \label{eq:ChuVan}
  (\ell_1 + \cdots + \ell_n-1)_{n-3} = (n-3)!
  \sum_{k_1,k_2,\ldots,k_n\geq 0\atop k_1+k_2+\cdots+k_n=n-3} \binom{\ell_1-1}{k_1} \binom{\ell_2}{k_2} \cdots \binom{\ell_n}{k_n},
\end{equation}
which is nothing but an equality between polynomials in
$\ell_1,\ldots,\ell_n$. This yields
\begin{equation}
  N_{0,n}(2m_1,\ldots,2m_n) = (n-3)!
  \sum_{k_1,k_2,\ldots,k_n\geq 0\atop k_1+k_2+\cdots+k_n=n-3} \hat{p}_{k_1}(m_1) \hat{q}_{k_2}(m_n) \cdots \hat{q}_{k_n}(m_n)
\end{equation}
where
$\hat{p}_k(m):=\sum_{\ell=0}^m B_{m,\ell} \binom{\ell-1}{k}
\binom{2\ell-1}{\ell}$ and
$\hat{q}_k(m):=\sum_{\ell=0}^m B_{m,\ell} \binom{\ell}{k}
\binom{2\ell-1}{\ell}$. We recover Theorem~\ref{thm:maintheorem}, with
the multivariate polynomial $p_{n-3}(m_1,m_2,\ldots,m_n)$ defined
via~\eqref{eq:defpkqkmultivariate}, provided that $\hat{p}_k(m)$ and
$\hat{q}_k(m)$ are respectively equal to the univariate polynomials
$p_k(m)$ and $q_k(m)$ defined in~\eqref{eq:defpkqk}. This is ensured
by the following:
\begin{lemma}
  \label{lem:conv}
  The univariate polynomials $p_k(m)$ and $q_k(m)$ defined
  in~\eqref{eq:defpkqk} satisfy
  \begin{equation}
    \label{eq:conv}
    \sum_{m=0}^\ell A_{\ell,m} p_k(m) = \binom{\ell-1}{k} \binom{2\ell-1}{\ell}, \qquad
    \sum_{m=0}^\ell A_{\ell,m} q_k(m) = \binom{\ell}{k} \binom{2\ell-1}{\ell}.
  \end{equation}
\end{lemma}

\begin{proof}
  These hypergeometric identities can be proved using algorithmic
  methods, see \cite{Petkovsek1996} and references therein.

  Alternatively, a bijective proof for $\ell \geq 1$ follows from
  Propositions~\ref{prop:pkinterpret} and \ref{prop:qkinterpret},
  themselves proved bijectively in Section~\ref{sec:oneface}. More
  precisely, the first identity is obtained by counting in two
  different ways plane trees with $\ell$ edges and $k+2$ distinct
  marked vertices, two of them distinguished and labeled. Namely, the
  left-hand side is obtained via the tight core decomposition and
  Proposition~\ref{prop:pkinterpret}, while the right-hand side is
  obtained by a direct enumeration: $\binom{2\ell-1}{\ell}$ is the
  number of plane trees with $\ell$ edges and two distinguished
  labeled vertices, and $\binom{\ell-1}{k}$ is the number of ways to
  choose the $k$ other marked vertices. The second identity is
  obtained similarly by counting in two different ways plane trees
  with $\ell$ edges, one of them marked, and $k+1$ distinct marked
  vertices, one of them distinguished.
\end{proof}

Note that, doing the above reasoning backwards, it is conversely
possible to recover Tutte's slicings formula from
Theorem~\ref{thm:maintheorem}, using~\eqref{eq:MNrel}. We now briefly
discuss the quasi-bipartite case: let $\ell_1,\ell_2$ be positive
half-integers, and $\ell_3,\ldots,\ell_n$ be non-negative integers,
$n \geq 3$. Then, by~\cite[Section~6]{Tutte1962}, the number
$M(\ell_1,\ldots,\ell_n)$ of planar quasi-bipartite maps with $n$
labeled boundaries of prescribed lengths $2\ell_1,\ldots,2\ell_n$
reads
\begin{equation}
  \label{eq:slicingsquasi}
  M(\ell_1,\ldots,\ell_n) = (\ell_1 + \cdots + \ell_n-1)_{n-3} \binom{2\ell_1-1}{\ell_1-\tfrac12} \binom{2\ell_2-1}{\ell_2-\tfrac12} \prod_{i=3}^n
  \binom{2\ell_i-1}{\ell_i}.
\end{equation}
The tight core decomposition works as before and we find
that~\eqref{eq:MNrel} still holds, upon understanding that the sums
over $m_1$ and $m_2$ should be now taken over half-integer values,
$A_{\ell,m}$ being still defined by~\eqref{eq:Adef} for $\ell,m$
half-integers. By a slight variant of the reasoning above, we may
deduce Theorem~\ref{thm:quasibiptheorem} with the multivariate
polynomial $\tilde{p}_k(m_1,m_2;m_3,\ldots,\allowbreak m_n)$ being given
by~\eqref{eq:defptildekmultivariate}. Namely, we modify the
expansion~\eqref{eq:ChuVan} of the falling factorial by replacing
$\ell_1-1$ and $\ell_2$ in the right-hand side by $\ell_1-\tfrac12$
and $\ell_2-\tfrac12$, respectively, and we make use of the identity
\begin{equation}
  \label{eq:convbis}
  \sum_{m\in\{\frac12,\frac32,\ldots,\ell\}} A_{\ell,m} \tilde{p}_k(m) = \binom{\ell-\tfrac12}{k} \binom{2\ell-1}{\ell-\tfrac12}
\end{equation}
valid for $\ell$ a positive half-integer. Again, this identity can be
proved either via algorithmic methods, or via a bijective argument: it
is obtained by counting in two different ways planar maps with one
face of odd degree $2\ell$, one face of degree one and $k+1$ distinct
marked vertices, one of them distinguished. The left-hand side is
obtained by the tight-core decomposition together with Proposition
\ref{sec:quasi-bipartite-case} (which will be derived bijectively in
Section \ref{sec:oneface}). As for the right-hand side, note that, by
collapsing the face of degree one, such maps correspond to rooted
plane trees with $\ell-\frac12$ edges: $\binom{2\ell-1}{\ell-1/2}$ is
the number of such trees with one distinguished vertex, and
$ \binom{\ell-1/2}{k}$ is the number of ways to choose the $k$ other
marked vertices.

\subsection{A non-bipartite slicings formula}\label{sec:non-bipart-slic}

By arguing similarly, we may use Theorem \ref{thm:quasimain} to obtain
a generalization of Tutte's slicings formula counting planar maps with
a prescribed degree sequence. The relevant identity to use, valid for
$\ell-\frac{s+1+\epsilon}{2}\in \Z$ and $\epsilon \in \{0,1\}$, is
\begin{equation}
  \label{eq:8}
  \sum_{m\in
    \Z_{\geq 0}/2}A_{\ell,m}\pi^{(\epsilon)}_{r,s}(m)=\binom{2\ell-1}{\ell-\frac{s+1+\epsilon}{2}\,
    ,\, 
    \ell-r-\frac{s+1-\epsilon}{2}\, ,\, r\, ,\, s}\, ,
\end{equation}
where in the left-hand side, we observe by \eqref{eq:pieps} that only
the terms for which 
$\ell-m\in \Z_{\geq 0}$ contribute, 
and in the right-hand side, we use a multinomial coefficient notation. 
To understand this formula, recall from Proposition \ref{prop:pieps}
that $\pi^{(\epsilon)}_{r,s}(m)$ counts tight petal trees, i.e.\ tight
maps with an exterior face of degree $2m$, $s+1+\epsilon$ petals,
$1+\epsilon$ of which are distinguished, and $r+1-\epsilon$ marked
vertices, $1-\epsilon$ of which are distinguished. By applying the
tight core decomposition, the left-hand side of \eqref{eq:8} expresses
the number of petal trees which are not necessarily tight, with an
exterior face of degree $2\ell$, and with the same number of
(distinguished) petals and (distinguished) marked
vertices as described in the preceding sentence. Checking that this
number equals the right-hand side of 
\eqref{eq:8} is a straightforward exercise 
based on the methods used in Section \ref{sec:petaltrees} to  prove
Proposition \ref{prop:pieps}, and is simpler due to the absence of the
tightness condition. 

Formula \eqref{eq:MNrel} remains unchanged if $\ell_1,\ldots,\ell_n$
are allowed to take half-integer values, except that the corresponding
sums should then run over half-integer $m_i$'s as well. Substituting
in \eqref{eq:MNrel} the formula of Theorem \ref{thm:quasimain} for
$N_{0,n}(2m_1,\ldots,2m_n)$, we obtain the following:

\begin{theorem}[A census of non-bipartite slicings]
  \label{thm:extendedslicings}
  The number $M(\ell_1,\ldots,\ell_n)$ of planar maps with $n$ labeled
  faces of degrees $2\ell_1,\ldots,2\ell_n$, at least four of which
  are odd, is given by
  \begin{multline}
    \label{eq:extendedslicings}
    M(\ell_1,\ldots,\ell_n)
    =\sum_{\left(\substack{\epsilon_1,\ldots,\epsilon_n\\r_1,\ldots,r_n\\s_1,\ldots
          ,s_n}\right)\in I_n} \left( \sum_{i=1}^n r_i \right)! \left(
      \sum_{i=1}^n \epsilon_i s_i \right) \left( \sum_{i=1}^n s_i -1
    \right)!
    \\
    \times\prod_{i=1}^n
    \binom{2\ell_i-1}{\ell_i-\frac{s_i+1+\epsilon_i}{2},
      \ell_i-r_i-\frac{s_i+1-\epsilon_i}{2},r_i,s_i}\, ,
  \end{multline}
  where $I_n$ is as in \eqref{eq:7}, and where it is understood that
  the multinomial coefficient vanishes whenever
  $\ell_i-\frac{s_i+1+\epsilon_i}{2}$ is not an integer. The formula
  makes sense when some $\ell_i$ vanish, upon understanding that a
  face of degree $0$ is in fact a vertex.
\end{theorem}

\begin{example}
  For $\ell_1,\ldots,\ell_4\in \Z_{\geq 0} + \frac12$, performing the
  sum in~\eqref{eq:extendedslicings} yields a number of planar maps
  with four labeled faces of odd degrees $2\ell_1,\ldots,2\ell_4$
  equal to
  \begin{equation}
    M(\ell_1,\ldots,\ell_4)=
    (\ell_1+\cdots+\ell_4-2)
    \prod_{i=1}^4 \binom{2\ell_i-1}{\ell_i-\frac12}.
  \end{equation}
  In particular, we find
  $M\left(\frac32,\frac12,\frac12,\frac12\right)=2$ which, after
  rooting each face, gives a number of slicings equal to $6$
  consistently with~\cite[Section~6]{Tutte1962}.
\end{example}

To our knowledge, this extension of Tutte's slicings formula for
general non-bipartite planar maps with prescribed degrees is new.
Note that~\eqref{eq:extendedslicings} does not hold when the number of
faces of odd degree is zero or two, see the discussion in
Remark~\ref{rem:unification} in the tight setting. Even though we
obtain Theorem~\ref{thm:extendedslicings} as a consequence of
Theorem~\ref{thm:quasimain}, the former could be proved directly by
the approach of Section~\ref{sec:extens-non-bipart}, forgetting about
the tightness constraint.

\section{Bijective proofs in the bipartite and quasi-bipartite cases}
\label{sec:bijective}
\subsection{Case of maps with one face}
\label{sec:oneface}
In this section, we first present a combinatorial proof of
Propositions~\ref{prop:pkinterpret} and \ref{prop:qkinterpret} by a
direct enumeration of the planar tight bipartite maps with one face
considered in these propositions.  Our approach is inspired from the
bijective interpretation of Narayana numbers given
in~\cite[Section~3.2]{Dershowitz1980}.

More precisely, we wish to enumerate planar bipartite maps with one face of degree $2m$, which are nothing but \emph{plane trees with $m$ edges} ($m\geq 1$), 
endowed with
a distinguished marked vertex labeled $1$, and either a second distinguished marked vertex labeled $2$ in the context of Proposition~\ref{prop:pkinterpret} or with a marked edge (the root edge) in the context of Proposition~\ref{prop:qkinterpret}. These two 
situations will be referred to respectively as ``case (p)''  and ``case (q)'' in the following,
to remind the reader that they concern the combinatorial interpretation of
$p_k(m)$ and $q_k(m)$ respectively. The trees are finally decorated by
the choice of $k$ additional unlabeled marked vertices ($k\geq 0$),
with the constraint that \emph{all the leaves of the tree are marked
  vertices} (either labeled or unlabeled), except possibly the
endpoints of the root edge in the case (q) if it happens that such an endpoint is a leaf.

\begin{figure}[t]
  \centering
  \fig{}{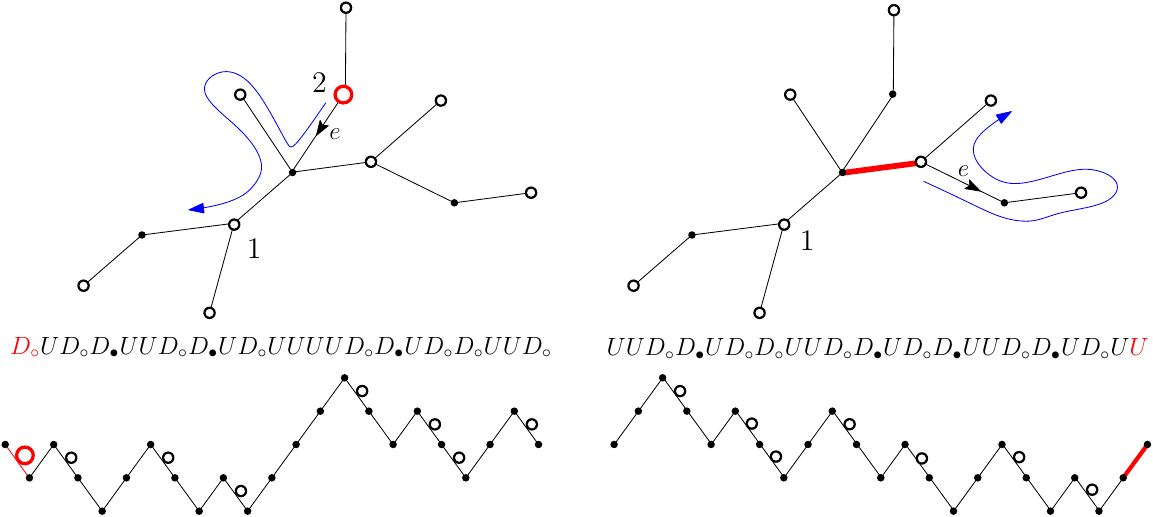}
  \caption{The coding of planar tight maps with one face with degree
    $2m$ (here $m=11$) by dressed words with letters $U$, $D_\circ$
    and $D_\bullet$. Left: in the context of
    Proposition~\ref{prop:pkinterpret}, such a map has $k+2$ (here
    $k=7$) marked vertices (represented as white, as in Figure
    \ref{fig:sampletightmap}), two of them distinguished as vertex $1$
    and $2$ (represented by a bigger red circle). The associated word,
    indicated under the map, starts with a $D_\circ$.  Right: in the
    context of Proposition~\ref{prop:qkinterpret}, the map has $k+1$
    marked vertices, one of them distinguished as vertex $1$, and a
    root edge (represented in thick red). The associated word ends
    with a $U$. We indicated under each word its visualization as a
    lattice path, where the vertex markings have been transferred to
    descending steps (as indicated by circles). We also indicated the
    starting oriented edge $e$ for the contour word of each map.}
  \label{fig:treetoword}
\end{figure}

Our proof is based on the classical coding of plane trees by their
contour word, here in terms of the letters $U$ (up) and $D$
(down). The following discussion is illustrated on
Figure~\ref{fig:treetoword}. Ignoring the $k$ unlabeled vertex
markings for now, the coding that we use here is adapted to trees
\emph{with both a marked vertex $1$ and a marked oriented edge $e$}
whose choice in cases (p) and (q) will be discussed below. For any
such tree, the contour word is obtained as follows: we start from the
right of $e$ and visit all the edge sides counterclockwise around the
tree. We then record a letter $U$ if we move away from vertex $1$ on
the tree, and a letter $D$ otherwise. After one turn around the tree,
we get a word of $2m$ letters containing $m$ occurrences of $U$ and
$m$ occurrences of $D$, since the two sides of any given edge give
rise to exactly one $U$ and one $D$. Viewing the successive $U$'s and
$D$'s as successive up and down steps, the coding may alternatively be
represented as a lattice path of length $2m$ starting and ending at
the same height (a so-called bridge in the lattice path terminology)
and whose nodes correspond to the successive visited corners around
the tree. In this representation, the corners at vertex $1$ are
associated with the nodes with minimal height. Setting the minimal
height to $0$, the height of a node is nothing but the graph distance
to vertex $1$ of the vertex incident to the associated corner. The
above coding by words/paths is clearly bijective.

To use this coding in the context of Propositions~\ref{prop:pkinterpret} and \ref{prop:qkinterpret} where the trees already have a distinguished vertex $1$, 
we need a canonical prescription for the choice of the oriented edge $e$ 
at which we start the contour. In case (p) where the tree is endowed
with a second distinguished vertex $2$, we take for $e$ the first edge
of the branch between 
vertex $2$ and vertex $1$, oriented towards $1$. Clearly the knowledge of $e$ and that of the vertex $2$ are equivalent, but we note that, by construction, the associated word necessarily \emph{starts with a $D$} in case (p).

In case (q), we first orient the root edge away from vertex $1$ and
take for $e$ the edge following it along the counterclockwise contour
around the tree. The edge $e$ is therefore incident to the endpoint of
the root edge further away from vertex $1$ and we orient it away from
that vertex\footnote{Note that it may happen that the edge $e$ be
  identical to the root edge itself, but with the opposite
  orientation, in which case the first letter of the word is a $D$. In
  all other cases the first letter is a $U$.}. Clearly the knowledge
of $e$ and that of the root edge are equivalent, but we note that, by
construction, the last visited edge side in the contour is that of the
root edge itself, going away from vertex $1$, hence the associated
coding word necessarily \emph{ends with a $U$} in case (q).

It remains to introduce the $k$ additional vertex markings. The markings may be 
recorded in the coding word as follows: every vertex $v$ distinct from vertex $1$ 
may be associated bijectively with a letter $D$ of the coding
word. Indeed $v$ is bijectively associated with the edge $e(v)$
incident to $v$ that belongs to the branch between $v$ and the vertex
$1$ and exactly one of the two sides of $e(v)$ is coded by the letter
$D$. If $v$ is a marked vertex, we transfer its marking to the
associated letter $D$, which we denote by $D_\circ$ to record the
marking. If $v$ is not a marked vertex, the associated letter $D$ will
be denoted by $D_\bullet$, so that the letter $D$ eventually appears
in two flavors $D_\circ$ and $D_\bullet$, leading to dressed words
made of the three letters $U$, $D_\circ$ and $D_\bullet$. In case (p),
we also transfer  the marking of vertex $2$, so that the first letter (which we know is originally a $D$) is now a $D_\circ$. The numbers of $U$, $D_\circ$ and $D_\bullet$ letters are therefore, respectively, $m$, $k+1$ and $m-k-1$ in case (p) and 
$m$, $k$ and $m-k$ in case (q). 

Apart from possibly vertex $2$ (which is marked anyway) in case (p) or possibly an endpoint of the root edge (which needs not being marked) in case (q), any leaf in the tree corresponds to a sequence $UD$ in the associated word.
Requiring that all leaves be marked boils down to  demanding that any $D$ following
a $U$ be marked, i.e.~the sequence $UD_\bullet$ is not allowed in the dressed words.

Altogether, a word coding for a tree in case (p) has the canonical form
  \begin{equation}
    \label{eq:UDform1}
    D_\circ D_\bullet^{a_1} U^{b_1} D_\circ D_\bullet^{a_2} U^{b_2} \cdots
    D_\circ D_\bullet^{a_{k+1}} U^{b_{k+1}}
  \end{equation}
  where the $a_i$ and $b_i$ are nonnegative integers such that
  $a_1 + a_2 + \cdots + a_{k+1} = m-k-1$ and
  $b_1 + b_2 + \cdots + b_{k+1} = m$. In other words, the $a_i$ and
  $b_i$ form weak compositions of $m-k-1$ and $m$, respectively, into
  $k+1$ summands. There are respectively $\binom{m-1}{k}$ and
  $\binom{m+k}{k}$ such compositions, hence the number of trees in case (p)
  is $\binom{m-1}{k}\binom{m+k}{k}=p_k(m)$ as wanted.
  
 Similarly, a word
  coding for a tree in case (q) has the canonical form 
  \begin{equation}
    \label{eq:UDform2}
     D_\bullet^{a_1} U^{b_1} D_\circ D_\bullet^{a_2} U^{b_2} D_\circ \cdots
    D_\bullet^{a_k} U^{b_k} D_\circ D_\bullet^{a_{k+1}} U^{b_{k+1}}U
  \end{equation}
  where the $a_i$ and $b_i$ form weak compositions of $m-k$ and
  $m-1$, respectively, into $k+1$ summands\footnote{Note that the case $a_1=b_1=0$ corresponds to a word starting with a $D_\circ$ which, cyclically, comes after
  the last letter $U$. This situation corresponds to the case where the endpoint
  of the root edge further away from vertex $1$ is a leaf and is marked. The case $a_1>0$ corresponds to the case where this vertex is a leaf and is unmarked and, finally, 
  the case $a_1=0$, $b_1>0$ corresponds to the case where this vertex is not a leaf.}. There are respectively $\binom{m}{k}$ and
  $\binom{m+k-1}{k}$ such compositions, hence the number of trees in case (q)
  is $\binom{m}{k}\binom{m+k-1}{k}=q_k(m)$ as wanted.

This ends the combinatorial proof of Propositions~\ref{prop:pkinterpret}
and \ref{prop:qkinterpret}.

\begin{figure}[t]
  \centering
  \fig{.5}{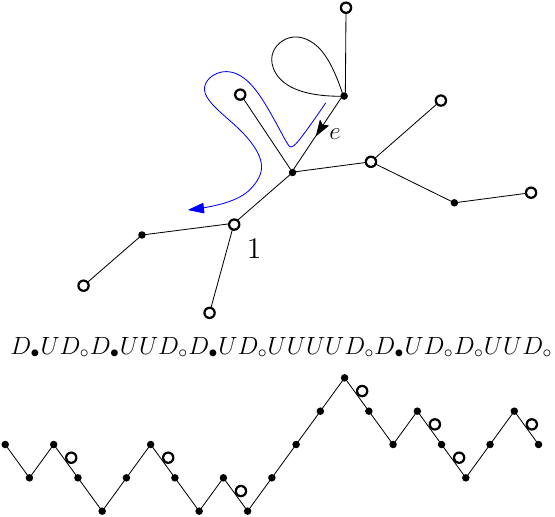}
  \caption{The coding of planar tight maps with one face with odd degree $2m$ (here $m=23/2$) and one face of degree $1$ by dressed words with letters $U$, $D_\circ$ and $D_\bullet$. As in Proposition~\ref{sec:quasi-bipartite-case} the map has $k+1$ (here $k=7$) marked vertices, one of them 
 distinguished as vertex $1$. We indicated under the map the associated word, of length $2m-1$, 
 and its visualization as a lattice path. We also indicated on the map the starting oriented edge $e$ for the contour word.}
  \label{fig:treetowordquasi}
\end{figure}
\medskip
\paragraph{Quasi-bipartite case. }

The proof of Proposition~\ref{sec:quasi-bipartite-case} is obtained
along similar lines, see Figure~\ref{fig:treetowordquasi} for an
example. Indeed, we may transform bijectively a planar map with one
face of odd degree $2m$ ($m\in \Z_{\geq 0}+3/2$), one face of degree
$1$ and a distinguished vertex $1$ into a plane tree with $m-1/2$
edges with both a marked vertex $1$ and a marked oriented edge
$e$. This is done by considering the unique vertex incident to the
loop formed by the degree $1$ face, by marking its incident edge $e$
lying immediately to the left of that loop, with $e$ oriented away
from the vertex and finally erasing the loop. We can now use our
coding of such pointed rooted trees by words with $2m-1$
letters. Taking the markings into account gives rise to a dressed word
with exactly $k$ occurrences of $D_\circ$, $m-1/2-k$ occurrences of
$D_\bullet$ and $m-1/2$ occurrences of $U$, with no occurrence of the
sequence $UD_\bullet$, hence with canonical form
  \begin{equation}
    \label{eq:UDform3}
    D_\bullet^{a_1} U^{b_1} D_\circ D_\bullet^{a_2} U^{b_2} \cdots
    D_\circ D_\bullet^{a_{k+1}} U^{b_{k+1}}
  \end{equation}
  where the $a_i$ and $b_i$ are nonnegative integers such that
  $a_1 + a_2 + \cdots + a_{k+1} = m-k-1/2$ and
  $b_1 + b_2 + \cdots + b_{k+1} = m-1/2$. In other words, the $a_i$ and
  $b_i$ form weak compositions of $m-k-1/2$ and $m-1/2$, respectively, into
  $k+1$ summands. There are respectively $\binom{m-1/2}{k}$ and
  $\binom{m+k-1/2}{k}$ such compositions, hence the number of maps at hand is $\binom{m-1/2}{k}\binom{m+k-1/2}{k}=\tilde{p}_k(m)$ as announced. This holds for $m\geq 3/2$. For $m=1/2$, the value $\tilde{p}_k(1/2)=\delta_{k,0}$ is consistent with the fact that there is a unique planar map with two (distinguished) faces of degree $1$ and one marked vertex labeled $1$, which is its unique vertex so that the map cannot host any other marked vertices when $k>0$.

\subsection{Case of maps with two faces}
\label{sec:twofaces}

Let us now provide a bijective proof of Proposition
\ref{prop:pk2interpret}. To this end, we will first need to
reinterpret slightly the objects counted by $p_k(m),q_k(m)$ that were
discussed in the preceding section.

\begin{definition}
  \label{sec:bipartite-case.-}
  For given integers
$a\geq 1$ and $b\geq 0$, an {\em $(a,b)$-forest} is a tight map with exactly two faces
$f,f_*$, such that: 
\begin{itemize}
\item $f_*$ is a simple face\footnote{We say that a face is
    \emph{simple} if its contour is a simple cycle, i.e. does not
    visit a same vertex several times.} of degree $a+b$ and one
  distinguished incident vertex $v_*$,
\item the $a$ vertices following and including $v_*$ in
  counterclockwise order around $f_*$ are not marked. 
\end{itemize}
\end{definition}

\begin{figure}[t]
  \centering
  \includegraphics[width=.8\textwidth]{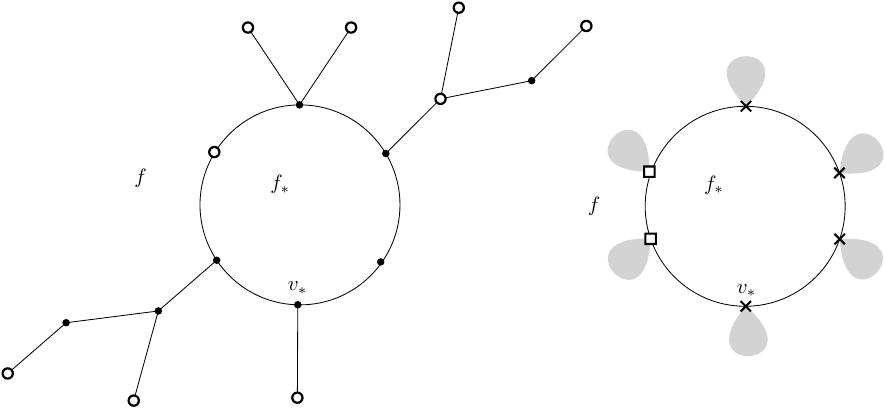}
  \caption{Left: an example of a $(4,2)$-forest, with $9$ marked
    vertices and size $28$. The marked
    vertices are represented as white, as in Figure \ref{fig:sampletightmap}. Note that 
    none of the first $4$ vertices arriving in counterclockwise order
    after $v_*$ are marked, while the following $2$ comprise one
    marked and one unmarked vertex. Right: a schematic, generic
    representation of a $(4,2)$-forest, where the grey blobs 
    represent tree components (that may be reduced to a single root
    vertex)
    and white squares represent markable roots, while crosses
    represent unmarkable roots.}
  \label{fig:sampleforest}
\end{figure}

We call the $a$ vertices referred to above as the {\em unmarkable}
vertices, and the other $b$ vertices are called the {\em markable}
vertices. In the illustrating figures, starting with
Figure~\ref{fig:sampleforest}, the latter will be represented by white
squares, while the former will be represented by crosses.
Equivalently, by removing the $a+b$ edges incident to $f_*$, we may
view an $(a,b)$-forest as a linearly ordered collection of $a+b$
rooted plane trees starting from the one rooted at $v_*$, whose leaves
(non-root vertices of degree $1$) are all marked, and such that the
roots of the first $a$ trees are unmarked while those of the remaining
$b$ trees may be marked or not.

The {\em size} of an $(a,b)$-forest is the degree of the face $f$. If we
view it as a collection of trees as above, then this size is equal to
$2e+a+b$ where $e$ is the total number of edges in the trees composing
the forest. 

\begin{figure}[t]
  \centering
  \includegraphics[width=.9\textwidth]{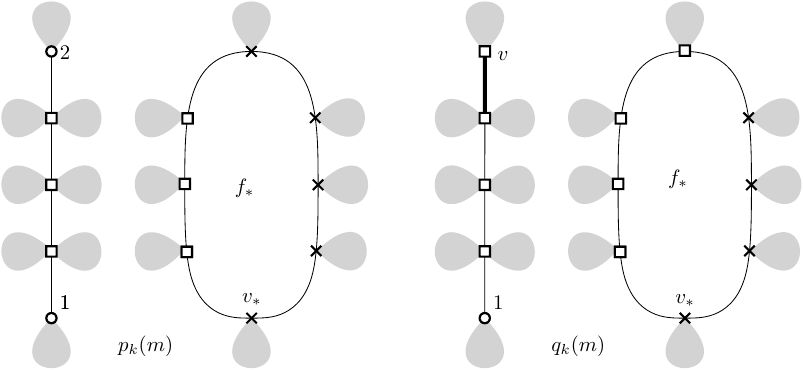}
  \caption{The cutting operation turning tight maps with one face into
    forests. On the left, tight maps counted by $p_k(m)$ are cut into
    $(d+1,d-1)$ forests, and on the right, rooted
    tight maps counted by $q_k(m)$ become $(d,d)$-forests. }
  \label{fig:unfoldtree}
\end{figure}

We now describe two simple bijections, illustrated on
Figure~\ref{fig:unfoldtree}, linking the numbers $p_k(m)$ and $q_k(m)$
to the forests discussed above. First, recall from Proposition
\ref{prop:pkinterpret} the interpretation of $p_k(m)$ as counting
tight bipartite planar maps with one face $f$ of degree $2m$ and $k+2$
distinct marked vertices, of which exactly two are labeled as $1$ and
$2$. There is a natural operation consisting in cutting open the
branch $\gamma$ linking the distinguished vertices, into a simple face
$f_*$ of degree $2d$, where $d\geq 1$ is the graph distance between
these vertices. In doing so, we duplicate the vertices lying on the
path $\gamma$, except its extremities, into ``left and right'' copies
(upon orienting $\gamma$ from vertex $1$ to vertex $2$), and in case
some of these vertices are marked, we always decide to transfer the
mark to the left copy. The vertex initially distinguished and labeled
as $1$ is then renamed as $v_*$ and seen as unmarked, while we remove
the mark and label on the vertex initially labeled $2$. The result is
then a $(d+1,d-1)$-forest. Conversely, given a $(d+1,d-1)$-forest for
some $d\geq 1$, we can glue together the $r$-th edge of $f_*$ in
counterclockwise order starting from $v_*$ with the opposite
$(2d-r+1)$-th one, for $r\in \{1,2,\ldots,d\}$, relabel $v_*$ as
vertex $1$ and the diametrally opposite vertex of $f_*$, lying at
distance $d$ from $v_*$, as vertex $2$, and finally, after gluing the
vertices in pairs along the contour of $f_*$, transferring to the
newly created vertices the marks carried by all markable vertices. By
construction, every markable vertex is matched to an unmarkable
vertex, and this operation is the inverse of the cutting procedure
described above. Finally, these operations preserve the number $k$ of
marked (unlabeled) vertices, and are size-preserving in the sense that
the degree $2m$ of the unique face of the map to be cut corresponds to
the size of the resulting forest.

Similarly, recall from Proposition \ref{prop:qkinterpret} that
$q_k(m)$ enumerates rooted tight bipartite planar maps with one face
$f$ of degree $2m$ and $k+1$ distinct marked vertices, exactly one of
them being distinguished and labeled as vertex $1$. Let $v$ be the
vertex incident to the root edge of such a map, and which is further
away from vertex $1$. Again, we cut open along the branch from vertex
$1$ to $v$, with length $d\geq 1$ say, creating a simple face $f_*$ of
degree $2d$. We transfer the marks along this path to the left copies
of the vertices created in the cutting operation, and if $v$ happens
to be marked, we keep this mark. Finally, we remove the mark on vertex
$1$ and rename it as $v_*$. This results in a $(d,d)$-forest, since
now the vertex diametrically opposite to $v_*$ in $f_*$ is a markable
vertex. This operation is clearly invertible by a similar gluing
operation as above, and it preserves the number of marked unlabeled
vertices as well as the degree of $f$. We may conclude with the
following statement.

\begin{proposition}
  \label{sec:case-maps-with}
For integers $m\geq 1$ and $k\geq 0$, the number $p_k(m)$ (resp.\
$q_k(m)$) is the cardinality of the set 
of $(d+1,d-1)$-forests (resp.\ $(d,d)$-forests) with size $2m$ and $k$
marked unlabeled vertices, where $d$ can take any value in $\Z_{>0}$.
\end{proposition}

Our bijective proof of Proposition \ref{prop:pk2interpret} will
consist in showing that a tight map with two faces $f_1,f_2$ of
respective degrees $2m_1,2m_2$, and with $k$ unlabeled marked
vertices and one extra distinguished marked vertex with label $1$ can
be decomposed uniquely and bijectively into a pair of tight maps consisting in
\begin{itemize}
\item a $(d_1+1,d_1-1)$-forest with $k_1$ marked vertices
 \item a $(d_2,d_2)$-forest with $k_2$ marked vertices
\end{itemize}
where $d_1,d_2\geq 1$ and $k_1+k_2=k$. By Proposition
\ref{sec:case-maps-with} and the definition
\eqref{eq:defpkqkmultivariate} of $p_k(m_1,m_2)$, this immediately
implies Proposition  \ref{prop:pk2interpret}.

Given an $(a,b)$-forest and an integer $c$ such that $1\leq c\leq \min(a-1,
b)$,
the {\em $c$-partial gluing} of the forest is the map obtained by
gluing the $r$-th edge following $v_*$ in counterclockwise order
around $f_*$ with the opposite $(a+b+1-r)$-th one, for $r\in \{1,2,\ldots,c\}$, 
and transferring any mark on the markable vertices to the resulting glued
vertices. The vertex inherited from $v_*$ in the new map is distinguished and labeled
as vertex $1$, while the last vertex to be glued, lying at distance
$c$ from $v_*$, is distinguished and called $v_{**}$. Since $c\leq \min(a-1,
b)$, some edges remain unglued, and the resulting
map still has two faces which we call $f,f_{**}$, where $f_{**}$ is
the ``remnant'' of $f_*$, of degree $a+b-2c$, and we do not change the name for the
exterior face since it has the same contour information as the
original face. Moreover, the assumption that $c\leq \min(a-1 ,b)$ implies
that every markable vertex is glued to an unmarkable vertex.

\begin{figure}[t]
  \centering
  \includegraphics[scale=.9]{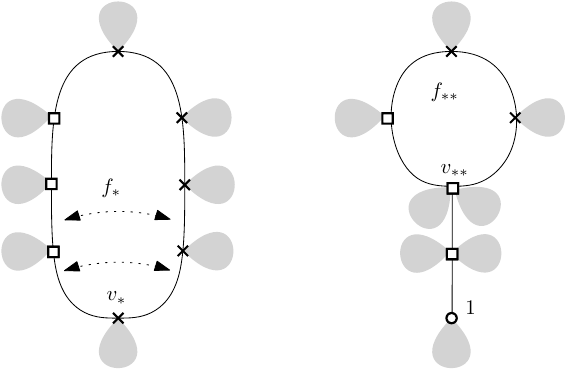}
  \caption{The $2$-partial gluing of a $(5,3)$-forest yielding a $(2,2)^*$-forest. }
  \label{fig:partialgluing}
\end{figure}

The obtained map is then an $(a',b')^*$-forest
with $a'=a-c-1$ and $b'=b-c+1$, according to the following definition,
similar to Definition \ref{sec:bipartite-case.-}: 

\begin{definition}\label{sec:bipartite-case.--1}
  For integers $a\geq 0$ and $b\geq 1$, an {\em $(a,b)^*$-forest}
  is a tight map with exactly two faces
$f,f_{**}$, such that: 
\begin{itemize}
\item $f_{**}$ is a simple face of degree $a+b$,
 \item there is an extra distinguished vertex labeled $1$ incident to the 
   face $f$, and we let $v_{**}$ be the vertex incident to $f_{**}$
   that is closest to $1$, 
\item the $a$ vertices following and excluding $v_{**}$ in
  counterclockwise order around $f_{**}$ are not marked. 
\end{itemize}
\end{definition}
The partial gluing operation is clearly invertible, by cutting along
the simple path of length $c$ from $v_{**}$ to the distinguished
labeled vertex $1$.  See Figure \ref{fig:partialgluing} for an
illustration.

\begin{figure}[htbp]
  \centering
  \includegraphics[scale=1]{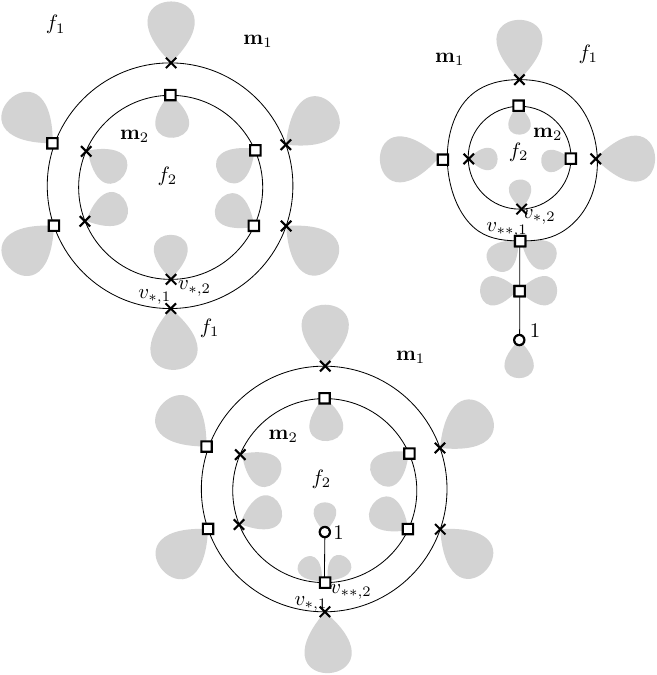}
  \caption{Decomposition of a tight map with two faces $f_1,f_2$, $k$
    marked unlabeled vertices and one extra distinguished labeled
    vertex, into two forests.  The top left, top right, and bottom
    pictures represent respectively the situations where the
    distinguished labeled vertex belongs to the common boundary of
    $f_1$ and $f_2$, is incident to $f_1$ but not $f_2$, and is
    incident to $f_2$ but not $f_1$. We denote by $2d$ the length of
    the cycle separating $f_1$ and $f_2$. In the first case, the map
    is decomposed into a $(d+1,d-1)$ forest glued to a $(d,d)$
    forest. In the second case, the map is decomposed into the
    $(d_1-d)$-partial gluing of a $(d_1+1,d_1-1)$-forest with a
    $(d,d)$-forest, and in the third case, it is instead a
    $(d+1,d-1)$-forest glued to the $(d_2-d)$-partial gluing of a
    $(d_2,d_2)$-forest. Note that, on this picture, the
    counterclockwise order along the face $f_{*,2}$, which is the
    ``exterior'' simple face of $\mathbf{m}_2$, appears to be
    clockwise, since $f_{*,2}$ is the unbounded face in the plane
    embedding. }
  \label{fig:matchingtwofaces}
\end{figure}

\begin{proof}[Proof of Proposition \ref{prop:pk2interpret}. ]
  Let $m_1,m_2,d_1,d_2$ be positive integers, and $k_1,k_2$ be
  non-negative integers.  
Suppose we are given a $(d_1+1,d_1-1)$ forest $\mathbf{f}_1$ with size $2m_1$ and
$k_1$ marked vertices, and a $(d_2,d_2)$ forest $\mathbf{f}_2$ with size $2m_2$ and
$k_2$ marked vertices. We let $v_{*,1},v_{*,2}$ be the distinguished
vertices in these maps. As explained above, we wish to use these
pieces to build a planar tight map with two faces of degrees $2m_1,2m_2$, with $k=k_1+k_2$ marked
unlabeled vertices, and with one extra distinguished vertex. There are three
possible situations, illustrated in Figure
\ref{fig:matchingtwofaces}. 

Suppose first that $d_1=d_2=d$. Then we can glue together the two simple
boundaries of $\mathbf{f}_1$ and $\mathbf{f}_2$, in such a way that
$v_{*,1}$ and $v_{*,2}$ are glued 
together into a single distinguished vertex labeled $1$. The next $d$
unmarkable vertices after $v_{*,1}$ (resp. the last $d-1$ markable
vertices) of $\mathbf{f}_1$ are then glued to the
last $d$ markable vertices (resp. the $d-1$
unmarkable vertices following $v_{*,2}$) of $\mathbf{f}_2$. The result is a tight map
with two faces of degrees $2m_1,2m_2$, and with $k_1+k_2$ marked
unlabeled vertices as well as a distinguished marked vertex labeled
$1$ lying on the boundary of both faces. 

Suppose next that $d_1>d_2$. In this case, we first perform the
$(d_1-d_2)$-partial gluing of $\mathbf{f}_1$, resulting in a
$(d_2,d_2)^*$-forest (with distinguished vertex called $v_{**,1}$),
which we glue along the simple face of $\mathbf{f}_2$ by identifying
$v_{**,1}$ and $v_{*,2}$. Note that each of the $d_2$ markable
vertices on either side of the gluing is matched with an 
unmarkable vertex on the other side. The resulting map has a distinguished
labeled vertex incident to $f_1$ but not to $f_2$. 

Finally, the case $d_2>d_1$ is similar, except that we now perform the
$(d_2-d_1)$-partial gluing of $\mathbf{f}_2$ first, resulting in a
$(d_1-1,d_1+1)^*$-forest, whose $d_1+1$ markable vertices and $d_1-1$
unmarkable vertices are matched with the $d_1+1$ unmarkable vertices
and $d_1-1$ markable vertices of $\mathbf{f}_1$.  The resulting map
has a distinguished labeled vertex incident to $f_2$ but not to $f_1$.

The above construction can clearly be inverted by the following
cutting operation. Start from a tight
map $\mathbf{m}$ with two faces $f_1,f_2$, $k$ marked vertices and one extra
distinguished vertex labeled $1$. We observe that such a map is
unicyclic, and therefore contains a unique simple cycle $\gamma$, of
length $2d\geq  2$ say. 
We cut along this cycle, separating $f_1$ and $f_2$. Formally, this
means that we associate with $\mathbf{m}$ the two maps
$\mathbf{m}_1,\mathbf{m}_2$ 
respectively obtained by  
  removing all edges and vertices that are incident to $f_2$ but not
  $f_1$ on the one hand, and $f_1$ but not $f_2$ on the other
  hand. Note that for $i\in \{1,2\}$, $\mathbf{m}_i$ is made of the
  initial face $f_i$ and has an ``exterior'' simple face $f_{*,i}$
  which is the remnant of the face $f_{3-i}$. All marked vertices of
  $\mathbf{m}$ that are not on $\gamma$ are naturally transferred to
  either $\mathbf{m}_1$ or $\mathbf{m}_2$, and we need a convention to
  transfer the marked vertices lying on $\gamma$. 

To this end,   
we distinguish the cases as above depending on whether the vertex
labeled $1$ is incident to both $f_1,f_2$, to 
$f_1$ but not $f_2$, or to $f_2$ but not $f_1$. In the first case, the
cutting operation splits the labeled vertex $1$ into two copies
$v_{*,1},v_{*,2}$, respectively belonging to $\mathbf{m}_1$ and
$\mathbf{m}_2$,  that we declare unmarkable. The $d$ vertices
following $v_{*,1}$ in counterclockwise order around $f_{*,1}$
are declared unmarkable, as well as the $d-1$
vertices following $v_{*,2}$ in counterclockwise order around
$f_{*,2}$, and all other vertices incident to $f_{*,1}$ and $f_{*,2}$
are declared markable. In this way, every vertex of $\gamma$ has been split into a
a markable/unmarkable pair in $\mathbf{m}_1$ and $\mathbf{m}_2$. 
We then transfer the marks that were
located on the vertices $\gamma$ to the unique associated markable duplicate. 
This gives the wanted pair
$(\mathbf{f}_1,\mathbf{f}_2)=(\mathbf{m}_1,\mathbf{m}_2)$ of
$(d+1,d-1)$- and $(d,d)$-forests, 
of sizes $2m_1$ and $2m_2$, 
which receive $k_1$ and $k_2$ marked unlabeled vertices with
$k_1+k_2=k$.

In the second case, we let $v_{**}$ be the vertex incident to $f_2$
that is closest to the distinguished vertex $1$. When cutting along
the cycle $\gamma$, this vertex is separated
into two copies, one called $v_{**,1}$ is incident to $f_{*,1}$ and is
declared markable, as well as the $d-1$ vertices {\em preceding it}
around $f_{*,1}$, the other called $v_{*,2}$ is incident to $f_{2,*}$ and is 
declared unmarkable, as well as the $d-1$ vertices following it around
$f_{*,2}$. The map $\mathbf{m}_2=:\mathbf{f}_2$ is then a $(d,d)$-forest, while we
further cut $\mathbf{m}_1$ along the
simple path of length $d'\geq 1$ from $v_{**,1}$ to the distinguished
vertex $1$, hence creating a
$(d+d'+1,d+d'-1)$-forest $\mathbf{f}_1$, attributing the marked vertices in the
natural way (this operation is the reverse of the
$d'$-partial gluing of the resulting forest). 

The situation in the third case is similar, with a slightly different
convention for the markable and unmarkable vertices, as illustrated in
Figure \ref{fig:matchingtwofaces}.
\end{proof}

\paragraph{Quasi-bipartite case. }

We now consider the quasi-bipartite case where one assumes that 
$m_1,m_2\in \Z_{\geq 0}+1/2$, that is, $2m_1$ and $2m_2$ are odd 
integers, and aim at proving Proposition
\ref{prop:tildepk2interpret}. A discussion parallel to
the above applies, except that the separating cycle between the two
faces of a tight map with faces of degrees $2m_1$ and $2m_2$ will have
an odd length, say $2d-1$ for some $d\geq 1$.

In this situation, unfolding the above argument {\em mutatis
  mutandis}, there is now a canonical
decomposition of a planar tight map with two faces of degrees $2m_1,2m_2$,
with $k$ marked unlabeled vertices and one extra distinguished vertex
labeled $1$ into a pair formed of a $(d_1,d_1-1)$-forest and a
$(d_2,d_2-1)$-forest, for some $d_1,d_2\geq 1$, respectively with sizes $2m_1$ and
$2m_2$ and with $k_1$ and $k_2$ marked
vertices, where $k_1+k_2=k$. The situation is therefore more symmetric
since the glued forests are of the same nature and have the same
numbers of markable vertices. 

By performing
the $(d_i-1)$-partial gluing of the $(d_i,d_i-1)$-forest of size $2m_i$ with $k_i$
marked vertices, we see that such objects are
in bijection with tight maps with one face of degree $2m_i$, one face of
degree $1$, and $k_i$ marked vertices, which are precisely counted by
$\tilde{p}_{k_i}(m_i)$, as discussed in Section \ref{sec:oneface}. Together with the above, this shows that
$\tilde{p}_k(m_1,m_2)$ indeed enumerates the wanted quasi-bipartite
planar tight maps with two faces, as stated in Proposition
\ref{prop:tildepk2interpret}.

\subsection{Case of pointed rooted maps via slices}
\label{sec:pointedrooted}

We now aim at proving Proposition~\ref{prop:pointedrooted},
interpreting $(n-1)!\,  q_{n-1}(m_1,m_2,\ldots,m_n)$ as the number of
pointed rooted planar tight bipartite maps with $n$ labeled boundaries
of respective lengths $2m_1,\ldots,2m_n$, where $m_1,\ldots,m_n$ are
integers not all equal to zero. To this end, we will need the slice
decomposition developed in \cite{hankel,irredmaps}. Here, we follow closely the
presentation of \cite[Section 2.2]{BouttierHDR} and adapt it to the tight setting.

A {\em slice} is a planar map with one distinguished {\em exterior
  face}, whose contour carries three distinguished (but not
necessarily distinct) corners $A$, $B$
and $C$ appearing in this counterclockwise order around the map, that split the
contour in three parts: 
\begin{itemize}
\item the contour segment $AB$, called the {\em left
    boundary}\footnote{Note that, in the accepted denominations ``left
    boundary'' and ``right boundary'', the term ``boundary'' has a
    meaning different from that in the rest of the paper.}, which is a
  geodesic path,
\item the contour segment $AC$, called the {\em right boundary},
  which is the unique geodesic path between its two endpoints, and
  intersects the left boundary only at $A$,
  \item the contour segment $BC$, called the {\em base}. 
\end{itemize}
The length of the base (i.e.\ the number of edges on the corresponding
contour segment, counted with multiplicity) is called the {\em width}
of the slice. The length of the left boundary is called the {\em
  depth}, and the depth minus the length of the right boundary is
called the {\em tilt}. The corner $A$ is called the \emph{apex}.

A slice of width $1$ is called {\em elementary}. The tilt of an elementary
slice is necessarily in $\{1,0,-1\}$. By the uniqueness property of the
right boundary, there is a unique elementary slice of tilt $-1$, called the {\em
trivial} slice, which consists of a single edge with extremities $A=B$
and $C$. The trivial slice differs from the {\em empty} slice
consisting in a single edge with extremities $B$ and $A=C$, which has
tilt $+1$. If we restrict our attention to bipartite maps, as is the
case in this section, there are no slices with tilt $0$. 

\begin{figure}[t]
  \centering
  \includegraphics[scale=1]{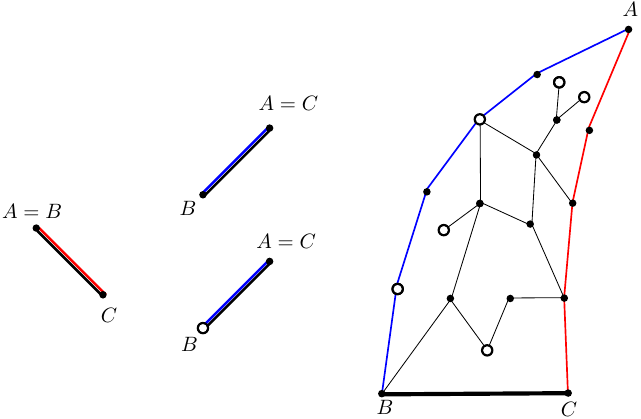}
  \caption{The different types of tight elementary slices, with marked
    vertices shown in white. 
    Left: the trivial slice. Center: the empty slice, in its unmarked and
    marked versions. Right: a non-empty tight
    elementary slice}
  \label{fig:sampletightslice}
\end{figure}

Finally, a {\em tight slice} is a slice, elementary or not, that carries some marked
vertices, in such a way that 
\begin{itemize}
\item 
all vertices of degree $1$ distinct from those
incident to $A,B$ and $C$ are marked,
\item  {\em the right boundary carries no marked
  vertices}.   
\end{itemize}
Note that the vertex incident to $B$ may possibly be marked, but not
those incident to $A$ and $C$, even if those vertices have degree
one. In particular, the empty slice comes with two tight versions,
depending on whether the vertex incident to $B$ is marked or not, and
we will call the marked version the {\em marked empty slice}, which
will play an important role later on. 
See Figure \ref{fig:sampletightslice} for an illustration of the
different types of tight elementary slices. 

\paragraph{Pointed rooted maps and elementary slices. }
There is a simple one-to-one correspondence between pointed rooted planar bipartite maps on the one hand\footnote{Here, as for
  tight maps, we use the slightly unusual convention that a rooted
  map is a map with a distinguished, non oriented edge.}, and non-empty,
bipartite elementary slices of tilt $1$ on the other hand. Starting
from a pointed rooted bipartite map $\mathbf{m}$ with root $e$ and
distinguished vertex $v$, we can perform the opening operation
$O(\mathbf{m},e)$ that opens the edge $e$ into an exterior face of degree
$2$. We let $b$ and $c$ be the two corners incident to this new face, where $c$ is
closest from $v$. We then cut open the map along the leftmost
geodesic\footnote{See for instance \cite[Figure 2.1]{BouttierHDR} for
  a careful definition of the leftmost geodesic from a corner to a vertex.} $\gamma$
from $c$ to $v$, hence enlarging the exterior face. The resulting map
is an elementary bipartite slice of tilt $1$, if we let $B$ be the
corner inherited from $b$, $A$ be the unique corner of the exterior
face incident to $v$, and $C$ be the corner immediately following $B$
as we walk with the exterior face on the right (so that $C$ is one of the two
duplicates of $c$ created after cutting). The fact that the right
boundary $AC$ is the unique geodesic between these two extremities
comes from the fact that $\gamma$ was chosen to be
leftmost, and the slice is non-empty because it has at least one
inner face, inherited from the map we started with.
Conversely, starting from a non-empty elementary bipartite slice of tilt $1$, we
may glue ``isometrically'' together the left and right boundaries starting from the
apex. This results in a bipartite map pointed at the vertex $v$ incident to the
apex, and with a face of degree $2$ whose contour is made of the base
and of the first edge of the left boundary incident to $B$, which are
necessarily distinct since the slice is non-empty. 
We may finally
glue these two edges together into a single edge $e$, at which we root
the resulting map.

These two operations are inverse of one another.  They specialize to a
correspondence between pointed rooted planar {\em tight} bipartite
maps and non-empty {\em tight} bipartite slices of tilt $1$, if we
take the convention that the marks of marked vertices in $\mathbf{m}$
that belong to the leftmost geodesic $\gamma$ considered above should
systematically be transferred to the left boundary of the slice.

\paragraph{Decomposing a slice into a path decorated with elementary
  slices. }
Next, we discuss the decomposition of a bipartite slice $\mathbf{s}$
of width $w\in \Z_{>0}$ and tilt $t\in \Z$ into a
collection of elementary slices. We list the corners of the base as
$c_0=B,c_1,\ldots,c_w=C$, walking with the exterior face on the right. We let
$\ell_i=d(c_0,A)-d(c_i,A)$ for $i\in \{0,1,\ldots,w\}$, where $d$ is
the graph distance in $\mathbf{s}$, and the distance between two
corners is defined as the distance between their incident
vertices. In particular, $\ell_w=t$ is the tilt of $\mathbf{s}$, and
$L=(\ell_0,\ell_1,\ldots,\ell_w)$ is a walk on $\Z$ with increments
$\ell_i-\ell_{i-1}\in \{-1,+1\}$, that we will systematically
identify with the lattice path made of the union of segments
$[(i-1,\ell_{i-1}),(i,\ell_i)],1\leq i\leq n$
in the plane. 

For every $i\in
\{0,1,\ldots,w\}$, we let $\gamma_i$ be the leftmost geodesic from
$c_i$ to the apex $A$. In particular, $\gamma_0$ and $\gamma_w$ are
respectively the left and right boundary of $\mathbf{s}$. 
For $i\in \{1,2,\ldots,w\}$, we let $v'_i$
be the first vertex common 
to $\gamma_{i-1}$ and $\gamma_i$. Then the map $\mathbf{s}'_i$ delimited by these two
geodesics is an elementary bipartite slice with base $c_{i-1}c_i$,
with apex $A'_i$ incident 
to $v'_i$, and with tilt $t_i=\ell_i-\ell_{i-1}$. If $t_i=-1$, in
which case we say that $i$ is a {\em down step}, then
$\mathbf{s}'_i$ is trivial, while if $t_i=+1$, in which case we call $i$
an {\em up step}, then $\mathbf{s}'_i$ is
non-trivial. It may however be the empty slice, precisely when the
geodesic $\gamma_{i-1}$ starts by following the base edge from
$c_{i-1}$ to $c_i$. 

With a bipartite slice $\mathbf{s}$ with width $w$ and tilt $t$, we have associated
a lattice path $L$ from $(0,0)$ to
$(w,t)$ with increments $\pm 1$, where each of the $(w+t)/2$ up steps
$i$ is decorated with a bipartite
elementary slice $\mathbf{s}'_i$ of tilt $1$, while all the $(w-t)/2$ down steps
are decorated with the trivial slice, so that these last decorations are
in fact irrelevant and can be omitted. 

We can invert this
decomposition: given a lattice
path from $(0,0)$ to $(w,t)$ whose up steps $i$ are decorated with
bipartite elementary slices $\mathbf{s}'_i$ of tilt $1$ (and where
$\mathbf{s}'_i$ is the trivial slice if $i$ is a down step), we may
associate a slice of width $w$ 
in the following way.
For every
down step $i$, we identify the segment
$s_i=[(i-1,\ell_{i-1}),(i,\ell_i)]$ of the lattice path with the
associated trivial slice $\mathbf{s}'_i$, hence color it in red as in Figure
\ref{fig:sampletightslice}. 
Next, for every up step $i$, we consider an 
embedding of $\mathbf{s}'_i$ in the plane in which the base edge is the segment $s_i$,
and so that the left boundary (resp.\ the right boundary) is
represented as a curve, monotone in its two coordinates, that starts from
$(i-1,\ell_{i-1})$ (resp.\ $(i,\ell_i)$), is entirely contained in $[i-1,i]\times
[\ell_{i-1},\infty)$ (resp.\ $[i-1,i]\times [\ell_i,\infty)$), and such that its $r$-th vertex, starting from
the base, has its ordinate equal to $\ell_{i-1}+r-1$ (resp.\ $\ell_i+r-1$). 
By convention, the edges of the left boundaries of the slices
$\mathbf{s}'_i$ are declared blue, while the edges of the right
boundaries are declared red. Note that if $\mathbf{s}'_i$ is an empty
slice, then the above operation simply consists in coloring the
segment $s_i$ in blue. 
Then, every red element (either an edge lying on the
right boundary of some slice, or the segment associated with a down
step of the lattice path), lying in some square $[i-1,i]\times
[\ell-1,\ell]$, attempts to be matched to the first available blue edge
in some square $[j-1,j]\times [\ell-1,\ell]$ for some $j>i$, and all
matched edges are glued together. After this gluing is performed, we
obtain a bipartite slice of width $w$ and tilt $t$, where 
the unmatched edges, i.e.\ the blue
edges which are not preceded by red edges at the same ordinate, and
the red edges that are not followed by a blue edge at the same
ordinate, form respectively the left and right boundaries.

Let us now discuss how this decomposition behaves with respect to the
tightness constraint. We first observe that
it
associates with a {\em  tight} bipartite slice $\mathbf{s}$ of width
$w\geq 1$ and tilt $t$ a lattice path $L$ from $(0,0)$ to $(w,t)$ with
$\pm1$ steps decorated with {\em tight} elementary bipartite slices. The
only ambiguity that should be lifted is how we transfer the marks of
marked vertices that belong to the union of leftmost geodesics
$\gamma_i$ defined above to exactly one of their duplicates. We choose
the duplicate that belongs to the left boundary of the slice
$\mathbf{s}'_i$, where $i$ is the maximal index such that the marked
vertex at hand belongs to $\gamma_{i-1}$. Note that such a maximal
index $i$ always exists, since, by definition, tight slices carry no
marked vertices on their right boundaries, and that the duplicate of
the vertex is different from the apex of $\mathbf{s}'_i$ by maximality
of $i$.  With these conventions, all the slices
$\mathbf{s}'_i,1\leq i\leq w$, with transferred marks, are tight
slices.

Moreover, there is an additional restriction on the family 
$L,(\mathbf{s}'_i,1\leq i\leq w)$ that guarantees that the original
slice be tight, i.e.\ that it contains no undesired unmarked vertices
of degree $1$.
Observe that, in the above correspondence, the vertex incident
to a corner $c_i,i\in \{1,\ldots,w-1\}$ of the base distinct from the extremities 
will have degree $1$ precisely in the situation where $i$ is a down
step and $i+1$ is an up step decorated with the empty slice. Indeed,
if the vertex $v_i$ incident to $c_i$ has degree $1$, then the vertices $v_{i-1}$ and $v_{i+1}$ incident to the
corners $c_{i-1}$ and $c_{i+1}$ are the same vertex, and therefore
the geodesics $\gamma_{i}$ and $\gamma_{i+1}$ delimit the empty slice,
since $\gamma_i$ meets $\gamma_{i+1}$ after one single step. 
Conversely, in the gluing procedure, if a down step $i$ is
immediately followed by an up step $i+1$, then 
the ``red'' segment
$[(i-1,\ell_{i-1}),(i,\ell_{i-1}-1)]$ associated with the 
down-step $i$ will be matched to the first edge of the left boundary
of the slice $\mathbf{s}'_{i+1}$. This will result in a vertex of degree
$1$ precisely when this first edge is equal to the base edge, and the
only elementary slice of tilt $1$ with this property is the empty
slice. Consequently, in a tight slice, every such up step $i+1$ 
is decorated with a slice $\mathbf{s}'_i$ that is either non-empty, or
is the {\em marked empty slice}, that is the empty slice with marked
base vertex $B$. 

By forgetting the redundant information of up steps that are
decorated with unmarked empty slices, that is, by letting $(\mathbf{s}_j,1\leq j\leq k)$ be the
sequence $(\mathbf{s}'_i,1\leq i\leq w)$ whose trivial and unmarked empty elements have
been removed,  we obtain the following
result. 

\begin{proposition}
  \label{sec:case-pointed-rooted}
  There is a one-to-one correspondence between tight bipartite slices of width
  $w$ and tilt $t$ on the one hand, and
  pairs of the form $(L,(\mathbf{s}_j,1\leq j\leq k))$ on the other
  hand, where: 
   \begin{itemize}
   \item $L$ is a 
   lattice path from $(0,0)$ to
   $(w,t)$ with $\pm1$ steps, that has $k$ marked up
   steps, in such a way that every up step immediately
   following a down step is marked, 
 \item for $1\leq j\leq k$, $\mathbf{s}_j$ is either the marked empty
   slice or a tight bipartite elementary slice with at least one inner
   face.
\end{itemize}
\end{proposition}

Now observe that if $\mathbf{s}$ is a non-empty bipartite slice with
tilt $1$, then the base edge is incident to an inner face. Calling
$2m>0$ the degree of this face, and after removing the base edge, we
obtain a bipartite slice with width $2m-1$ and tilt $1$. This simple
operation preserves the tight characters of the maps at hand, which
implies the following: 
\begin{corollary}
  \label{sec:case-pointed-rooted-1}
   There is a one-to-one correspondence between non-trivial, non-empty tight bipartite
   elementary slices, whose inner face incident to
   the base edge has degree $2m>0$ on the one hand, 
   and pairs of the form $(L,(\mathbf{s}_j,1\leq j\leq k))$ on the
   other hand, where $k$ is some non-negative integer and: 
   \begin{itemize}
   \item $L$ is a 
   lattice path from $(0,0)$ to
   $(2m-1,1)$ with $\pm1$ steps, that has $k$ marked up
   steps, in such a way that every up step immediately
   following a down step is marked, 
 \item for $1\leq j\leq k$, $\mathbf{s}_j$  is either the marked empty
   slice or a tight bipartite elementary slice with at least one inner
   face.
   \end{itemize}
\end{corollary}

Finally, by convention, we extend the above correspondence by
associating with the marked empty slice the pair
$(\{(0,0)\},\varnothing)$ consisting of the trivial lattice path of
length zero, with
no marks. 

\begin{figure}[t]
  \centering
  \includegraphics[width=\textwidth]{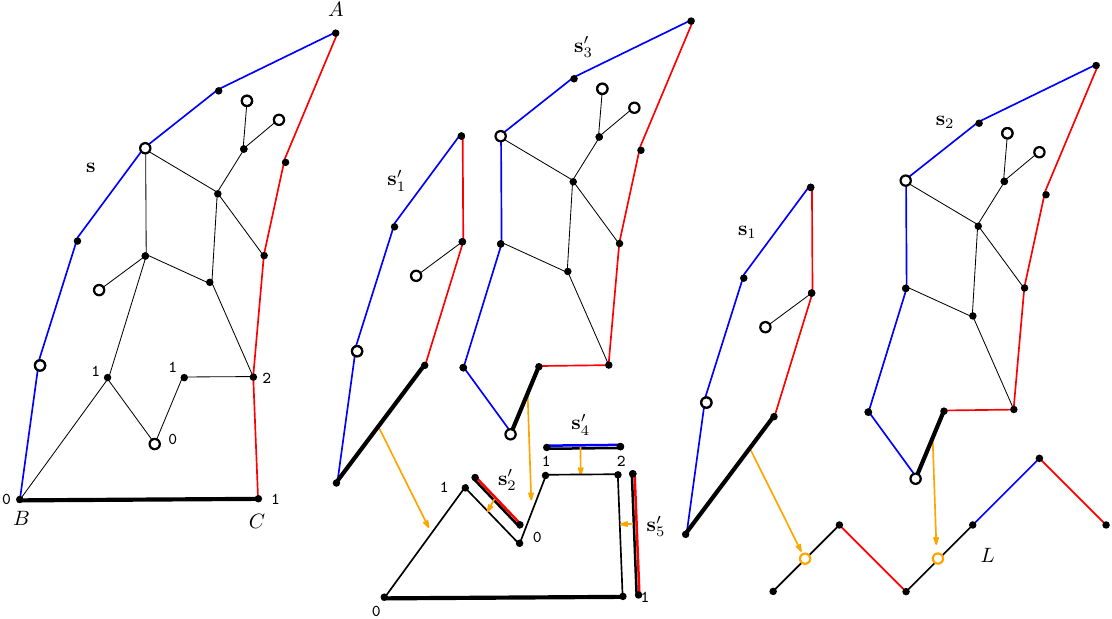}
  \caption{Illustration of the decomposition of a non-empty elementary tight
    slice of tilt $1$ into a lattice path decorated with elementary
    slices. The middle picture shows the
    decomposition along the leftmost geodesics started from the
    corners of the inner face incident to the base. The right-hand
    side shows the result of the decomposition. Note that the trivial
    and empty
    slices $\mathbf{s}'_2,\mathbf{s}'_4$ and $\mathbf{s}'_5$ may be
    discarded, at the price of transferring their color to the
    corresponding step of the lattice path, while the two other slices
  are non-empty and therefore correspond to marked up steps. Note also
  that the second up step, which is consecutive to the first down step,
  is marked, as is required.  
}
  \label{fig:cactus-6}
\end{figure}

\begin{figure}[t]
    \centering
    \includegraphics[scale=.8]{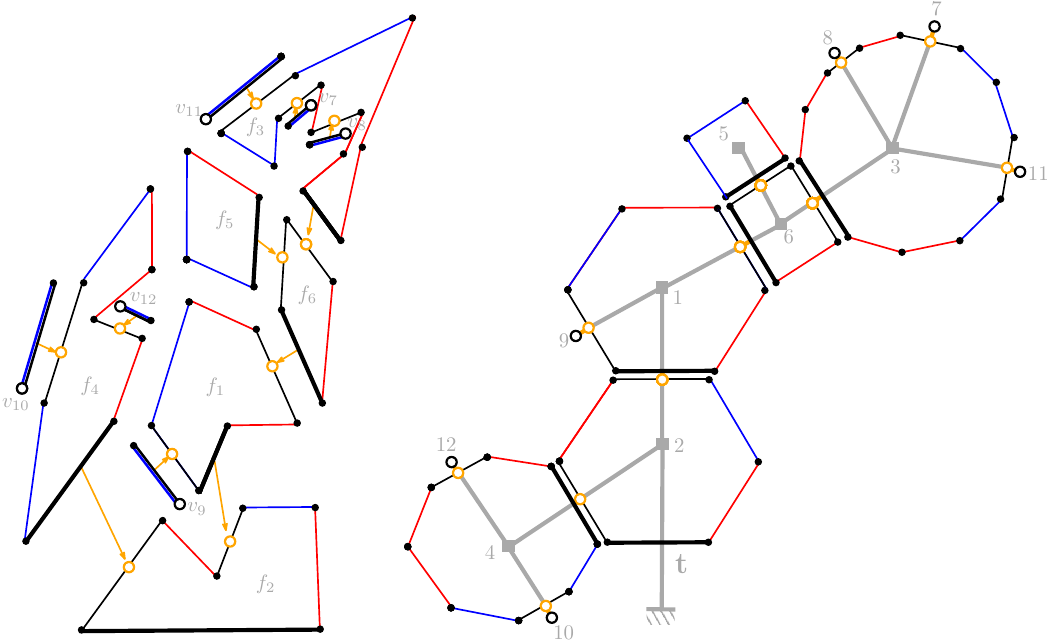}
    \caption{Iterating the slice decomposition yields a plane tree
      whose vertices are decorated by marked paths. On this picture,
      we do not represent the lattice paths, but encode them using a
      color code on polygons drawn around the vertices of
      $\mathbf{t}$, where black/blue correspond to 
      marked/unmarked up steps and red to 
      down steps. Note that some of the terminal nodes of $\mathbf{t}$
      correspond to non-trivial paths with no marked steps, while some
      others correspond to the trivial path, represented as white
      dots, and are associated in turn to
      the marked vertices of the original tight slice. }
    \label{fig:cactus-9}
  \end{figure}

By iterating the decomposition of this corollary, i.e.\ inductively replacing each
non-empty elementary slice in the above decomposition
by a lattice path with some marked up steps, and an ordered family of 
as many elementary slices, 
we obtain a plane tree (where the plane order is
induced by the order of the up steps to which the slices are
connected).  See Figure \ref{fig:cactus-9} for an illustration. 
For integers $m_1,\ldots,m_n$ not all
equal to $0$, let $\mathcal{T}(m_1,\ldots,m_n)$ be the
family of pairs $(\mathbf{t},(L_i)_{1\leq i\leq n})$ where: 
\begin{itemize}
\item 
$\mathbf{t}$ is a rooted plane tree with vertices labeled by
$\{1,2,\ldots,n\}$, 
and, denoting by $k_i$ the number of children of the vertex labeled
$i$ in $\mathbf{t}$, 
\item if $m_i>0$ then 
  $L_i$ is a lattice path from $(0,0)$ to $(2m_i-1,1)$ with $k_i$
  marked up steps, such that all up steps immediately following a down
  step are marked, 
  \item if $m_i=0$ then $L_i$ is the trivial path $\{(0,0)\}$, in
    which case necessarily $k_i=0$. 
  \end{itemize}

\begin{corollary}
  \label{sec:decomp-slice-into}
 For any non-negative integers $m_1,\ldots,m_n$ not all
 equal to $0$, the iterated slice decomposition yields a one-to-one correspondence between pointed rooted
  planar tight maps with $n$ labeled boundaries of respective lengths
  $2m_1,\ldots,2m_n$, and the set $\mathcal{T}(m_1,\ldots,m_n)$. 
\end{corollary}

We are now in position to prove Proposition
\ref{prop:pointedrooted}. 
Note that the number of lattice paths of
length $2m-1$ from $0$ to $1$ that has exactly $k$ marked up steps, including all up steps
immediately following down steps, is precisely the number
$q_k(m)$. Indeed, they are exactly
counted by the words discussed in Section \ref{sec:oneface}, with
$m-1$ letters $U$, $m-k$ letters $D_\bullet$ and $k$ letters
$D_\circ$, with forbidden subword $UD_\bullet$ (up
to flipping upside down the lattice paths to better match the
interpretation of the letters $U,D$). In fact, this even
holds for $m=k=0$ since in this case $q_k(0)=\delta_{k0}$, so we
interpret $q_0(0)$ combinatorially as counting the unique
marked empty slice.   

We now proceed to enumerating the elements of the set
$\mathcal{T}(m_1,\ldots,m_n)$. By~\cite[Section 5,
Equation~(18)]{Bernardi2014}, for a given $n$-uple $(k_1,\ldots,k_n)$
of non-negative integers, if $k_1+\cdots+k_n=n-1$ then there are
exactly $(n-1)!$ rooted labeled plane trees on the vertex set
$\{1,2,\ldots,n\}$ such that vertex $i$ has $k_i$ children for all
$i=1,\ldots,n$, and there are no such tree otherwise. For any such
tree $\mathbf{t}$, there are $\prod_{i=1}^nq_{k_i}(m_i)$ possible
choices of marked lattice paths $L_i,1\leq i\leq n$ such that
$(\mathbf{t},(L_i)_{1\leq i\leq n})$ belongs to
$\mathcal{T}(m_1,\ldots,m_n)$. This finally explains the wanted
formula
\begin{equation}
  \label{eq:8bis}
  \begin{split}
    N_{0,n+2}(2m_1,\ldots,2m_n,2,0)&=(n-1)!\sum_{\substack{k_1,\ldots,k_n\geq 0\\
      k_1+\cdots+k_n=n-1}}\prod_{i=1}^nq_{k_i}(m_i)\\
                                 &=(n-1)!\, q_{n-1}(m_1,\ldots,m_n)\, ,
                               \end{split}
                             \end{equation}
                               which concludes the proof. 

\subsection{General case}
\label{sec:generalbip}

We will now prove Theorem~\ref{thm:maintheorem} in all generality, as
well as its quasi-bipartite analog,
Theorem~\ref{thm:quasibiptheorem}. As the case of maps with one face
was already treated in Section~\ref{sec:oneface}, it suffices to treat
the case where the first and second boundaries are faces, that is when
$m_1,m_2>0$. To this end, we will combine the ideas of
Sections~\ref{sec:twofaces} and~\ref{sec:pointedrooted} via the
following:

\begin{proposition}
  \label{prop:generalbij}
  Let $m_1,m_2$ be positive integers or half-integers and let
  $m_3,\ldots,m_n$ be non-negative integers ($n \geq 3$). Then, there
  is a bijection between the set of planar tight maps with $n$
  boundaries labeled from $1$ to $n$ with respective lengths
  $2m_1,2m_2,\ldots,2m_n$, and the set of pairs
  $(\mathbf{m}_{12},\mathbf{s})$ such that there exists
  $k \in \{0,\ldots,n-3\}$ for which:
  \begin{itemize}
  \item $\mathbf{m}_{12}$ is a planar tight two-face map, with two
    faces of respective degrees $2m_1$ and $2m_2$, and $k+1$ distinct
    marked vertices, one of them distinguished,
  \item $\mathbf{s}=(\mathbf{s}_1,\ldots,\mathbf{s}_{k+1})$ is a
    $(k+1)$-tuple of slices such that:
    \begin{itemize}
    \item for each $j=1,\ldots,k+1$, $\mathbf{s}_j$ is either the
      marked empty slice or a tight bipartite elementary slice with at
      least one inner face, whose inner faces and marked vertices are
      labeled by integers in $\{3,\ldots,n\}$,
    \item each $i \in \{3,\ldots,n\}$ appears in exactly one
      $\mathbf{s}_j$ and labels an inner face of degree $2m_i$ for
      $m_i>0$, or a marked vertex for $m_i=0$,
  \item the label $3$ appears in the first slice $\mathbf{s}_1$.
    \end{itemize}
  \end{itemize}
\end{proposition}

Before proving this proposition, let us see how it implies
Theorems~\ref{thm:maintheorem} and \ref{thm:quasibiptheorem}. We start
with the former: our purpose is to enumerate the pairs
$(\mathbf{m}_{12},\mathbf{s})$ of the proposition when $m_1$ and $m_2$
are integers. For a fixed $k \in \{0,\ldots,n-3\}$, the number of
possible maps $\mathbf{m}_{12}$ is equal to $p_k(m_1,m_2)$ by
Proposition~\ref{prop:pk2interpret}, which we proved in
Section~\ref{sec:twofaces}. As for the number of possible
$\mathbf{s}$, it is given by a slight variant of the reasoning in
Section~\ref{sec:pointedrooted}. Indeed, by recursively decomposing
each slice $\mathbf{s}_j$ into a tree of lattice paths, we see that
the set of possible $(k+1)$-tuples $\mathbf{s}$ is in bijection with
the set $\mathcal{F}_{k+1}(m_3,\ldots,m_n)$ defined as the set of
pairs $(\mathbf{f},(L_i)_{3\leq i\leq n})$ where:
\begin{itemize}
\item $\mathbf{f}$ is a plane forest with $k+1$ connected
  components, i.e.\ a $(k+1)$-tuple of rooted plane trees, whose
  vertices are labeled by $\{3,\ldots,n\}$, the label $3$ appearing in
  the first component,
  \item denoting by $k_i$ the number of children of
  the vertex labeled $i$ in $\mathbf{f}$: 
  \begin{itemize}
\item if $m_i>0$ then $L_i$ is a lattice path from $(0,0)$ to
  $(2m_i-1,1)$ with $k_i$ marked up steps, such that all up steps
  immediately following a down step are marked,
\item if $m_i=0$ then $L_i$ is the trivial path $\{(0,0)\}$, in which
  case necessarily $k_i=0$.
\end{itemize}
\end{itemize}
Note that $\mathcal{F}_1(m_3,\ldots,m_n)$ is nothing but the set
$\mathcal{T}(m_3,\ldots,m_n)$ as defined in
Section~\ref{sec:pointedrooted}. By
Proposition~\ref{prop:onetypeforest} of
Appendix~\ref{sec:enumeration-two-type} below, for a given $n$-uple
$(k_3,\ldots,k_n)$, if $k_3+\cdots+k_n=n-3-k$ then there are exactly
$(n-3)!$ plane forests on the vertex set $\{3,\ldots,n\}$ with $k+1$
components, the first of which contains the label $3$, and such that
vertex $i$ has $k_i$ children for all $i=3,\ldots,n$, and there are no
such forests otherwise. For any such forest $\mathbf{f}$, there are
$\prod_{i=3}^nq_{k_i}(m_i)$ possible choices of marked lattice paths
$L_i,3\leq i\leq n$ such that $(\mathbf{f},(L_i)_{3\leq i\leq n})$
belongs to $\mathcal{F}_{k+1}(m_3,\ldots,m_n)$. This gives
\begin{equation}
  \begin{split}
    \mathrm{Card}\left(\mathcal{F}_{k+1}(m_3,\ldots,m_n)\right)
    &=(n-3)!\sum_{
\substack{k_3,\ldots,k_n\geq 0 \\ k_3+\cdots+k_n=n-3-k}}\prod_{i=3}^nq_{k_i}(m_i)\\
    &=(n-3)!\, q_{n-3-k}(m_3,\ldots,m_n).
  \end{split}
\end{equation}
Thus, multiplying by $p_k(m_1,m_2)$, summing over $k$ and using
\eqref{eq:defpkqkmultivariate}, we get
\begin{equation}
  \begin{split}
    N_{0,n}(2m_1,2m_2,2m_3\ldots,2m_n)&=(n-3)! \sum_{k=0}^{n-3} p_k(m_1,m_2) q_{n-3-k}(m_3,\ldots,m_n) \\
    &= (n-3)!\, p_{n-3}(m_1,m_2,m_3,\ldots,m_n)
  \end{split}
\end{equation}
as wanted. This concludes the proof of Theorem~\ref{thm:maintheorem}
assuming Proposition~\ref{prop:generalbij}.

If we now assume that $m_1,m_2$ are half-integers, the only change we
have to do in the above reasoning is that, for a fixed $k$, the number
of possible maps $\mathbf{m}_{12}$ is now equal to
$\tilde{p}_k(m_1,m_2)$ by Proposition~\ref{prop:tildepk2interpret}
(also proved in Section~\ref{sec:twofaces}). Thus, by \eqref{eq:defpkqkmultivariate} and \eqref{eq:defptildekmultivariate}, we now have
\begin{equation}
  \begin{split}
    N_{0,n}(2m_1,2m_2,2m_3\ldots,2m_n)&=(n-3)! \sum_{k=0}^{n-3} \tilde{p}_k(m_1,m_2) q_{n-3-k}(m_3,\ldots,m_n) \\
    &= (n-3)!\, \tilde{p}_{n-3}(m_1,m_2;m_3,\ldots,m_n)
  \end{split}
\end{equation}
which establishes Theorem~\ref{thm:quasibiptheorem}.

The remainder of this section is devoted to the proof of
Proposition~\ref{prop:generalbij}. It uses the slice decomposition of
annular maps which was introduced in \cite[Section~9.3]{irredmaps},
see also \cite[Section 2.2]{BouttierHDR} for another exposition. Here
we will give yet another, modernized, exposition which,
following~\cite{triskell}, makes use of the key notion of
\emph{Busemann function}.

\paragraph{Decomposing an annular map into a pair of paths decorated
  with elementary slices.}

\begin{figure}
  \centering \includegraphics[width=\textwidth]{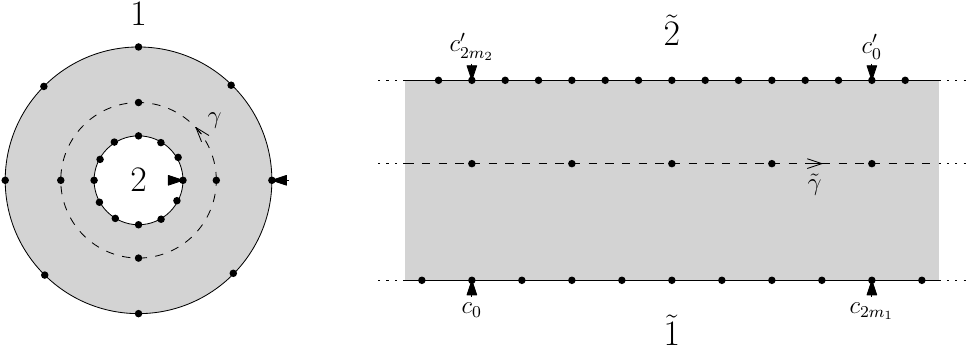}
  \caption{Sketch of an annular map $\mathbf{m}$ (left) and of its
    universal cover $\tilde{\mathbf{m}}$ (right). The external face
    $1$ and the central face $2$ of $\mathbf{m}$ lift respectively to
    the lower face $\tilde{1}$ and to the upper face $\tilde{2}$ of
    $\tilde{\mathbf{m}}$, which have infinite degrees (marked corners
    are represented by arrows). All the other faces form the grey
    region. The innermost minimal separating cycle $\gamma$ lifts to a
    biinfinite geodesic $\tilde{\gamma}$.}
  \label{fig:annulift}
\end{figure}

Let $\mathbf{m}$ be an \emph{annular} map, that is a planar map with
two distinguished faces labeled $1$ and $2$. For now, we do not assume
that $\mathbf{m}$ is a tight map, but we assume that every face other
than $1$ and $2$ has even degree. Denoting by $2m_1$ and $2m_2$ the
respective degrees of $1$ and $2$, $\mathbf{m}$ is either bipartite
when $m_1$ and $m_2$ are integers, or quasi-bipartite when they are
half-integers.

Let us choose a representation of $\mathbf{m}$ in the complex plane
such that face $1$ is the unbounded face, and such that the origin
(point of affix zero) lies in face $2$. We then consider the preimage
$\tilde{\mathbf{m}}$ of $\mathbf{m}$ by the mapping
$z \mapsto \exp(2i \pi z)$: it is an infinite map which we call the
\emph{universal cover} of $\mathbf{m}$ (upon viewing $0$ and $\infty$
as ``punctures''). We refer to~\cite{triskell} for a detailed
discussion of the properties of the universal cover of a map drawn on
the triply-punctured sphere, it can be adapted without difficulty to
the simpler case considered here of a map drawn on the
doubly-punctured sphere. For our purposes, we simply note that the
faces $1$ and $2$ lift in $\tilde{\mathbf{m}}$ to unique faces
$\tilde{1}$ and $\tilde{2}$, which have infinite degrees, while all
the other faces and vertices of $\mathbf{m}$ lift to infinitely many
preimages in $\tilde{\mathbf{m}}$ with the same finite degree. In
particular, $\tilde{\mathbf{m}}$ is bipartite.  See
Figure~\ref{fig:annulift} for an illustration.

Let $s$ denote the \emph{separating girth} of $\mathbf{m}$, that is
the minimal length (number of edges) of a closed path in $\mathbf{m}$
winding around the origin. We call such path a \emph{separating cycle}
and say that it is \emph{minimal} if it has length $s$. Let then
$\gamma$ be the \emph{innermost} minimal separating cycle, which we
define as follows. Consider the interiors of all minimal separating cycles.
It is straightforward to check that their intersection is a simply
connected region containing the origin, and that its boundary
oriented clockwise is still a minimal separating cycle: this is the
innermost minimal separating $\gamma$ that we are looking for. We consider
the path $\tilde{\gamma}=(\tilde{\gamma}(t))_{t \in \Z}$ obtained by
following $\gamma$ counterclockwise infinitely many times and lifting
this biinfinite path in $\tilde{\mathbf{m}}$. The parametrization of
$\tilde{\gamma}$ depends on a choice of a vertex $\tilde{\gamma}(0)$
whose projection belongs to $\gamma$, but we will see that the outcome
of our construction does not depend on it.
By~\cite[Proposition~2.5]{CdVEr10}, $\tilde{\gamma}$ is a
\emph{biinfinite geodesic} in $\tilde{\mathbf{m}}$, which by definition means that
\begin{equation}
  \tilde{d}(\tilde{\gamma}(t),\tilde{\gamma}(t')) = |t-t'|
\end{equation}
for any $t,t' \in \Z$, where $\tilde{d}$ denotes the graph distance in
$\tilde{\mathbf{m}}$. We now define the \emph{Busemann function}
$B_{\tilde{\gamma}}$ from the vertex set $\tilde{V}$ of
$\tilde{\mathbf{m}}$ to $\Z$ by
\begin{equation}
  B_{\tilde{\gamma}}(v) := \lim_{t \to \infty} \tilde{d}(v,\tilde{\gamma}(t))-t, \qquad
  v \in \tilde{V}.
\end{equation}
By the triangle inequality, one easily checks that the
limit indeed exists and is attained for $t$ large enough (for instance $B_{\tilde{\gamma}}(\tilde{\gamma}(s))=-s$ for every $s$ in $\Z$). Furthermore,
for any two adjacent vertices $v,v'$ we have
$B_{\tilde{\gamma}}(v)-B_{\tilde{\gamma}}(v')=\pm 1$ (since
$\tilde{\mathbf{m}}$ is bipartite), and every vertex $v$ has a
neighbor $v'$ for which this difference is $+1$ (i.e.
$B_{\tilde{\gamma}}$ has no local minimum). Let us denote by $T$ the
automorphism of $\tilde{\mathbf{m}}$ corresponding to the translation
$z \mapsto z+1$ of the complex plane (it corresponds
to making one turn counterclockwise around the origin in
$\mathbf{m}$): we then have
$B_{\tilde{\gamma}}(Tv)=B_{\tilde{\gamma}}(v)-s$ as a consequence of
the relation $T\tilde{\gamma}(t)=\tilde{\gamma}(t+s)$.

Now, let us denote by $(c_i)_{i \in \Z}$ and $(c'_i)_{i \in \Z}$ the
successive corners incident to the faces $\tilde{1}$ and $\tilde{2}$,
respectively, as we follow their contours walking with the face at
hand to the right. This depends on a choice for the corners $c_0$ and
$c'_0$, whose influence will be discussed later. We have
$Tc_i=c_{i+2m_1}$ and $Tc'_i=c'_{i-2m_2}$ for all $i$.

We then set, for any $i \in \Z$,
$\ell_i:=B_{\tilde{\gamma}}(c_0)-B_{\tilde{\gamma}}(c_i)$ and
$\ell'_i:=B_{\tilde{\gamma}}(c'_0)-B_{\tilde{\gamma}}(c'_i)$ (these
differences do not depend on the choice of $\tilde{\gamma}(0)$, as
claimed). Observe that the definition of $\ell_i$ is analogous to that
used in Section~\ref{sec:pointedrooted} for the decomposition of a
slice into a path decorated with elementary slices: the distance to
the apex $d(\cdot,A)$ is just replaced by the Busemann function
$B_{\tilde{\gamma}}(\cdot)$. The sequences $\ell=(\ell_i)_{i \in \Z}$
and $\ell'=(\ell'_i)_{i \in \Z}$ form infinite lattice paths (with
increments $\pm 1$) which are periodic: indeed, we have
$\ell_{i+2m_1}=\ell_i+s$ and $\ell'_{i+2m_2}=\ell'_i-s$ for all
$i$. Thus, these sequences are entirely determined by their data in a
fundamental domain: $L:=(\ell_0,\ell_1,\ldots,\ell_{2m_1})$ and
$L':=(\ell'_0,\ell'_1,\ldots,\ell'_{2m_2})$ form lattice paths
connecting $(0,0)$ to respectively $(2m_1,s)$ and $(2m_2,-s)$.

With each corner $c_i$, we associate the \emph{leftmost} infinite
geodesic $\gamma_i=(\gamma_i(t))_{t\geq 0}$ defined inductively as
follows. We let $\gamma_i(0)$ be the vertex incident to $c_i$ and,
assuming that $\gamma_i(t)$ is known, we let $\gamma_i(t+1)$ be the
leftmost vertex such that
$B_{\tilde{\gamma}}(\gamma_i(t+1))=B_{\tilde{\gamma}}(\gamma_i(t))-1$,
using the edge $\gamma_i(t-1)\gamma_i(t)$ (or the corner $c_i$ for
$t=0$) as a reference. We define in the same way the leftmost infinite
geodesic $\gamma'_i$ starting at $c'_i$. We may check that each of
these leftmost geodesics eventually merges with $\tilde{\gamma}$ (for
$\gamma'_i$, this uses the fact that we chose $\gamma$ to be the
innermost minimal separating cycle).

\begin{figure}
  \centering \includegraphics{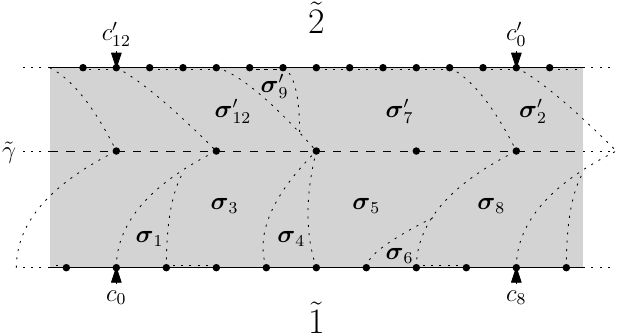}
  \caption{Sketch of the decomposition of the universal cover of an
    annular map into slices. The leftmost infinite geodesics (shown in
    dotted lines) eventually merge with the biinfinite geodesic
    $\tilde{\gamma}$ (shown in dashed lines). These geodesics delimit
    slices $\boldsymbol\sigma_i$ and $\boldsymbol\sigma'_i$,
    $i \in \Z$. In the situation displayed here, the slices
    $\boldsymbol\sigma_2$, $\boldsymbol\sigma_7$,
    $\boldsymbol\sigma'_1$, etc, are trivial slices (reduced to an
    edge).}
  \label{fig:annuliftcut}
\end{figure}
  
Let $\boldsymbol\sigma_i$ (resp.\ $\boldsymbol\sigma'_i$) be the map
delimited by $\gamma_{i-1}$ and $\gamma_i$ (resp.\ $\gamma'_{i-1}$ and
$\gamma'_i$), which we stop at their first common vertex $v_i$ (resp.\
$v'_i$). It is an elementary bipartite slice with base $c_{i-1}c_i$
(resp.\ $c'_{i-1}c'_i$). As in Section~\ref{sec:pointedrooted}, the
slice is trivial whenever $i$ corresponds to a down step of $\ell$
(resp.\ $\ell'$), and is non-trivial otherwise (but may be empty). By
periodicity, we have $\boldsymbol\sigma_i=\boldsymbol\sigma_{i+2m_1}$
and $\boldsymbol\sigma'_i=\boldsymbol\sigma'_{i-2m_2}$ (in the sense
of equality as maps). Note that every finite face, edge and vertex of
$\tilde{\mathbf{m}}$ belongs to exactly one slice
$\boldsymbol\sigma_i$ or $\boldsymbol\sigma'_i$ deprived of its right
boundary.

Let us finally discuss the roles of the reference corners $c_0$ and
$c'_0$. Changing the reference corner $c_0$ amounts to translating the
lattice path $\ell$, and to reparametrizing the sequence
$\boldsymbol\sigma=(\boldsymbol\sigma_i)_{i \in \Z}$ by a translation
of $i$. In particular, both $\ell$ and $\boldsymbol\sigma$ are
invariant if we change $c_0$ into one its translates $T^k c_0$, thus
they only depend on the choice of a corner incident to face $1$ in
$\mathbf{m}$. Similarly, $\ell'$ and
$\boldsymbol\sigma'=(\boldsymbol\sigma'_i)_{i \in \Z}$ only depend on
the choice of a corner incident to face $2$ in $\mathbf{m}$.

\label{page:label3}
Let us now assume that the map $\mathbf{m}$ has a third distinguished
element (face or vertex) labeled $3$, which is indeed the case in
the setting of Proposition~\ref{prop:generalbij}. Then it is possible to choose corners incident to
faces $1$ and $2$ in a canonical way. Namely, we pick a preimage
$\tilde{3}$ of $3$ in $\tilde{\mathbf{m}}$: $\tilde{3}$ belongs to
precisely one slice $\boldsymbol\sigma_i$ or $\boldsymbol\sigma'_i$
deprived of its right boundary. If $\tilde{3}$ belongs to a
$\boldsymbol\sigma_i$, then by changing the reference corner $c_0$ we
can ensure that it belongs to $\boldsymbol\sigma_1$. Then, we choose
$c'_0$ in such way that
$B_{\tilde{\gamma}}(c'_0)=B_{\tilde{\gamma}}(c_0)$: this is possible
because the function $B_{\tilde{\gamma}}(\cdot)$ has $\pm 1$
increments along the contour of $\tilde{2}$ and decreases by $s$ over
a translation $T$, hence it is surjective. There might exist several
such $c'_0$, in which case we pick the last one, so that $B_{\tilde{\gamma}}(c'_i)\neq B_{\tilde{\gamma}}(c'_0)$ for every
$i>0$. 
If $\tilde{3}$ belongs to a $\boldsymbol\sigma'_i$, we proceed in the
same way upon exchanging the roles of $c_0$ and $c'_0$.

\begin{lemma}
  \label{sec:decomp-an-annul}
The above construction is a bijection between: 
\begin{itemize}
\item 
the set of maps $\mathbf{m}$ with two marked faces of degrees
$2m_1,2m_2\in \Z_{>0}$ and all other faces of even degree, with a third 
distinguished face or vertex, and with
separating girth $s\in \Z_{>0}$,
\item
 the set of  quadruples
$(L,L',\boldsymbol\sigma,\boldsymbol\sigma')$, where $L$ and $L'$ are
lattice paths connecting $(0,0)$ to respectively $(2m_1,s)$ and
$(2m_2,-s)$, where
$\boldsymbol\sigma$ and $\boldsymbol\sigma'$ are sequences of
elementary slices of respective periods $2m_1$ and $2m_2$, such that
for every $i=1,\ldots,2m_1$ (resp.\ $i=1,\ldots,2m_2$) the tilt of
$\boldsymbol\sigma_i$ and (resp.\ $\boldsymbol\sigma'_i$) is equal to
$\ell_i-\ell_{i-1}$ (resp.\ $\ell'_i-\ell'_{i-1}$), where either $\boldsymbol\sigma_1$ 
or $\boldsymbol\sigma'_1$ carries a distinguished face or vertex, and
where $L'$ (resp.\ $L$) does not return to height $0$ if $\boldsymbol\sigma_1$
(resp.\ $\boldsymbol\sigma'_1$) 
carries the distinguished element. 
\end{itemize}
\end{lemma}

To prove this lemma, we now describe the
inverse procedure which consists in assembling an annular map
$\mathbf{m}$ from a quadruple
$(L,L',\boldsymbol\sigma,\boldsymbol\sigma')$.

\paragraph{Assembling an annular map from a pair of paths decorated
  with elementary slices.}

Using the lattice path $L$ and the sequence of slices
$\boldsymbol\sigma$ restricted to a period, we may use the inverse of
the decomposition described in Section~\ref{sec:pointedrooted} to
obtain a slice $\boldsymbol\Sigma$ of width $2m_1$ and tilt $s$. We
then perform the operation described in~\cite[Section~7.1]{irredmaps}
and called \emph{wrapping} in \cite[Section 2.2]{BouttierHDR}. It
consists in gluing ``isometrically'' the left and right boundaries of
$\boldsymbol\Sigma$ together, but unlike the gluing performed in
Section~\ref{sec:pointedrooted} to obtain a pointed rooted map, we now
start from the two endpoints of the base of $\boldsymbol\Sigma$ which
we identify together, and perform the gluing from there (see the
illustrations in the aforementioned references). As the tilt $s$ is
positive, the $s$ edges of the left boundary closest to the apex
remain unglued. This produces an annular map $\mathbf{m}_1$ whose two
distinguished faces have respective degrees $2m_1$ and $s$, the latter
resulting from the unglued edges, the former resulting from the base,
whose two identified endpoints yield a distinguished corner denoted
$c$.

Similarly, the lattice path $L'$ and the sequence of slices
$\boldsymbol\sigma'$ yield a slice $\boldsymbol\Sigma'$ of width
$2m_2$ and tilt $-s$, whose wrapping produces an annular map
$\mathbf{m}_2$ with two distinguished faces of respective degrees
$2m_2$ and $s$ (since the tilt of $\boldsymbol\Sigma'$ is negative, it
is now the $s$ edges closest to the apex on the \emph{right} boundary
that form this latter face), the former having a distinguished corner
denoted $c'$.

The annular map $\mathbf{m}$ is then obtained by gluing the two
annular maps $\mathbf{m}_1$ and $\mathbf{m}_2$ together along the
contours of their distinguished faces of degree $s$. The contour then
becomes a separating cycle $\gamma$, which can be shown to have
minimal length. Note that there are a priori $s$ ways to glue
$\mathbf{m}_1$ and $\mathbf{m}_2$ together, however only one way is
compatible with the prescription that, in the universal cover
$\tilde{\mathbf{m}}$ of $\mathbf{m}$, the corners $c$ and $c'$ admit
respective lifts $c_0$ and $c'_0$ such that
$B_{\tilde{\gamma}}(c_0)=B_{\tilde{\gamma}}(c'_0)$ where
$\tilde{\gamma}$ is the infinite geodesic constructed by lifting
$\gamma$.

We may check that the assembling procedure is indeed the inverse of
the decomposition described above. The most subtle point is to show
that any quadruple $(L,L',\boldsymbol\sigma,\boldsymbol\sigma')$ is
left invariant if we assemble it then decompose the resulting annular
map. This requires to check that the boundaries of the slices become,
in the annular map, precisely the leftmost geodesics that we use in
the decomposition. More details can be found in \cite[Section~7]{irredmaps}.

\paragraph{Application to tight maps and end of the proof of
  Proposition~\ref{prop:generalbij}.}

So far, all our discussion holds without the assumption that
$\mathbf{m}$ is tight. Adding this constraint amounts to two
restrictions on the corresponding quadruple
$(L,L',\boldsymbol\sigma,\boldsymbol\sigma')$, which are similar to
those encountered in Section~\ref{sec:pointedrooted}. First, every
$\boldsymbol\sigma_i$ and $\boldsymbol\sigma'_i$ must be a tight
elementary slice. Here, we use the natural convention that the mark of
a marked vertex is transferred to its copy in the unique slice
deprived of its right boundary that contains it (in particular, marked
vertices lying on the separating cycle $\gamma$ have their marks
transferred to $\boldsymbol\sigma_i$'s and not to
$\boldsymbol\sigma'_i$'s). Second, the lattice paths $\ell$ and
$\ell'$ (which are the periodic extensions of $L$ and $L'$) must be
such that every up step immediately following a down step cannot be
decorated with an unmarked empty slice.

To complete the proof of Proposition~\ref{prop:generalbij}, we assume
that $\mathbf{m}$ has $n\geq 3$ boundaries, labeled $1$ to $n$.  By
the slice decomposition, the labels $3,\ldots,n$ get distributed among
the $\boldsymbol\sigma_i$ and $\boldsymbol\sigma'_i$. Let
$i_1 \leq \cdots \leq i_j$ be the indices $i$ between $1$ and $2m_1$
such that $\boldsymbol\sigma_i$ is neither the unmarked empty slice
nor the trivial slice, i.e.\ contains at least one label. Let
$i'_1 \leq \cdots \leq i'_{j'}$ be similarly the indices $i$ between
$1$ and $2m_2$ such that $\boldsymbol\sigma'_i$ contains at least one
label. We set $k=j+j'-1$ and note that $0 \leq k \leq n-3$ since there
are $n-2$ labels in total to distribute among the slices. Recall that,
by the aforementioned prescription for choosing the reference corners,
the label $3$ is either in $\boldsymbol\sigma_1$ or
$\boldsymbol\sigma'_1$. If it is in $\boldsymbol\sigma_1$, we set
$(\mathbf{s}_1,\ldots,\mathbf{s}_{k+1}):=(\boldsymbol\sigma_{i_1=1},\ldots,\boldsymbol\sigma_{i_j},\boldsymbol\sigma'_{i'_1},\ldots,\boldsymbol\sigma'_{i'_{j'}})$. Otherwise,
we set
$(\mathbf{s}_1,\ldots,\mathbf{s}_{k+1}):=(\boldsymbol\sigma'_{i'_1=1},\ldots,\boldsymbol\sigma'_{i'_{j'}},\boldsymbol\sigma_{i_1},\ldots,\boldsymbol\sigma_{i_j})$. This
defines the $(k+1)$-tuple of slices $\mathbf{s}$ of the proposition.

As for the two-face map $\mathbf{m}_{12}$, it is obtained as
follows. Let $\boldsymbol\omega$ and $\boldsymbol\omega'$ be the
sequences obtained from respectively $\boldsymbol\sigma$ and
$\boldsymbol\sigma'$ by replacing every slice different from the
unmarked empty slice and from the trivial slice with the marked
empty slice. Then, $\mathbf{m}_{12}$ is the annular map obtained by
assembling the quadruple
$(L,L',\boldsymbol\omega,\boldsymbol\omega')$. It is by construction a
tight two-face map with $k+1$ marked vertices, and we distinguish the
marked vertex coming from the marked empty slice which replaces the
slice containing the label $3$.

We check that the mapping
$\mathbf{m} \mapsto (\mathbf{m}_{12},\mathbf{s})$ is a bijection by
exhibiting the inverse bijection. Let $(\mathbf{m}_{12},\mathbf{s})$
be a pair as in the proposition, and let
$(L,L',\boldsymbol\omega,\boldsymbol\omega')$ be the quadruple
obtained by decomposing the two-face map $\mathbf{m}_{12}$, its
distinguished marked vertex playing the role of the third
distinguished element labeled $3$ used in the construction of
page~\pageref{page:label3}. By construction the sequences
$\boldsymbol\omega$ and $\boldsymbol\omega'$ consist only of empty or
trivial slices, with a number $k+1$ of marked empty slices.  Replacing
these marked empty slices with $\mathbf{s}_1,\ldots,\mathbf{s}_{k+1}$,
we obtain two sequences $\boldsymbol\sigma$ and $\boldsymbol\sigma'$
such that the assembling of the quadruple
$(L,L',\boldsymbol\sigma,\boldsymbol\sigma')$ gives the tight map
$\mathbf{m}$ we are looking for. This ends the proof of Proposition~\ref{prop:generalbij}. \qed

\begin{remark}
  It might seem more direct to attempt to enumerate the
  quadruples $(L,L',\boldsymbol\sigma,\boldsymbol\sigma')$. But, since
  the paths $L$ and $L'$ have a height variation depending on the
  separating girth $s$, this leads to an expression involving a sum
  over $s$. The trick of ``recombining'' $L$ and $L'$ into a two-face
  map allows to circumvent this issue.
\end{remark}

\section{Bijective proofs for non necessarily bipartite maps}\label{sec:extens-non-bipart}

In this section, we explain how the bijective approach of Section~\ref{sec:bijective} may be extended so
as to enumerate planar tight maps which are not necessarily bipartite.

Recall from Section \ref{sec:quasipolmain}
that a \emph{petal} is a face of degree one. The key idea to
extend our construction  
to the non bipartite case is to realize that petals play a role very
similar to marked vertices. 
Mimicking the organization of Section~\ref{sec:bijective}, we will first enumerate
tight maps with a single non-petal face in Section~\ref{sec:petaltrees}, then tight maps with just two non-petal faces
in Section~\ref{sec:petalnecklaces}, 
and finally tight maps with an arbitrary number of non-petal faces in Section~\ref{sec:generalnonbip}.
A prerequisite to this latter enumeration will be that of tight non
necessarily bipartite slices in Section~\ref{sec:generalslices}.

Before we start our discussion, let us recall from Section \ref{sec:quasipolmain} the definition of the
polynomials
\begin{equation}
   p_{k,e}(m):=\frac{1}{(k!)^2}\,
   \prod_{i=1}^k\left(m^2-\left(i-\frac{e}2\right)^2\right) =
   \binom{m+\frac{e}2-1}{k} \binom{m-\frac{e}2+k}{k}
\end{equation}
for $k\in \Z_{\geq 0}$.
For $m - \frac{e}2 \in \Z$, we may interpret $p_{k,e}(m)$ as counting words of
the form~\eqref{eq:UDform1} where the $a_i$ and $b_i$ are nonnegative
integers such that $a_1 + a_2 + \cdots + a_{k+1} = m+\frac{e}2-k-1$
and $b_1 + b_2 + \cdots + b_{k+1} = m -\frac{e}2$ (i.e.\ the word has
length $2m$, $k+1$ occurrences of $D_\circ$ and $e$ more $D$'s than
$U$'s). This interpretation holds a priori only for $m\geq \max\left(k+1-\frac{e}2,\frac{e}2\right)$
but this domain may be extended to $m\geq \max\left(1-\frac{e}2,\frac{e}2-k\right)$ since $p_{k,e}(m)$ vanishes 
in the additional domain.

\begin{remark}
Even though we will not use it in the sequel, let us mention that, by a variation of the arguments of Section~\ref{sec:oneface}, we may
show that, for $m$ as above, $p_{k,e}(m)$ is the number of
two-face tight maps with one face of degree $2m$, one simple face of
degree $|e|$ and $k+1$ marked vertices, one of them distinguished,
with the condition that for $e<0$ no marked vertex is incident to the
face of degree $e$. Note that such maps are closely related with the notion of
$(a,b)^*$-forests of Definition~\ref{sec:bipartite-case.--1}: for
$e>0$, $p_{k,e}(m)$ counts $(0,e)^*$-forests with size $2m$ and $k+1$
marked vertices including the distinguished vertex. 
\end{remark}

\subsection{Petal trees}
\label{sec:petaltrees}

We recall that a \emph{petal tree} is a planar map having an
\emph{exterior} face of arbitrary degree, and such that every other
face is a petal. A tight
petal tree is just a petal tree with marked vertices, which is tight as
a map, see again Figure~\ref{fig:petaltree}.

We also recall, for $r,s\in \Z_{\geq 0}$, $\epsilon\in \Z$
and $m\in \Z/2$, the definition of the quasi-polynomial
  \begin{equation}
    \label{eq:pieps2eme}
    \pi_{r,s}^{(\epsilon)}(m) :=
    \begin{cases}
      \binom{r+s}s p_{r+s,s+1+\epsilon}(m) & \text{if $m-\frac{s+1+\epsilon}2 \in \Z$,} \\
      0 & \text{otherwise.}
    \end{cases}
  \end{equation}
Our goal in this section is to prove Proposition \ref{prop:pieps},
showing that, for $m \in \Z_{> 0}/2$ and $\epsilon \in \{-1,0,1\}$, $\pi^{(\epsilon)}_{r,s}(m)$ enumerates tight petal trees
with one exterior face of degree $2m$, $s+1+\epsilon$ petals
(excluding the exterior face when $m=1/2$), $1+\epsilon$ of
which are distinguished, and $r+1-\epsilon$ marked vertices, $1-\epsilon$ of which
are distinguished. 
The reason for the $\epsilon$-dependence in the  statement is that we
must distinguish two elements  
among marked vertices and petals. For this reason, there are three
situations to consider: we may distinguish two marked vertices ($\epsilon=-1$), two
petals ($\epsilon=1$), or one of each type ($\epsilon=0$).

We note that, upon ``labeling'' the 
undistinguished elements, Proposition~\ref{prop:pieps} is equivalent to the more
symmetric statement: 
\begin{proposition}
  \label{prop:Pkeid}
  For $0 \leq e \leq k+2$ and $m\in \Z_{> 0}/2$, the number of tight petal
  trees with an exterior face of degree $2m$, with $e$ petals and
  $k+2-e$ marked vertices, all labeled, is given by:
  \begin{equation}
    \label{eq:Pkeid}
    N_{0,k+3}(2m,\underbrace{1,\ldots,1}_{e},\underbrace{0,\ldots,0}_{k+2-e}) =
    \begin{cases}
      k! p_{k,e}(m) & \text{if $m-\frac{e}2 \in \Z$,} \\
      0 & \text{otherwise.}
    \end{cases}
  \end{equation}
\end{proposition}
More precisely, Proposition~\ref{prop:Pkeid} is recovered from Proposition~\ref{prop:pieps} by setting $k=r+s$ and $e=s+1+\epsilon$
(hence $k+2-e=r+1-\epsilon$). The petal trees considered in Proposition~\ref{prop:pieps} have $s$ (respectively $r$) non-distinguished
petals (respectively marked vertices) which we may label in $s!$ (respectively $r!$) ways. Using
$s! r! \binom{r+s}s=(r+s)!=k!$, we obtain \eqref{eq:Pkeid} for $1+\epsilon\leq e\leq k+1+\epsilon$, hence for the whole range $0 \leq e \leq k+2$ by letting in $\epsilon$ vary in $\{-1,0,1\}$.

\begin{proof}[Proof of Proposition~\ref{prop:pieps}] Let us first discuss how we may code petal trees using words, 
or equivalently lattice paths, as we did in Section~\ref{sec:oneface} for ordinary trees. Note first that petal trees have 
two types of edges: \emph{``tree-type''} edges whose both sides are incident to the exterior face and \emph{``petal-type''} edges 
with one edge side incident to a petal and the other to the exterior face. In particular, the tree-type edges have distinct
endpoints and they form a plane tree (i.e.~a map with a single face), which we call the \emph{``wood''} of the petal tree. 
Here the coding that we shall use applies to \emph{rooted} petal trees, i.e.~petal trees where we distinguish a corner (the root corner)
\emph{in the exterior face} (note that this also induces a rooting of the wood tree). 
Starting from this root corner and following the contour of the exterior face going counterclockwise around 
the tree  (i.e.~with the exterior face on the right), we record a letter $U$ (respectively $D$) for each tree-type edge visited while going away from (respectively towards) the root and a letter $E$ for each visited petal-type edge. The obtained three-letter word may alternatively be visualized as a lattice path with three types of elementary steps: up, down and horizontal, associated
respectively to the letters $U$, $D$ and $E$. This path is a Motzkin path of length equal to the degree $2m'$ of the exterior face, 
i.e.~it goes from $(0,0)$ to $(2m',0)$ and \emph{stays above the
  $x$-axis}. Indeed, the ordinates of the path are \emph{non-negative} since they record graph distances
to the \emph{root vertex} incident to the root corner.
 The number $s'$ of occurrences 
of the letter $E$ in the coding word is nothing but the number of petals, and is such that $m'-\frac{s'}2$, which is the number of $U$'s (or equivalently of $D$'s), is a nonnegative integer. 

Let us now consider a rooted petal tree with marked vertices, \emph{where
the root vertex is unmarked}, and such that all non-root leaves are marked. Every non-root vertex 
being bijectively associated with its parent tree-type edge in the wood tree, hence with a letter $D$, we
record the markings as in Section~\ref{sec:oneface} by replacing each $D$ by a $D_\circ$ if the associated vertex is marked
and by a $D_\bullet$ otherwise. We end up with a word made of four letters, $U$, $D_\circ$, $D_\bullet$ and $E$, associated
with a ``dressed'' Motzkin path, 
with $k'$ occurrences of $D_\circ$ and $m'-\frac{s'}2-k'$ occurrences
of $D_\bullet$ if the petal tree has $k'$ marked vertices.
As before, requiring that all non-root leaves be marked simply amounts to forbidding the sequence $UD_\bullet$ in the coding
words.

 \begin{figure}[t]
    \centering
    \includegraphics[scale=.9]{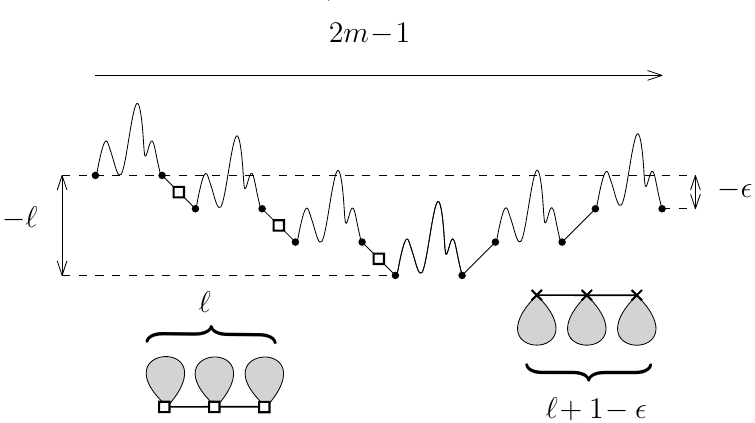}
    \caption{Schematic picture of the decomposition of a lattice path of length $2m-1$, height difference $-\epsilon$, 
and minimal height $-\ell $ into a sequence of $2\ell+1-\epsilon$ Motzkin paths, where the squares
indicate that the underlying descending step may be marked (if associated with a $D_\circ$) or not (if associated
with a $D_\bullet$). For $j=1,2,\cdots,\ell$, the $j$-th Motzkin path is followed by the first elementary down step reaching height
 $-j$. These $\ell$ Motzkin paths code for an ordered set $\mathcal{F}_\square$ of $\ell$ rooted petal trees whose root
vertex are marked or not according to the nature of the following elementary down step. For $j=\ell+2,\ell+3,\cdots,2\ell+1-\epsilon$,
the $j$-th Motzkin path is preceded by the last elementary up step reaching height $j-1-2\ell$.
Together with the $(\ell+1)$-th Motzkin path, these Motzkin paths code for an ordered set $\mathcal{F}_\times$ of $\ell+1-\epsilon$ 
rooted petal trees whose root vertex is unmarked. The ordered sets $\mathcal{F}_\square$ and $\mathcal{F}_\times$ are drawn
here with additional edges connecting the successive roots of their petal tree components, represented schematically 
by grey blobs with a black boundary, with a root vertex drawn as a square if markable and a cross if not. }
    \label{fig:piepsilon}
  \end{figure}

\begin{figure}[t]
    \centering
    \includegraphics[scale=.9]{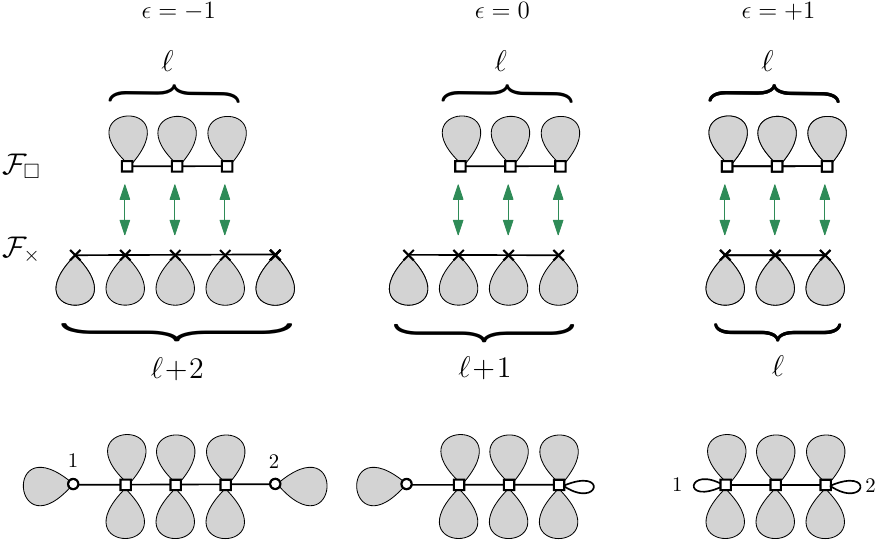}
    \caption{The gluing of a sequence $\mathcal{F}_\square$ of $\ell$ petal trees with markable roots and a sequence
    $\mathcal{F}_\times$ of $\ell+1-\epsilon$ petal trees with unmarked roots into a tight petal tree with $1-\epsilon$ 
    distinguished marked vertices (represented by circles) and $1+\epsilon$ distinguished petals (represented
    by small empty loops), for $\epsilon=-1$ (left), $0$ (center) and $1$ (right).}
    \label{fig:piepsilon-101}
  \end{figure}

Returning to the setting of Proposition~\ref{prop:pieps} and assuming that $m-\frac{s+1+\epsilon}2 \in \Z$, recall that 
$p_{r+s,s+1+\epsilon}(m)$ counts three-letter words of
the form~\eqref{eq:UDform1} with $r+s+1$ occurrences of $D_\circ$, $m-r-\frac{s+1-\epsilon}2$
occurrences of $D_\bullet$ and $m-\frac{s+1+\epsilon}2$ occurrences of $U$, so that the height difference is $-(s+1+\epsilon)$. As already mentioned, the property holds for 
the extended range $m\geq \max\left(1-\frac{s+1+\epsilon}2,1-r-\frac{s+1-\epsilon}2\right)$, hence for all positive $m$ 
in the current setting where $r,s\geq 0$ and $\epsilon\in \{-1,0,1\}$. We now transform the three-letter word into a four-letter word coding for a \emph{sequence of petal trees} 
as follows: we first remove the first letter $D_\circ$, which leaves us with $r+s$ occurrences of $D_\circ$, 
and pick $s$ of them that we transform into $E$'s. This transformation can be done in $\binom{r+s}s$ ways, hence the obtained 
four-letter words are counted by $\binom{r+s}s p_{r+s,s+1+\epsilon}(m)$. By construction, these words have $s$ occurrences 
of $E$, $r$ occurrences of $D_\circ$, $m -\frac{s+1+\epsilon}2$ occurrences of $U$ and  $m-r-\frac{s+1-\epsilon}2$ 
occurrences of $D_\bullet$. These correspond to lattice  paths of length $2m-1$ and height difference $-\epsilon$, 
say from $(0,0)$ to $(2m-1,-\epsilon)$, hence with minimal height $-\ell $ for some 
$\ell\geq \max(0,\epsilon)$. Any such path is canonically decomposed into a sequence of $2\ell+1-\epsilon$ dressed 
Motzkin paths
(possibly of length $0$), obtained by cutting out the first elementary down steps reaching height $-j$ for $j=1,2,\ldots,\ell$ and the last elementary 
up steps reaching height $-j'+1$  for $j'=\ell,\ell-1,\ldots,\epsilon+1$, see Figure~\ref{fig:piepsilon}. Each
Motzkin path component codes for a marked petal tree whose root vertex is unmarked, and we decide to mark it if
this Motzkin path is followed by a $D_\circ$ in the original path. By doing so, we both ensure that the total number of marked 
vertices in the petal tree sequence is $r$ and that the four-letter word can be recovered bijectively from the petal tree sequence
(since we know the nature $D_\bullet$, $D_\circ$ or $U$ of all the removed steps). Note that by construction,
only the first $\ell$ petal trees may have their root vertex marked.
All in all, $\binom{r+s}s p_{r+s,s+1+\epsilon}(m)$ enumerates pairs $(\mathcal{F}_\square,\mathcal{F}_\times)$ made of a sequence 
$\mathcal{F}_\square$ of $\ell$ petal trees
with marked vertices whose root vertex is markable, and a sequence $\mathcal{F}_\times$ of $\ell+1-\epsilon$ petal trees
with marked vertices whose root vertex is unmarked, for some arbitrary integer $\ell\geq \max(0,\epsilon)$, with a
total of $s$ petals, $r$ marked vertices, and $m-\ell-\frac{s+1-\epsilon}2$ tree-type edges, and with all the non-root leaves in the petal trees marked.
As a final step, the pairs $(\mathcal{F}_\square,\mathcal{F}_\times)$ are transformed into the desired tight petal trees by attaching the roots 
of the $\ell$ (respectively $\ell+1-\epsilon$) petal trees in the sequence $\mathcal{F}_\square$ (respectively $\mathcal{F}_\times$)
by additional edges and gluing them along a spine of length $\ell-\epsilon$ as displayed in Figure~\ref{fig:piepsilon-101}.
Note that a markable root vertex is always matched to an unmarked one: the markings of root vertices may thus be transferred without ambiguity after gluing, with no risk of double markings along the spine. Let us detail the three possibilities
$\epsilon=-1,0,1$.

For $\epsilon=-1$, the sequence
$\mathcal{F}_\times$ has two more petal trees than $\mathcal{F}_\square$ so that the root vertices of its two extremal petal trees
remain unmatched.
We decide to mark these vertices and distinguish them as vertices $1$ and $2$: 
the resulting object is a tight petal tree with $s$ petals and $r+2$ marked vertices, two of which distinguished. 
The construction is clearly reversible by cutting along the branch 
between vertices $1$ and $2$, and holds for any $\ell\geq 0$ (for $\ell=0$, vertices $1$ and $2$ are incident).
The degree of the exterior face is easily found equal to $2(m-\ell-\frac{s+2}2)+s+2(\ell+1)=2m$ as wanted.

For $\epsilon=0$, the sequence
$\mathcal{F}_\times$ has one more petal tree than $\mathcal{F}_\square$ so that the root vertex of its first petal tree remains unmatched.
We decide to mark this vertex and to add an additional, distinguished petal to mark the corner inbetween the last two glued 
petal trees at
the end of the spine: the resulting object is a tight petal tree with $s+1$ petals, one of which distinguished, and $r+1$ marked vertices, one of which distinguished. The construction is clearly reversible by cutting along the branch 
between the distinguished marked vertex and the distinguished petal, eventually removing this
petal, and holds for any $\ell\geq 0$ (for $\ell=0$, the distinguished vertex is incident to the distinguished petal). 
The degree of the exterior face is $2(m-\ell-\frac{s+1}2)+s+1+2\ell=2m$ as wanted.

Finally, for $\epsilon=1$, the sequences
$\mathcal{F}_\times$ and $\mathcal{F}_\square$ have the same number $\ell\geq 1$ of petal tree components
and their gluing is made reversible by adding a petal inbetween the glued petal trees at each extremity of the spine,
that we distinguish as petals $1$ and $2$: the resulting object is a tight petal tree with $s+2$ petals, two of which distinguished, and $r$ marked vertices. The construction is clearly reversible by cutting along the branch 
between petals $1$ and $2$, eventually removing these
petals, and holds for any $\ell\geq 1$ (for $\ell=1$, the distinguished petals
are incident to the same vertex). The degree of the exterior face is $2(m-\ell-\frac{s}2)+s+2+2(\ell-1)=2m$ as wanted.
This ends the proof of Proposition~\ref{prop:pieps}.
\end{proof}

\subsection{Petal necklaces} 
\label{sec:petalnecklaces}
We now consider \emph{petal necklaces},
namely planar maps having two distinguished faces of arbitrary
degrees, and such that any other face is a petal. Again, a tight petal necklace is just a petal necklace with
marked vertices, which is tight as a map. We have the following 
enumeration result:
  \begin{figure}[htbp]
    \centering
    \includegraphics[width=\textwidth]{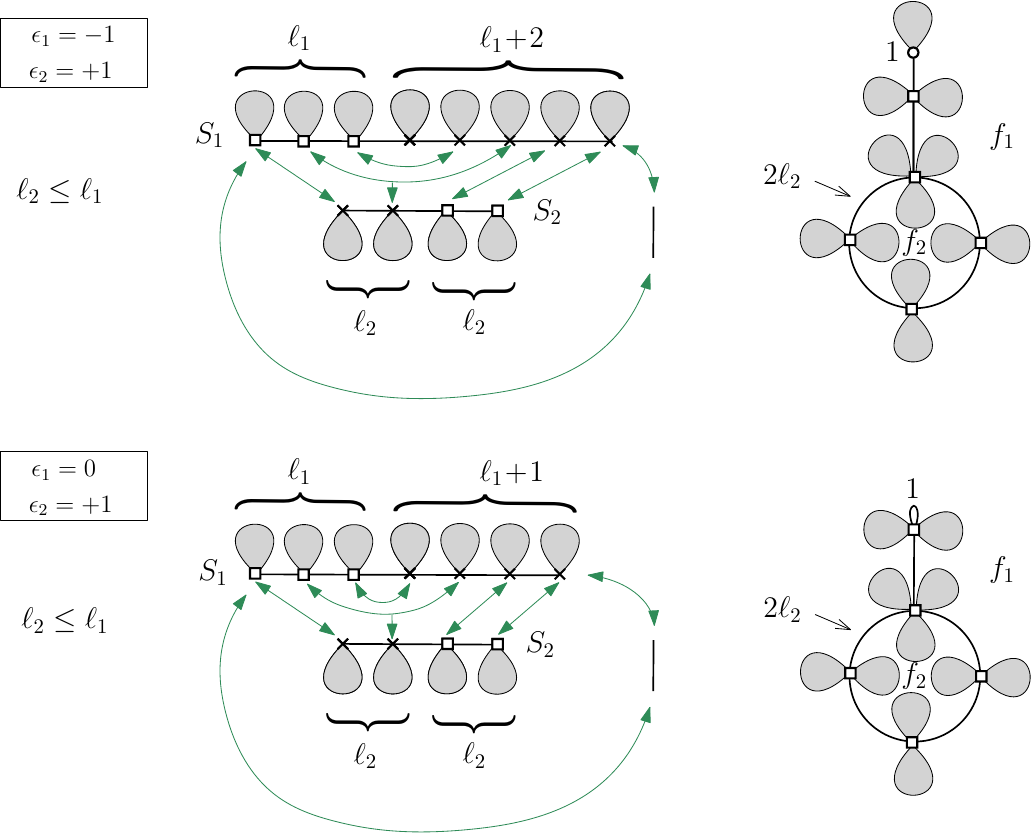}
    \caption{The partial and mutual gluings of the segments $S_1$ (oriented from left to right) and $S_2$ (oriented from right to left) connecting 
    the petal tree components of the sequences $\mathcal{F}_1$ and 
    $\mathcal{F}_2$ respectively, here for $\ell_2\leq \ell_1$. The vertex identifications are indicated by arrows.
    Top: in the case $\epsilon_1=-1$, $\epsilon_2=+1$, the partial gluing  of $S_1$ amounts to pulling in its $(\ell_1+1)$-th vertex,
     creating a branch of length $\ell_1-\ell_2+1$ in $f_1$ with a distinguished vertex $1$ at its end. Bottom: in the case $\epsilon_1=0$, $\epsilon_2=+1$, the partial gluing of $S_1$ amounts to pulling in its $\ell_1$-th edge,
creating a branch of length $\ell_1-\ell_2$ in $f_1$ with a distinguished petal $1$ at its end. In both cases, the mutual gluing generates a loop
    of even length $2\ell_2$. 
     }
    \label{fig:twofaces-1101}
  \end{figure}
\begin{proposition}
  \label{prop:deuxfacesnonbip}
  For $r,s$ nonnegative integers not both zero and $m_1,m_2 \in \Z_{>0}/2$, let
  \begin{equation}
    \Pi_{r,s}(m_1,m_2):=N_{0,r+s+2}(2m_1,2m_2,\underbrace{0,\ldots,0}_{r},\underbrace{1,\ldots,1}_{s})
  \end{equation}
  be the number of tight petal necklaces with two distinguished faces
  of degrees $2m_1$ and $2m_2$, with $r$ marked vertices and $s$
  petals, all labeled.  Then, for $r \geq 1$ we have
  \begin{equation}
    \label{eq:Pirs1}
    \Pi_{r,s}(m_1,m_2) = (r-1)!s! \sum_{\substack{r_1,r_2,s_1,s_2\geq
        0 \\
        r_1+r_2=r-1 \\
        s_1+s_2=s}}
    \left( \pi_{r_1,s_1}^{(-1)}(m_1)
      \pi_{r_2,s_2}^{(1)}(m_2) + \pi_{r_1,s_1}^{(0)}(m_1)
      \pi_{r_2,s_2}^{(0)}(m_2) \right)
  \end{equation}
  and, for $s \geq 1$,
  \begin{equation}
    \label{eq:Pirs2}
    \Pi_{r,s}(m_1,m_2) = r! (s-1)! \sum_{\substack{r_1,r_2,s_1,s_2\geq
        0 \\ r_1+r_2=r \\ s_1+s_2=s-1}}
    \left( \pi_{r_1,s_1}^{(0)}(m_1)
      \pi_{r_2,s_2}^{(1)}(m_2) + \pi_{r_1,s_1}^{(1)}(m_1)
      \pi_{r_2,s_2}^{(0)}(m_2) \right)
  \end{equation}
  where $\pi_{r,s}^{(\epsilon)}$ is the univariate quasi-polynomial defined in~\eqref{eq:pieps}.
\end{proposition}

This proposition admits a ``partially unlabeled'' equivalent, which
will be useful later on.

\begin{proposition}
  \label{prop:deuxfacesnonbip-unlab}
  Given $r_0,s_0 \in \Z_{\geq 0}$ and $m_1,m_2 \in \Z_{>0}/2$, the
  number of tight petal necklaces with two distinguished labeled faces of
  degrees $2m_1,2m_2$, one distinguished \underline{marked vertex},
  $r_0$ other unlabeled marked vertices, and $s_0$ unlabeled petals,
  is equal to
  \begin{equation}
    \label{eq:pirs1}
    \pi^{(0)}_{r_0,s_0}(m_1,m_2) :=  
    \sum_{\substack{r_1,r_2,s_1,s_2\geq
        0 \\ r_1+r_2=r_0 \\ s_1+s_2=s_0}}
    \left( \pi_{r_1,s_1}^{(-1)}(m_1)
      \pi_{r_2,s_2}^{(1)}(m_2) + \pi_{r_1,s_1}^{(0)}(m_1)
      \pi_{r_2,s_2}^{(0)}(m_2) \right),
  \end{equation}
  while the number of tight petal necklaces with two distinguished
  labeled faces of degrees $2m_1,2m_2$, one distinguished
  \underline{petal}, $s_0$ other unlabeled petals, and $r_0$ unlabeled
  marked vertices, is equal to
  \begin{equation}
    \label{eq:pirs2}
    \pi^{(1)}_{r_0,s_0}(m_1,m_2) := \sum_{\substack{r_1,r_2,s_1,s_2\geq
        0 \\ r_1+r_2=r_0 \\ s_1+s_2=s_0}}
    \left( \pi_{r_1,s_1}^{(0)}(m_1)
      \pi_{r_2,s_2}^{(1)}(m_2) + \pi_{r_1,s_1}^{(1)}(m_1)
      \pi_{r_2,s_2}^{(0)}(m_2) \right).
  \end{equation}
\end{proposition}

Propositions~\ref{prop:deuxfacesnonbip}
and~\ref{prop:deuxfacesnonbip-unlab} are indeed equivalent to one
another, since the petal necklaces considered in the latter have three
distinguished boundaries hence no symmetries. Namely, we pass
from~\eqref{eq:pirs1} to \eqref{eq:Pirs1} by taking $r=r_0+1$ and
$s=s_0$, and from~\eqref{eq:pirs2} to \eqref{eq:Pirs2} by taking
$r=r_0$ and $s=s_0+1$. In both cases, we multiply by $r_0! s_0!$ to
label the unlabeled marked vertices/petals in all possible ways.
For later use we record the following:

\begin{remark}
  \label{rem:deuxfacesnonbip-lab}
  Given a petal necklace $\mathbf{m}_{12}$ with a distinguished
  marked vertex/petal as considered in
  Proposition~\ref{prop:deuxfacesnonbip-unlab}, it is possible to
  label the marked vertices and petals in a \emph{canonical} way from
  $1$ to $r_0+s_0+1$, with the label $1$ assigned to the distinguished
  marked vertex/petal. The precise labeling procedure is
  irrelevant, as long as it is deterministic\footnote{For instance, calling
  $f_1$ and $f_2$ the two distinguished faces, a possible algorithm
  consists in picking the vertex $v$ on the boundary between $f_1$ and
  $f_2$ lying on the same branch (possibly of length $0$) in the petal
  necklace as the distinguished marked vertex/petal. Following the
  contour of $f_1$ clockwise, starting from the leftmost corner at $v$
  in this face, and doing then the same for $f_2$, we label the marked
  vertices and petals, including the distinguished element, by
  successive integers from $1$ to $r_0+s_0+1$ at the \emph{first
    visit} of such marked vertex (encountered via an incident corner)
  or petal (encountered via a petal-type edge). We then perform a
  cyclic permutation of the labels so that the distinguished element
  receives the label $1$.}.
\end{remark}

\begin{proof}[Proof of Proposition \ref{prop:deuxfacesnonbip-unlab}. ]
Our proof is a direct generalization of that of Proposition
\ref{prop:pk2interpret} given in Section~\ref{sec:twofaces}. Still, rather than
recoursing to $(a,b)$- and $(a,b)^*$-forests or generalizations
thereof, we will instead use the  
related notion of pairs of petal tree sequences\footnote{There is
  indeed a clear correspondence between $(a,b)$-forests and
  collections of trees with $a$ unmarkable roots and $b$ markable
  ones, and this could be generalized to our setting, with petal
  trees replacing trees.}, as introduced in the proof of Proposition~\ref{prop:pieps}. 
More precisely, we have seen there that, for $r_1,s_1$ nonnegative integers, $m_1\in \Z_{> 0}/2$, and
 $\epsilon_1=-1,0,1$, the quantity $\pi_{r_1,s_1}^{(\epsilon_1)}(m_1)$ defined in \eqref{eq:pieps} enumerates 
pairs $(\mathcal{F}_\square,\mathcal{F}_\times)$ of sequences of respectively $\ell_1$ and 
$\ell_1+1-\epsilon_1$ rooted petal trees for some arbitrary integer $\ell_1\geq \max(0,\epsilon_1)$. 
The petal trees have a total of $r_1$ marked vertices, with all their non-root leaves marked, a total of $s_1$ petals, and 
a total of $m_1-\ell_1-\frac{s_1+1-\epsilon_1}2$ tree-type edges. Finally, only the petal trees in $\mathcal{F}_\square$ may have their root vertex marked. The enumeration statement above holds even if $m_1-\frac{s+1+\epsilon_1}2$ is a half-integer, as the number of sequence pairs is trivially $0$ in this case, and so is $\pi_{r_1,s_1}^{(\epsilon_1)}(m_1)$ by definition. 

To construct a petal necklace, we first merge the two sequences $\mathcal{F}_\square$ and $\mathcal{F}_\times$
into a single sequence $\mathcal{F}_1$ of $2\ell_1+1-\epsilon_1$ petal trees, which we transform into a single connected object
by adding edges connecting the successive roots of the petal tree components and forming a linear segment $S_1$ of length $2\ell_1-\epsilon_1$  -- see Figure~\ref{fig:twofaces-1101} for a schematic representation. This first connected 
object will code for the face of degree $2m_1$ of our necklace. The coding of the face of degree $2m_2$ is performed
via a second sequence $\mathcal{F}_2$, taken now in the set enumerated by $\pi_{r_2,s_2}^{(\epsilon_2)}(m_2)$. 
The sequence $\mathcal{F}_2$ has the same properties,
mutatis mutandis, as the sequence $\mathcal{F}_1$, with now $2\ell_2+1-\epsilon_2$ petal trees, 
and a total of $r_2$ marked vertices and $s_2$ petals. Again $\mathcal{F}_2$ is transformed into a connected object by
adding a linear segment $S_2$ of $2\ell_2-\epsilon_2$ edges connecting its successive petal tree roots. 
To obtain the desired necklace, we will simply, by some ``partial gluing'' process reminiscent of that for
$(a,b)$-forests, squeeze the largest of the two 
segments $S_1$ and $S_2$ so as to get two segments of the same length $\min(2\ell_1-\epsilon_1,2\ell_2-\epsilon_2)$ 
which can then be glued together, head to tail, into a map with a spine that, after closing by some additional edge, forms 
the cycle separating the two distinguished faces of the necklace. Let us discuss in detail how we perform the partial and
mutual gluing processes. An important property of both processes is that we always identify a markable vertex
to an unmarked one. By transferring the possible marking of their markable copy, this guarantees that all the vertices obtained 
by gluing are markable and marked or not without ambiguity (and without double markings).  

Consider first the case where, say $\epsilon_1=-1$ and $\epsilon_2=+1$
so that $S_2$ has length $2\ell_2-1$, with some $\ell_2\geq 1$ and
$S_1$ has length $2\ell_1+1$ for some $\ell_1\geq 0$. Assume first
that $\ell_2\leq \ell_1$: we may then squeeze $S_1$ by ``pulling in''
its $(\ell_1+1)$-th vertex, namely by gluing the $(\ell_1+1-i)$-th
(markable) vertex along $S_1$ to the $(\ell_1+1+i)$-th (unmarked) one
for $i=1,\ldots, \ell_1-\ell_2+1$. This results in a segment $S'_1$ of
the same length as $S_2$, with $\ell_2$ markable vertices followed by
$\ell_2$ unmarked ones, together with a branch (attached to the
$\ell_2$-th vertex of $S'_1$) of length $\ell_1-\ell_2+1$ with all its
vertices markable but the last one, corresponding to the former
$(\ell_1+1)$-th vertex along $S_1$. This vertex was unmarked and we
decide to mark and distinguish it, say with the label $1$. We finally
glue the two segments $S'_1$ and $S_2$, head to tail, into a linear
spine that we close with an additional edge into a cycle of length
$2\ell_2$ separating two faces $f_1$ and $f_2$ -- see
Figure~\ref{fig:twofaces-1101}-top.  Again every markable vertex of
$S_2$ is matched to an unmarked vertex of $S'_1$ and vice versa, so
that all the vertices along the cycle are markable. It is easily
checked that $f_1$ and $f_2$ have respective degrees $2m_1$ and $2m_2$
and we end up with a petal necklace \emph{with an additional marked
  vertex} incident to $f_1$ but not to $f_2$.  For
$\ell_2\geq \ell_1+2$, a similar construction consisting now in
squeezing the segment $S_2$ by pulling in its $(\ell_2+1)$-th vertex
generates a necklace where the separating cycle has length
$2(\ell_1+1)$ and where the additional marked vertex $1$ is incident
to $f_2$ but not to $f_1$.  Finally, if $\ell_2=\ell_1+1$, $S_1$ and
$S_2$ have the same length and no squeezing is necessary: the faces
$f_1$ and $f_2$ are then separated by a cycle of length
$2\ell_2=2(\ell_1+1)$ and the vertex $1$ (obtained by gluing the
$(\ell_1+1)$-th vertex of $S_1$ to the $(\ell_2+1)$-th vertex of
$S_2$) lies along this cycle, i.e.~is incident to both faces.  The
construction is clearly reversible by a cutting procedure along the
separating cycle and along the branch leading from this cycle to
vertex $1$. We deduce that
$\pi_{r_1,s_1}^{(-1)}(m_1)\pi_{r_2,s_2}^{(1)}(m_2)$ counts petal
necklaces with $s_1$ unlabeled petals and $r_1$ unlabeled marked
vertices in its distinguished face $f_1$ of degree $2m_1$, $s_2$
unlabeled petals and $r_2$ unlabeled marked vertices in its
distinguished face $f_2$ of degree $2m_2$, and with an additional
marked vertex (distinct from all the others) distinguished as vertex
$1$ and lying anywhere in the map (the map hence has a total of
$r_1+r_2+1\geq 1$ marked vertices).  The necklaces are easily seen to
be tight, since the petal tree components have no unmarked leaves and
no unmarked leaf was created in the process. Still, as apparent from
the above discussion, a restriction applies to these necklaces since,
by construction, the length of the cycle separating their
distinguished faces is always even. It is easily seen that the missing
set of necklaces is enumerated by
$\pi_{r_1,s_1}^{(0)}(m_1)\pi_{r_2,s_2}^{(0)}(m_2)$: indeed, repeating
the above construction on the corresponding pairs of petal tree
sequences $\mathcal{F}_1$ (with $2\ell_1+1$ petal trees) and
$\mathcal{F}_2$ (with $2\ell_2+1$ petal trees), we get necklaces whose
separating cycle now has the odd length $2\min(\ell_1,\ell_2)+1$ (for
some arbitrary $\ell_1,\ell_2\geq 0$), again with an additional marked
vertex distinguished as vertex $1$. We finally obtain~\eqref{eq:pirs1}
by summing over $r_1$ and $r_2$ with $r_1+r_2=r_0$, and over $s_1$ and
$s_2$ with $s_1+s_2=s_0$, to account for the dispatching of the marked
vertices and petals among $f_1$ and $f_2$.

An alternative construction consists in starting from $\epsilon_1=0$ and $\epsilon_2=+1$ so that $S_2$ still has length $2\ell_2-1$, with some $\ell_2\geq 1$
while $S_1$ has length $2\ell_1$ for some $\ell_1\geq 0$. As before, we assume first that $\ell_2\leq \ell_1$: we may then squeeze $S_1$, now by ``pulling in'' its $\ell_1$-th connecting edge, namely by gluing the 
$(\ell_1+1-i)$-th (markable) vertex along $S_1$ to the $(\ell_1+i)$-th
(unmarked) one for $i=1,\ldots, \ell_1-\ell_2+1$. This results in a segment $S'_1$ of the same length as $S_2$, with $\ell_2$ markable vertices followed by $\ell_2$
unmarked ones, together with a branch (attached to the $\ell_2$-th vertex of $S'_1$) of length $\ell_1-\ell_2$ with all its vertices markable, having now a petal at its end corresponding to the former $\ell_1$-th edge along $S_1$. This newly created 
petal is adjacent to $f_1$  but not to $f_2$ and we decide to distinguish it with the label $1$. Gluing $S'_1$ and $S_2$ and closing, we end up with a cycle of even length 
$2\ell_2$ separating two faces $f_1$ and $f_2$, see the bottom of Figure~\ref{fig:twofaces-1101}.
For  $\ell_2\geq \ell_1+1$, squeezing now the segment $S_2$ by pulling in its $\ell_2$-th edge generates 
a necklace with a separating cycle
of odd length $2\ell_1+1$ and where the additional petal is now adjacent to $f_2$. To get the missing necklaces,
namely those with an odd separating cycle and a distinguished petal adjacent to $f_1$, and those with an even 
separating cycle and a distinguished petal adjacent to $f_2$, we have, as clear by symmetry, to supplement
the above family of necklaces enumerated by $\pi_{r_1,s_1}^{(0)}(m_1)\pi_{r_2,s_2}^{(1)}(m_2)$ by that enumerated by
$\pi_{r_1,s_1}^{(1)}(m_1)\pi_{r_2,s_2}^{(0)}(m_2)$.
We finally obtain~\eqref{eq:pirs2}
by summing over $r_1$ and $r_2$ with $r_1+r_2=r_0$, and over $s_1$ and
$s_2$ with $s_1+s_2=s_0$, to account for the dispatching of the marked
vertices and petals among $f_1$ and $f_2$.
\end{proof}
\begin{remark}
The consistency of the two expressions \eqref{eq:Pirs1} and \eqref{eq:Pirs2} for $\Pi_{r,s}(m_1,m_2)$ when
$r,s\geq 1$ may be checked directly by setting $s!=(s-1)!(s_1+s_2)$ in \eqref{eq:Pirs1}, $r!=(r-1)!(r_1+r_2)$ 
in \eqref{eq:Pirs2}, and upon using the identifications 
\begingroup
\allowdisplaybreaks
\begin{align}
 \sum_{\substack{r_1,r_2,s_1,s_2\geq
        0 \\ r_1+r_2=r-1 \\ s_1+s_2=s}}
    s_1  \pi_{r_1,s_1}^{(-1)}(m_1)
      \pi_{r_2,s_2}^{(1)}(m_2)&=  \sum_{\substack{r_1,r_2,s_1,s_2\geq
        0 \\ r_1+r_2=r \\ s_1+s_2=s-1}}
    r_1  \pi_{r_1,s_1}^{(0)}(m_1)
      \pi_{r_2,s_2}^{(1)}(m_2)\\
       \sum_{\substack{r_1,r_2,s_1,s_2\geq
        0 \\ r_1+r_2=r-1 \\ s_1+s_2=s}}
    s_2  \pi_{r_1,s_1}^{(-1)}(m_1)
      \pi_{r_2,s_2}^{(1)}(m_2)&\overset{(*)}{=}  \sum_{\substack{r_1,r_2,s_1,s_2\geq
        0 \\ r_1+r_2=r-1 \\ s_1+s_2=s}}
    s_2  \pi_{r_1,s_1}^{(1)}(m_1)
      \pi_{r_2,s_2}^{(-1)}(m_2)\\
      &=  \sum_{\substack{r_1,r_2,s_1,s_2\geq
        0 \\ r_1+r_2=r \\ s_1+s_2=s-1}}
    r_2  \pi_{r_1,s_1}^{(1)}(m_1)
      \pi_{r_2,s_2}^{(0)}(m_2)\\
       \sum_{\substack{r_1,r_2,s_1,s_2\geq
        0 \\ r_1+r_2=r-1 \\ s_1+s_2=s}}
    s_1  \pi_{r_1,s_1}^{(0)}(m_1)
      \pi_{r_2,s_2}^{(0)}(m_2)&=  \sum_{\substack{r_1,r_2,s_1,s_2\geq
        0 \\ r_1+r_2=r \\ s_1+s_2=s-1}}
    r_1  \pi_{r_1,s_1}^{(1)}(m_1)
      \pi_{r_2,s_2}^{(0)}(m_2)\\
       \sum_{\substack{r_1,r_2,s_1,s_2\geq
        0 \\ r_1+r_2=r-1 \\ s_1+s_2=s}}
    s_2  \pi_{r_1,s_1}^{(0)}(m_1)
      \pi_{r_2,s_2}^{(0)}(m_2)&=  \sum_{\substack{r_1,r_2,s_1,s_2\geq
        0 \\ r_1+r_2=r \\ s_1+s_2=s-1}}
    r_2  \pi_{r_1,s_1}^{(0)}(m_1)
      \pi_{r_2,s_2}^{(1)}(m_2)\ .
\end{align}
\endgroup
Here all equalities follow from the identity $s  \pi_{r,s}^{(\epsilon)}(m)=(r+1)\pi_{r+1,s-1}^{(1+\epsilon)}(m)$ for $s\geq 1$,
apart from that marked with a $(*)$ which follows from the identity  $\pi_{r,s}^{(1)}(m)=\pi_{r,s}^{(-1)}(m)
+\pi_{r-1,s}^{(-1)}(m)$ for any $r\geq 0$ (with the convention
$\pi_{-1,s}^{(-1)}(m)=0$), see \eqref{eq:trans21}-\eqref{eq:trans24} in Appendix~\ref{app:symmetry}
for similar so-called transmutation relations.
\end{remark}

\subsection{Tight slices}
\label{sec:generalslices}

As in the case of bipartite maps, a key ingredient in the derivation of a general
formula for the number of tight general non-bipartite maps is the enumeration
of tight slices. Recall the basic definitions pertaining to (tight) slices from Section
\ref{sec:pointedrooted}, and note that we do not assume anymore that the face
degrees be even integers. In this general context, it is still
true that the only elementary slice of tilt $-1$ is the trivial slice,
but a major difference is that the set of elementary slices of tilt
$0$ is now non-empty.

As discussed in Section \ref{sec:pointedrooted}, a (tight) slice $\mathbf{s}$ with width $w\geq 1$ and tilt $t$
can be decomposed into a collection of $w$ (tight) elementary
slices. This discussion applies without change in our general context: we cut the slice
along the leftmost geodesics started from the consecutive corners $c_0,c_1,\ldots,c_w$
incident to the base, and record the lattice path
$L=(\ell_0,\ell_1,\ldots,\ell_w)$ where $\ell_i=d(c_0,A)-d(c_i,A)$,
and $A$ is the apex of $\mathbf{s}$. The only difference is that the tilt $t_i=\ell_i-\ell_{i-1}$ of
the slice $\mathbf{s}'_i$ with base $c_{i-1}c_i$ can take any value $\{-1,0,1\}$ rather than only
$\{-1,1\}$. As before, the tightness condition requires that, for
$i\in \{2,\ldots,w\}$, if 
$t_{i}=+1$ and $t_{i-1}=-1$, then the slice $\mathbf{s}'_{i}$ is either the marked empty slice, or
a tight elementary slice with at least one inner face. For $i\in
\{1,\ldots,w\}$, if $t_i=-1$
then the slice $\mathbf{s}'_i$ is necessarily trivial, while if $t_i=0$ then
the slice $\mathbf{s}'_i$ automatically has at least one inner face of
odd degree. We may
forget the redundant information of all down steps, as well as up steps decorated with unmarked
empty slices, by letting $\mathbf{s}_j,1\leq j\leq k$ be the sequence
$(\mathbf{s}'_i,1\leq i\leq w)$ whose trivial and unmarked empty elements
have been removed. This gives:

\begin{proposition}
  \label{sec:gener-non-bipart-1}
  There is a one-to-one correspondence between tight slices of width
  $w$ and tilt $t$ on the one hand, and
  pairs of the form $(L,(\mathbf{s}_j,1\leq j\leq k))$ on the other
  hand, where: 
   \begin{itemize}
   \item $L=(\ell_0,\ldots,\ell_w)$ is a 
   lattice path from $(0,0)$ to
   $(w,t)$ with increments in $\{-1,0,1\}$, that has $k$ marked 
   steps $1\leq i_1<i_2<\cdots<i_k\leq w$ with
   $\ell_{i_j}-\ell_{i_j-1}\in \{0,1\}$, in such a way that every up step immediately
   following a down step is marked, and every horizontal step is marked,  
 \item for $1\leq j\leq k$, $\mathbf{s}_j$ is a tight elementary slice
   of tilt $\ell_{i_j}-\ell_{i_j-1}$, and in the case where this tilt
   is $1$, $\mathbf{s}_j$ is either the marked empty slice, or has at least
   one inner face.  
\end{itemize}
\end{proposition}

We obtain the following generalization of Corollary
\ref{sec:case-pointed-rooted-1}.

\begin{corollary}
  \label{sec:gener-non-bipart}
  For $m \in \Z_{>0}/2$, there is a one-to-one correspondence between
  non-empty tight elementary slices of tilt $\epsilon\in \{0,1\}$,
  whose inner face incident to the base edge has degree $2m$ on the
  one hand, and pairs of the form $(L,(\mathbf{s}_j,1\leq j\leq k))$
  on the other hand, where $k$ is a non-negative integer, and: 
   \begin{itemize}
   \item $L=(\ell_0,\ldots,\ell_{2m-1})$ is a 
   lattice path from $(0,0)$ to
   $(2m-1,\epsilon)$ with increments in $\{-1,0,1\}$, that has $k$ marked 
   steps $1\leq i_1<i_2<\cdots<i_k\leq 2m-1$ with
   $\ell_{i_j}-\ell_{i_j-1}\in \{0,1\}$, in such a way that every up step immediately
   following a down step is marked, and every horizontal step is
   marked, 
 \item for $1\leq j\leq k$, $\mathbf{s}_j$ is a tight elementary slice
   of tilt $\ell_{i_j}-\ell_{i_j-1}$, and in the case where this tilt
   is $1$, $\mathbf{s}_j$ is either the marked empty slice, or has at least
   one inner face.  
\end{itemize}
\end{corollary}

By iterating the decomposition of this corollary, we obtain an
encoding of tight elementary slices by plane trees, where every vertex
may be associated with a slice of tilt in $\{0,1\}$. Let
$m_1,\ldots,m_n$ be integers or half integers, not all equal to $0$,
and $\epsilon\in \{0,1\}$. Let $\mathcal{T}_\epsilon(m_1,\ldots,m_n)$ be the
family of pairs $(\mathbf{t},(L_i)_{1\leq i\leq n})$, satisfying the
following conditions. 
\begin{itemize}
\item $\mathbf{t}$ is a rooted plane tree with $n$ vertices
  labeled by $\{1,2,\ldots,n\}$. 
  \item Each vertex receives a type
  in $\{0,1\}$, and the root has type $\epsilon$. We let 
  $\epsilon_i\in \{0,1\}$ be the type of vertex $i$ and 
  $(w^{(i)}_1,\ldots,w^{(i)}_{k_i})$ be the types of the $k_i$
  children of vertex $i$ in $\mathbf{t}$ in planar order. 
\item For $i\in \{1,2,\ldots,n\}$ such that 
   $m_i>0$, 
  $L_i=(\ell^{(i)}_0,\ldots,\ell^{(i)}_{2m_i-1})$ is a lattice path from $(0,0)$ to $(2m_i-1,\epsilon_i)$ with
  increments in $\{-1,0,1\}$, and with  $k_i$ marked steps $1\leq j_1<\cdots
  <j_{k_i}\leq 2m_i-1$, such that 
  $\ell^{(i)}_{j_r}-\ell^{(i)}_{j_r-1}=w^{(i)}_r$ for every $r\in
  \{1,\ldots,k_i\}$. In particular, only up or horizontal steps may be
  marked. We further require that 
   every up step of $L_i$ immediately following
  a down step is marked, and every horizontal step of $L_i$ is marked.
\item  For $i\in \{1,2,\ldots,n\}$ such that  $m_i=0$, we have
  $\epsilon_i=1$, $k_i=0$  and 
  $L_i$ is the trivial path $\{(0,0)\}$.
\end{itemize}

This leads to the following generalization of Corollary
\ref{sec:decomp-slice-into}. 

\begin{corollary}
  \label{sec:gener-non-bipart-2}
  For $\epsilon\in \{0,1\}$ and every choice of integers  or
  half-integers $m_1,\ldots,m_n$ not all
  equal to $0$, the iterated slice decomposition yields a bijection
  between the set of elementary tight slices with $n$ labeled
  boundaries of respective lengths $2m_1,\ldots,2m_n$
  with tilt $\epsilon$ and
the set  $\mathcal{T}_{\epsilon}(m_1,\ldots,m_n)$. 
\end{corollary}

It turns out that, in order to prove Theorem~\ref{thm:quasimain} in
Section~\ref{sec:generalnonbip}, we
will need an extension of this corollary associating sequences of
slices with sequences of trees, namely forests. Still, as a warm-up,
let us first consider the enumeration of 
$\mathcal{T}_{\epsilon}(m_1,\ldots,m_n)$, which is done based on 
that of {\em two-type trees}, obtained as a special case of Proposition
\ref{prop:twotypeforest} of Appendix
\ref{sec:enumeration-two-type}. Note that the types $1$ and $0$ considered in
the present section correspond to types $A$ and $B$ respectively in
the notation of 
the appendix. 
We count the elements $(\mathbf{t},(L_i)_{1\leq i\leq n})$ of
$\mathcal{T}_{\epsilon}(m_1,\ldots,m_n)$ by 
fixing the types $\epsilon_1,\ldots,\epsilon_n$ of the vertices
$1,\ldots,n$ of $\mathbf{t}$ (so that $\epsilon_i=w^{(i)}_0$ in the notation
of Appendix \ref{sec:enumeration-two-type}) with the constraint that
the root vertex should have type $\epsilon$, 
as well as the sequence of types
$w^{(i)}_1,\ldots,w^{(i)}_{k_i}$ of the consecutive children (in
planar order) of vertex $i$. These sequences must satisfy the 
consistency conditions \eqref{eq:Aconsis} and \eqref{eq:Bconsis}
of Appendix \ref{sec:enumeration-two-type}, which in this context are rewritten as 
\begin{equation}
  \label{eq:2}
  \sum_{i=1}^n\epsilon_i=\epsilon+\sum_{i=1}^nr_i \qquad \mbox{ and }\quad\sum_{i=1}^n\bar{\epsilon}_i=\bar{\epsilon}+\sum_{i=1}^ns_i\, ,
\end{equation}
where we let $\bar{\eta}=1-\eta$ for $\eta\in \{0,1\}$, and where
$r_i=\sum_{j=1}^{k_i}w^{(i)}_j$ and
$s_i=\sum_{j=1}^{k_i}\bar{w}^{(i)}_j$ are the numbers of type-$1$ and
type-$0$ children of
vertex $i$. These consistency
relations simply express in two different ways the number of
type $1$ (resp.\ $0$) vertices in the tree $\mathbf{t}$. Note that the number
of type-$0$ (resp.\ type-$1$) vertices that are the children of
type-$1$ (resp.\ type-$0$) vertices (those are respectively the
numbers $b^A$ and $a^B$ considered in Appendix \ref{sec:enumeration-two-type}) is
given by
\begin{equation}
  \label{eq:4}
  \sum_{i=1}^n\epsilon_i s_i\qquad
  \left(\mbox{resp. }\ \sum_{i=1}^n\bar{\epsilon}_i r_i\right)\, .
\end{equation}
Finally, by Proposition~\ref{prop:twotypeforest} (in the
notation therein, we have $a^O=\epsilon$ and $b^O=\bar{\epsilon}$,
expressing the fact that we are counting ``forests'' made of only one
tree with root of type $\epsilon$), we see that the
number of possible trees $\mathbf{t}$ contributing to
$\mathcal{T}_\epsilon(m_1,\ldots,m_n)$, with types given by
$\epsilon_1,\ldots,\epsilon_n$ and given a consistent type array $(w^{(i)}_j)_{j=1,\ldots,
 k_i}^{i=1\ldots, n}$, is equal to 
\begin{equation}
  \label{eq:5}
  \left(\epsilon
    \sum_{i=1}^n\epsilon_is_i+\bar{\epsilon}\sum_{i=1}^n\bar{\epsilon}_ir_i\right)\left(\sum_{i=1}^n\epsilon_i-1\right)!\left(\sum_{i=1}^n\bar{\epsilon}_i-1\right)!
\end{equation}
where $\sum_{i=1}^n\epsilon_i$ and $\sum_{i=1}^n\bar{\epsilon}_i$ are the numbers of
vertices of types $1$ and $0$ respectively. In the case where either
one of these numbers is equal to $0$, in accordance with
\eqref{eq:twotypeforestconst}, we should replace the whole
formula by $(n-2)!$.

It now remains to enumerate, for a given consistent array 
$(w^{(i)}_j)^{i=1,\ldots,n}_{j=1,\ldots,k_i}$, and for a given tree
$\mathbf{t}$ as above, the number of possible
lattice paths $L_i,1\leq i\leq n$ so that $(\mathbf{t},(L_i)_{1\leq
  i\leq n})$ is an element of
$\mathcal{T}_{\epsilon}(m_1,\ldots,m_n)$.  As before, we let 
$r_i=\sum_{j=1}^{k_i}w^{(i)}_j$ and
$s_i=\sum_{j=1}^{k_i}\bar{w}^{(i)}_j=k_i-r_i$. 
If $m_i=0$ then by definition of
$\mathcal{T}_\epsilon(m_1,\ldots,m_n)$, this requires 
$L_i=\{(0,0)\}$, $r_i=s_i=0$ and $\epsilon_i=1$: by 
\eqref{eq:piepsm0}, this is precisely counted by
$\pi^{(\epsilon_i)}_{r_i,s_i}(0)=\delta_{r_i, 0}\delta_{s_i,
  0}\delta_{\epsilon_i, 1}$. 
If $m_i>0$, $L_i$ should be
a lattice path from $(0,0)$ to $(2m_i-1,\epsilon_i)$, with increments in
$\{-1,0,1\}$, in which $k_i$ up or horizontal steps are marked, say $1\leq j_1<\cdots<j_{k_i}\leq
2m_i-1$, with increments at these steps respectively  given by
$w^{(i)}_1,\ldots,w^{(i)}_{k_i}$, and in such a way that all horizontal steps
are marked, as well as all the up steps immediately following a down
step. As explained in the proof of Proposition
\ref{prop:pieps} (and up to changing the up steps into down steps and
vice versa), the number of such paths is equal to 
$p_{r_i+s_i,s_i+1+\epsilon_i}(m_i)$ if
$m_i-\frac{s_i+1+\epsilon_i}{2}\in \Z$, and
to zero otherwise. In particular, this
number depends on the array
$(w^{(i)}_j)^{i=1,\ldots,n}_{j=1,\ldots,k_i}$ only through the values
of $r_i,s_i$.
Noting that for a given
value of $r_i,s_i$, there are $\binom{r_i+s_i}{s_i}$ possible choices
of $(w^{(i)}_j,1\leq j\leq k_i)$ inducing these values, and recalling
that
$\pi^{(\epsilon_i)}_{r_i,s_i}=\binom{r_i+s_i}{s_i}p_{r_i+s_i,s_i+1+\epsilon_i}(m_i)$,
we finally obtain the following result. 

\begin{proposition}
  \label{sec:gener-non-bipart-3}
For $\epsilon=0,1$ and for every choice of integers or
half-integers $m_1,\ldots,\allowbreak m_n$ not all equal to $0$,
the cardinality of $\mathcal{T}_{\epsilon}(m_1,\ldots,m_n)$, and hence
the number of elementary tight slices with $n$ labeled boundaries of
respective lengths $2m_1,\ldots,2m_n$ and tilt $\epsilon$, is given
by
\begin{multline}
  \label{eq:6}
\mathrm{Card}\left(\mathcal{T}_\epsilon(m_1,\ldots,m_n)\right)= \\
  \sum_{\substack{\epsilon_1,\ldots,\epsilon_n=0,1\\r_1,\ldots,r_n,s_1,\ldots,s_n\geq 0\\
    \epsilon_1+\cdots+\epsilon_n=\epsilon+r_1+\cdots+r_n\\
    \bar{\epsilon}_1+\cdots+\bar{\epsilon}_n=\bar{\epsilon}+s_1+\cdots+s_n}}\!\!\!\!\!\!
  \left(\epsilon
  \sum_{i=1}^n\epsilon_is_i+\bar{\epsilon}\sum_{i=1}^n\bar{\epsilon}_ir_i\right)
\left(\sum_{i=1}^n\epsilon_i-1\right)!\left(\sum_{i=1}^n\bar{\epsilon}_i-1\right)!\prod_{i=1}^n\pi^{(\epsilon_i)}_{r_i,s_i}(m_i)\, ,
\end{multline}
where the first three
factors in the sum should be replaced by $(n-2)!$ whenever the
argument in either of the two factorials equals $-1$. 
\end{proposition}

This proposition has interesting consequences for certain evaluations
of lattice count polynomials. Proceeding similarly to Section
\ref{sec:pointedrooted}, we may indeed associate bijectively tight
maps with  elementary tight slices of tilt in $\{0,1\}$. Precisely,
starting from an elementary tight slice of tilt $0$, we may glue
isometrically the left and right boundaries to obtain a tight map with
one marked petal (delimited by the base edge of the slice) and one
marked vertex (given by the apex). This construction can be inverted
by cutting along the leftmost geodesic from the marked petal to the
pointed vertex. Similarly, starting from an elementary tight slice of tilt $1$, we may glue
isometrically the left and right boundaries \emph{starting from the
  base}\footnote{This means in particular that the two
endpoints of the base are glued together. If we proceed as in Section \ref{sec:pointedrooted} and glue the left and right boundaries of a slice of tilt 1 starting from the apex, we do obtain a pointed rooted tight map, but this is not the most general such map, since by construction the two ends of the root edge (inherited from the base) are at different distances from the distinguished vertex (inherited from the apex).} to obtain a tight map with two marked
petals (one being delimited by the base edge of the slice, and the
other one by the edge of the left boundary incident to the apex). The
inverse construction is more involved, but is in fact a particular
case of the slice decomposition of annular maps.

These bijections imply that for every
$m_1,\ldots,m_n\in \Z_{\geq 0}/2$ not all equal to zero, and for
$\epsilon\in \{0,1\}$, we have 
  \begin{equation}
      \label{eq:9}
    N_{0,n}(2m_1,\ldots,2m_n,1,\epsilon)=\mathrm{Card}(\mathcal{T}_\epsilon(m_1,\ldots,m_n))\,
    . 
  \end{equation}

Finally, from the discussion of \cite[Appendix A]{hankel}, we
obtain the identity
\begin{multline}
  \label{eq:10}
  N_{0,n}(2m_1,\ldots,2m_n,2,0)=\mathrm{Card}\left(\mathcal{T}_1(m_1,\ldots,m_n)\right)
  \\
  +\sum_{\{1\}\subset I\subsetneq
  \{1,\ldots,n\}}\mathrm{Card}\left(\mathcal{T}_0\left((m_i)_{i\in I}\right)\right)
\mathrm{Card}\left(\mathcal{T}_0\left((m_i)_{i\notin
      I}\right)\right)\, .
\end{multline}

\subsection{General maps}
\label{sec:generalnonbip}

We are now ready to get a general enumeration formula for planar
tight maps.
Our approach will be parallel to that of Section~\ref{sec:generalbip},
and we start by stating the following analog of
Proposition~\ref{prop:generalbij}:
 \begin{proposition}
  \label{prop:generalbijnonbip}
  Let $m_1,\ldots,m_n$ be non-negative integers or half-integers
  ($n \geq 3$) with $m_1,m_2>0$. Then, there is a bijection between
  the set of planar tight maps with $n$ boundaries labeled from $1$ to
  $n$ with respective lengths $2m_1,2m_2,\ldots,2m_n$, and the set of
  pairs $(\mathbf{m}_{12},\mathbf{s})$ such that there exist
  $r_0,s_0\geq 0$ with $r_0+s_0 \in \{0,\ldots,n-3\}$ for which the
  following holds.
  \begin{itemize}
  \item $\mathbf{m}_{12}$ is a tight petal necklace with two distinguished
    faces of degrees $2m_1,2m_2$, with one extra distinguished element being either 
    a marked vertex or a petal, and with $r_0$ other marked vertices and $s_0$ other petals.
  \item $\mathbf{s}=(\mathbf{s}_1,\ldots,\mathbf{s}_{r_0+s_0+1})$ is a
    $(r_0+s_0+1)$-tuple of slices such that:
    \begin{itemize}
    \item each $\mathbf{s}_j$ ($j=1,\ldots,r_0+s_0+1$) is a tight
      elementary slice of tilt $0$ or $1$ containing at least one
      inner face or marked vertex,
    \item the inner faces and marked vertices of these slices are
      labeled by integers in $\{3,\ldots,n\}$,
    \item each $i \in \{3,\ldots,n\}$ appears in exactly one
      $\mathbf{s}_j$ and labels an inner face of degree $2m_i$ for
      $m_i>0$, or a marked vertex for $m_i=0$,
  \item the label $3$ appears in the first slice $\mathbf{s}_1$.
    \end{itemize}
  \item $\mathbf{m}_{12}$ and $\mathbf{s}$ are compatible in the sense
    that, if we label the marked vertices and petals of
    $\mathbf{m}_{12}$ in a canonical way from $1$ to $r_0+s_0+1$ as in
    Remark~\ref{rem:deuxfacesnonbip-lab}, and for
    $j=1,\ldots,r_0+s_0+1$ we set $w^{(0)}_j=1$ (resp.\ $w^{(0)}_j=0$)
    if label $j$ is on a marked vertex (resp.\ petal), then
    $\mathbf{s}_j$ has tilt $w^{(0)}_j$. Note that
    $r_0=\sum_{j=2}^{r_0+s_0+1}w^{(0)}_j$.
  \end{itemize}
\end{proposition}

The proof of this proposition follows exactly the same lines as that
of Proposition~\ref{prop:generalbij} in
Section~\ref{sec:generalbip}. In particular,
Lemma~\ref{sec:decomp-an-annul} admits a direct non-bipartite
extension in which the maps $\mathbf{m}$ may have all their faces of
arbitrary degrees, the lattice paths $L,L'$ may contain horizontal
steps, and consistently the sequences
$\boldsymbol\sigma,\boldsymbol\sigma'$ may contain slices of tilt
$0$. Further details are left to the reader.  We deduce the following
enumerative result:

\begin{proposition}
  \label{prop:gennonbip}
  For $n \geq 3$ and for non-negative integers or half-integers
  $m_1,\ldots,m_n$ with, say, $m_1,m_2>0$, we have
  \begingroup
  \allowdisplaybreaks
  \begin{align}
    \label{eq:gennonbip}
    N_{0,n}(2m_1,\ldots,2m_n) =& \nonumber \\
    \sum_{\substack{\epsilon_3,\ldots,\epsilon_n=0,1 \\ r_1,\ldots,r_n,s_1,\ldots,s_n\geq 0 \\
        \sum_{i=3}^n \epsilon_i=\sum_{i=1}^n r_i+1 \\
        \sum_{i=3}^n (1\!-\!\epsilon_i)=\sum_{i=1}^n s_i}} 
    &\left( \sum_{i=1}^n r_i \right)! \left(\epsilon_3 (s_1+s_2) + \sum_{i=3}^n \epsilon_i s_i \right) \left( \sum_{i=1}^n s_i -1 \right)!  \\ \nonumber
    \times &\left( \pi_{r_1,s_1}^{(-1)}(m_1)
      \pi_{r_2,s_2}^{(1)}(m_2) + \pi_{r_1,s_1}^{(0)}(m_1)
      \pi_{r_2,s_2}^{(0)}(m_2) \right) \prod_{i=3}^n \pi_{r_i,s_i}^{(\epsilon_i)}(m_i) \\ \nonumber
    + \!\!\!\!\! 
    \sum_{\substack{\epsilon_3,\ldots,\epsilon_n=0,1 \\ r_1,\ldots,r_n,s_1,\ldots,s_n\geq 0 \\
        \sum_{i=3}^n \epsilon_i=\sum_{i=1}^n r_i \\
        \sum_{i=3}^n (1\!-\!\epsilon_i)=\sum_{i=1}^n s_i+1}} \!\!\!\!\!
    &\left( \sum_{i=1}^n s_i \right)!
    \left((1-\epsilon_3) (r_1+r_2) + \sum_{i=3}^n (1-\epsilon_i) r_i \right)
    \left( \sum_{i=1}^n r_i - 1\right)! \\ \nonumber
    \times  &\left( \pi_{r_1,s_1}^{(0)}(m_1)
      \pi_{r_2,s_2}^{(1)}(m_2) + \pi_{r_1,s_1}^{(1)}(m_1)
      \pi_{r_2,s_2}^{(0)}(m_2) \right) \prod_{i=3}^n \pi_{r_i,s_i}^{(\epsilon_i)}(m_i)
  \end{align}
  \endgroup where $\pi_{r,s}^{(\epsilon)}$ is the univariate
  quasi-polynomial defined in~\eqref{eq:pieps} and where it is
  understood that
  $(\epsilon_3 (s_1+s_2) + \sum_{i=3}^n \epsilon_i s_i) ( \sum_{i=1}^n
  s_i -1)!$ is equal to $1$ if all the $s_i$ are zero, and that
  $((1-\epsilon_3) (r_1+r_2) + \sum_{i=3}^n (1-\epsilon_i) r_i ) (
  \sum_{i=1}^n r_i - 1)!$ is equal to $1$ if all the $r_i$ are zero.
\end{proposition}

\begin{proof}
  We need to enumerate the compatible pairs
  $(\mathbf{m}_{12},\mathbf{s})$ of
  Proposition~\ref{prop:generalbijnonbip}.  Note first that, given
  $r_0$ and $s_0$, the number of possible petal necklaces
  $\mathbf{m}_{12}$ is given by
  Proposition~\ref{prop:deuxfacesnonbip-unlab}.
  Turning now to the number of possible $\mathbf{s}$ compatible with a
  given petal necklace $\mathbf{m}_{12}$, it is given by a direct extension of
  Corollary~\ref{sec:gener-non-bipart-2} as follows: by decomposing
  recursively each elementary tight slice $\mathbf{s}_j$ (with tilt
  $w^{(0)}_j$) into a tree of lattice paths, we see that the set of
  possible $(r_0+s_0+1)$-tuples $\mathbf{s}$ is in bijection with the
  set
  $\mathcal{F}_{(w^{(0)}_1,\ldots,w^{(0)}_{r_0+s_0+1})}(m_3,\ldots,m_n)$
  defined as the set of pairs $(\mathbf{f},(L_i)_{1\leq i\leq n})$
  satisfying the following conditions.

\begin{itemize}
\item $\mathbf{f}$ is a plane forest with $r_0+s_0+1$ connected
  components, i.e.\ a $(r_0+s_0+1)$-tuple of rooted plane trees, and with a
  total of $n-2$ vertices labeled 
  by $\{3,\ldots,n\}$, the label $3$ appearing in
  the first component.
\item Each vertex receives a type in $\{0,1\}$, and the root vertex of
  the $j$-th tree component in the forest has type $w^{(0)}_j$.  For
  $i\in \{3,\ldots,n\}$, we let $\epsilon_i\in \{0,1\}$ be the type of
  vertex $i$, $k_i$ be its number of children, and
  $(w^{(i)}_1,\ldots,w^{(i)}_{k_i})$ be the types of these children
  (numbered in planar order in the rooted tree component at hand).  We
  also set $r_i:=\sum_{j=1}^{k_i}w^{(i)}_j$ the numbers of those
  children which are of type $1$ and $s_i:=k_i-r_i$ the numbers of
  those children which are of type $0$.
\item For $i\in \{3,\ldots,n\}$:
  \begin{itemize}
  \item if $m_i>0$, $L_i=(\ell^{(i)}_0,\ldots,\ell^{(i)}_{2m_i-1})$ is
    a lattice path from $(0,0)$ to $(2m_i-1,\epsilon_i)$ with
    increments in $\{-1,0,1\}$, and with $k_i$ marked steps
    $1\leq j_1<\cdots <j_{k_i}\leq 2m_i-1$, such that
    $\ell^{(i)}_{j_r}-\ell^{(i)}_{j_r-1}=w^{(i)}_r$ for every
    $r\in \{1,\ldots,k_i\}$. In particular, only up or horizontal
    steps may be marked, and we further require that every horizontal
    step of $L_i$ is marked, as well as every up step immediately
    following a down step,
  \item if $m_i=0$, we have $\epsilon_i=1$, $k_i=r_i=s_i=0$, and $L_i$ is
    the trivial path $\{(0,0)\}$.
    \end{itemize}
\end{itemize}
We may now obtain the number of elements of $\mathcal{F}_{(w^{(0)}_1,\ldots,w^{(0)}_{r_0+s_0+1})}(m_3,\ldots,m_n)$ 
from the results of Appendix~\ref{sec:enumeration-two-type} for the enumeration of two-type forests,
where the types $A$ and $B$ therein correspond to types $1$ and $0$ in the present setting. 
For fixed $\epsilon_i$ and $k_i$, $i=3,\ldots,n$ and for fixed $(w^{(i)}_j)^{i=3,\ldots,n}_{j=1,\ldots,k_i}$, 
the number of forests $\mathbf{f}$
is non-zero only if the two consistency relations \eqref{eq:Aconsis} and \eqref{eq:Bconsis} are satisfied, namely
\begin{equation}
\label{eq:checkconsitency01}
\sum_{i=3}^n \epsilon_i=w^{(0)}_1+r_0+\sum_{i=3}^n r_i\quad\hbox{and}\quad 
\sum_{i=3}^n \bar{\epsilon}_i=\bar{w}^{(0)}_1+s_0+\sum_{i=3}^n s_i
\end{equation}
where,
as before, we use the shorthand notation $\bar{\eta}:=1-\eta$ for $\eta\in\{0,1\}$. Again, these identities simply express in two different
ways the number of vertices of type $1$ (first identity) and of type $0$ (second identity), corresponding 
respectively to the quantities denoted by $a$ and $b$ in Appendix~\ref{sec:enumeration-two-type}.
When these conditions are satisfied, we may use the constrained enumeration result \eqref{eq:twotypeforestconst},
with the correspondence $a^O=w^{(0)}_1+r_0$, $b^O=\bar{w}^{(0)}_1+s_0$, $a^B=\sum_{i=3}^n\bar{\epsilon}_i r_i$
and $b^A= \sum_{i=3}^n \epsilon_i s_i$:
if the vertex $3$ is of type $1$, i.e.\ $\epsilon_3=1$, the number of forests $\mathbf{f}$ is given by
\begin{equation}
      \label{eq:twotypeforestconst01}
      \begin{cases}
        \left(s_0 + \sum\limits_{i=3}^n \epsilon_i s_i\right)   \left(r_0+\sum\limits_{i=3}^n r_i\right)!   \left(s_0+\sum\limits_{i=3}^n s_i-1\right)! & \text{if $w_1^{(0)}=1$,}\\
         \left(\sum\limits_{i=3}^n\bar{\epsilon}_i r_i\right)   \left(r_0+\sum\limits_{i=3}^n r_i-1\right)!   \left(s_0+\sum\limits_{i=3}^n s_i\right)! & \text{if $w_1^{(0)}=0$,}\\
            \end{cases}
    \end{equation}
with, in the first line, the convention  $(s_0 + \sum_{i=3}^n \epsilon_i s_i)\times (s_0+\sum_{i=3}^n s_i-1)!=1$ if $s_0+\sum_{i=3}^n s_i=0$, i.e.\ when
$s_0$ and all the $s_i$'s for $i\geq3$ are zero. Note that, in the second line, the quantity $r_0+\sum_{i=3}^n r_i-1$ is always nonnegative
since, from the first consistency relation, it is equal to $\sum_{i=4}^n \epsilon_i$.
By symmetry, if the vertex $3$ is of type $0$, i.e.\ $\epsilon_3=0$, the number of forests $\mathbf{f}$ is given by
\begin{equation}
      \label{eq:twotypeforestconst01bis}
      \begin{cases}
          \left(r_0 + \sum\limits_{i=3}^n \bar{\epsilon}_i r_i\right)   \left(s_0+\sum\limits_{i=3}^n s_i\right)!   \left(r_0+\sum\limits_{i=3}^n r_i-1\right)! & \text{if $w_1^{(0)}=0$,}\\
          \left(\sum\limits_{i=3}^n\epsilon_i s_i\right)   \left(s_0+\sum\limits_{i=3}^n s_i-1\right)!   \left(r_0+\sum\limits_{i=3}^n r_i\right)! & \text{if $w_1^{(0)}=1$,}\\
            \end{cases}
    \end{equation}
with, in the first line, the convention  $(r_0 + \sum_{i=3}^n \bar{\epsilon}_i r_i)\times (r_0+\sum_{i=3}^n r_i-1)!=1$ if $r_0+\sum_{i=3}^n r_i=0$, i.e.\ when
$r_0$ and all the $r_i$'s for $i\geq3$ are zero. Again, in the second line, the quantity $s_0+\sum_{i=3}^n s_i-1$ is always nonnegative
since, from the second consistency relation, it is equal to $\sum_{i=4}^n \bar{\epsilon}_i$.
Both cases $\epsilon_3=1$ or $0$ above may be summarized into the enumeration formulas 
 \begin{equation}
      \label{eq:twotypeforestconst01ter}
      \begin{cases}
          \left(\epsilon_3 s_0 + \sum\limits_{i=3}^n \epsilon_i s_i\right)    \left(s_0+\sum\limits_{i=3}^n s_i-1\right)!   \left(r_0+\sum\limits_{i=3}^n r_i\right)!& \text{if $w_1^{(0)}=1$,}\\
          \left(\bar{\epsilon}_3 r_0+\sum\limits_{i=3}^n\bar{\epsilon}_i r_i\right)   \left(r_0+\sum\limits_{i=3}^n r_i-1\right)!   \left(s_0+\sum\limits_{i=3}^n s_i\right)! & \text{if $w_1^{(0)}=0$,}\\
            \end{cases}
    \end{equation}
with the conventions that $(\epsilon_3 s_0 + \sum_{i=3}^n \epsilon_i s_i) (s_0+\sum_{i=3}^n s_i-1)!=1$ if
$s_0$ and all the $s_i$'s for $i\geq3$ are zero and that  $(\bar{\epsilon}_3 r_0 + \sum_{i=3}^n \bar{\epsilon}_i r_i) (r_0+\sum_{i=3}^n r_i-1)!=1$ if
$r_0$ and all the $r_i$'s for $i\geq3$ are zero.

It remains to enumerate, for a given array 
$(w^{(i)}_j)^{i=3,\ldots,n}_{j=1,\ldots,k_i}$ satisfying the consistency relations \eqref{eq:checkconsitency01},
and for a given two-type forest $\mathbf{f}$ as above, the number of families of
lattice paths $(L_i)_{3\leq i\leq n}$ so that $(\mathbf{f},(L_i)_{3\leq i\leq n})$ is an element of
$\mathcal{F}_{(w^{(0)}_1,\ldots,w^{(0)}_{r_0+s_0+1})}(m_3,\ldots,m_n)$.  Repeating the counting argument in the proof of
Proposition~\ref{sec:gener-non-bipart-3}, this number is simply equal to $\prod_{i=3}^n\pi^{(\epsilon_i)}_{r_i,s_i}(m_i)$.
Combining with \eqref{eq:twotypeforestconst01ter}, and recalling the consistency relations
\eqref{eq:checkconsitency01}, we obtain in the case $w^{(0)}_1=1$
\begin{equation}
    \label{eq:Fw1}
    \begin{split}
    &\mathrm{Card}\left(\mathcal{F}_{(w^{(0)}_1=1,w^{(0)}_2,\ldots,w^{(0)}_{r_0+s_0+1})}(m_3,\ldots,m_n)\right)=\\
    &\quad \sum_{\substack{\epsilon_3,\ldots,\epsilon_n=0,1 \\ r_3,\ldots,r_n,s_3,\ldots,s_n\geq 0 \\
        \sum_{i=3}^n \epsilon_i=r_0+\sum_{i=3}^n r_i+1 \\
        \sum_{i=3}^n \bar{\epsilon}_i=s_0+\sum_{i=3}^n s_i}} 
    \left(r_0+ \sum_{i=3}^n r_i \right)! \left(\epsilon_3 s_0 + \sum_{i=3}^n \epsilon_i s_i \right) \left(s_0+ \sum_{i=3}^n s_i -1 \right)!  \prod_{i=3}^n \pi_{r_i,s_i}^{(\epsilon_i)}(m_i)\ ,
    \end{split}
    \end{equation}
    and in the case $w^{(0)}_1=0$
    \begin{equation}
    \label{eq:Fw2}
    \begin{split}
    &\mathrm{Card}\left(\mathcal{F}_{(w^{(0)}_1=0,w^{(0)}_2,\ldots,w^{(0)}_{r_0+s_0+1})}(m_3,\ldots,m_n)\right)=\\
    &\quad \sum_{\substack{\epsilon_3,\ldots,\epsilon_n=0,1 \\ r_3,\ldots,r_n,s_3,\ldots,s_n\geq 0 \\
        \sum_{i=3}^n \epsilon_i=r_0+\sum_{i=3}^n r_i \\
        \sum_{i=3}^n \bar{\epsilon}_i=s_0+\sum_{i=3}^n s_i+1} }
    \left(s_0+ \sum_{i=3}^n s_i \right)! \left(\bar{\epsilon}_3 r_0 + \sum_{i=3}^n \bar{\epsilon}_i r_i \right) \left(r_0+ \sum_{i=3}^n r_i -1 \right)!  \prod_{i=3}^n \pi_{r_i,s_i}^{(\epsilon_i)}(m_i)\ . 
    \end{split}
    \end{equation}
Note that these expressions
do not depend on the precise sequence $(w^{(0)}_1,\ldots,w^{(0)}_{r_0+s_0+1})$ but only on the values of
$r_0=\sum_{j=2}^{r_0+s_0+1}w^{(0)}_j$, of $s_0=\sum_{j=2}^{r_0+s_0+1}\bar{w}^{(0)}_j$ and of $w^{(0)}_1$.
This is expected since permuting the terms of the sequence $(w^{(0)}_j)_{j=2,\ldots,n}$ at fixed $r_0$ and $s_0$ 
simply amounts to changing the order of the last $r_0+s_0$ trees in the forest $\mathbf{f}$, 
a harmless operation as far as enumeration is concerned. 
We may therefore use the counting formulas \eqref{eq:pirs1} and \eqref{eq:pirs2} for petal necklaces $\mathbf{m}_{12}$
with fixed $r_0$, $s_0$ and with a fixed value $w^{(0)}_1=1$
and $w^{(0)}_1=0$ respectively. The expression \eqref{eq:gennonbip} for $N_{0,n}(2m_1,\ldots,2m_n)$
is obtained by inserting \eqref{eq:Fw1} into \eqref{eq:pirs1} and \eqref{eq:Fw2} into \eqref{eq:pirs2}, adding these two
contributions and finally summing over $r_0$ and $s_0$, which
trivially amounts to replacing in \eqref{eq:Fw1} and \eqref{eq:Fw2} 
each occurrence of $r_0$ by $r_1+r_2$ and each occurrence of $s_0$ by
$s_1+s_2$, and summing over $r_1,r_2,s_1,s_2\geq 0$. This ends the proof of 
Proposition~\ref{prop:gennonbip}.
\end{proof}

Even though it is not apparent, the right-hand side of \eqref{eq:gennonbip} turns out to be, as expected,
symmetric upon permuting the $m_i$'s for $i$ in the whole set $\{1,\ldots,n\}$. This property is shown in 
Appendix~\ref{app:symmetry} where, after some algebraic manipulations, this quantity
is given a manifestly symmetric form, see for instance \eqref{eq:Pifin}  or \eqref{eq:Piinter}. Using this result,
we arrive at the following unified theorem encompassing
all the main theorems of the paper (Theorems~\ref{thm:maintheorem}, \ref{thm:quasibiptheorem} and \ref{thm:quasimain}):

\begin{theorem}
  \label{thm:themaintheorem}
  For $n \geq 3$ and $m_1,m_2,\ldots,m_n \in \Z_{\geq 0}/2$, not all equal to zero, the number
  $N_{0,n}(2m_1,2m_2,\ldots,2m_n)$ of planar tight maps with $n$
  boundaries labeled from $1$ to $n$ with respective lengths
  $2m_1,2m_2,\ldots,2m_n$ is given by the symmetric quasi-polynomial
  \begin{equation}
    \label{eq:themaintheorem}
    \begin{split}
          N_{0,n}(2m_1&,2m_2,\ldots,2m_n)= \\
    &\sum_{\left(\substack{\epsilon_1,\ldots,\epsilon_n\\r_1,\ldots,r_n\\s_1,\ldots
       ,s_n}\right)\in I_n}
    \left( \sum_{i=1}^n r_i \right)! \left( \sum_{i=1}^n \epsilon_i s_i \right) \left( \sum_{i=1}^n s_i -1 \right)!
     \prod_{i=1}^n \pi_{r_i,s_i}^{(\epsilon_i)}(m_i)\\
     +&(n-3)!  \sum_{\substack{(r_1,\ldots,r_n)\in \Z_{\geq 0}^n \\ r_1+\cdots+r_n=n-3}} \sum_{1\leq j<\ell \leq n} \pi_{r_j,0}^{(0)}(m_j) \pi_{r_\ell,0}^{(0)}(m_\ell)\  \prod_{\substack{i=1\\
     i\neq j,\ell}}^n \pi_{r_i,0}^{(1)}(m_i)\\
  +&(n-3)! \sum_{\substack{(r_1,\ldots,r_n)\in \Z_{\geq 0}^n \\ r_1+\cdots+r_n=n-3}}  \pi_{r_1,0}^{(-1)}(m_1)\ \prod_{i=2}^n \pi_{r_i,0}^{(1)}(m_i). 
   \end{split}
  \end{equation}
  with $\pi_{r,s}^{(\epsilon)}(m)$ as in \eqref{eq:pieps} and $I_n$ as in \eqref{eq:7}. Note that the sum in the last line is equal to the symmetric
polynomial $\nobreak{(n-3)! p_{n-3}(m_1,\ldots,m_n)}$ when all $m_i$ are
integers, and to zero otherwise, hence it is symmetric.
\end{theorem}
\begin{proof}
Since $N_{0,n}(2m_1,2m_2,\ldots,2m_n)$ is (by definition) symmetric in $m_1,m_2,\ldots,m_n$
and so is the right-hand side of \eqref{eq:themaintheorem}, we may assume without loss of
  generality that $m_1>0$ and either $m_2>0$ or
  $m_2=\cdots=m_n=0$. In the first case, Proposition~\ref{prop:gennonbip} and the identification of
  the right-hand sides of \eqref{eq:gennonbip} and \eqref{eq:themaintheorem} proved in Proposition \ref{prop:Pifin} of Appendix~\ref{app:symmetry}
  allows to conclude. In the second case, both sides of the equation reduce to the univariate quasi-polynomial $(n-3)!\pi^{(-1)}_{n-3,0}(m_1)$, equal to
  $(n-3)! p_{n-3}(m_1)$ if $m_1$ is a positive integer, and to zero otherwise.
\end{proof}
Let us now explain how to recover Theorems~\ref{thm:maintheorem}, \ref{thm:quasibiptheorem} and \ref{thm:quasimain}
from Theorem~\ref{thm:themaintheorem}. As in Remark~\ref{rem:unification}, let us denote by $k$ the number of half-integers
among $m_1,\ldots,m_n$. From the general property that $\pi_{r_i,s_i}^{(\epsilon_i)}(m_i)$ is non-zero only if
 $m_i-\frac{s_i+1+\epsilon_i}2 \in \Z$, we see that the sum in the third line of the right-hand side of \eqref{eq:themaintheorem} is non-zero
 if and only if all the $m_i$'s are integers, i.e.~the map is bipartite ($k=0$), in which case it reduces precisely to
 $(n-3)! p_{n-3}(m_1,\ldots,m_n)$. Similarly, the sum in the second line is non-zero
 if and only if exactly two of the $m_i$'s, say for $i=i_1$ and $i=i_2$ ($1\leq i_1< i_2 \leq n$), are half-integers, i.e.~the map is 
 quasi-bipartite ($k=2$), in which case it reduces to
 $(n-3)! \tilde{p}_{n-3}(m_{i_1},m_{i_2};m_1,\ldots,\check{m}_{i_1},\ldots,\check{m}_{i_2},\ldots,m_n)$
 where $\check{m}_{i_\ell}$ means that the argument $m_{i_\ell}$ is omitted.
 Therefore, for $k=1$ and $k\geq 3$, the only possibly non-zero term is the sum in the first line, which matches precisely  
the right hand side of \eqref{eq:quasimain}. This proves Theorem~\ref{thm:quasimain} (where in practice, only even values of $k$ yield a non-zero result, as it should). In order to recover 
Theorems~\ref{thm:maintheorem} and \ref{thm:quasibiptheorem}, it only remains to check that the sum in the first line of the right-hand side
of \eqref{eq:themaintheorem} vanishes in the bipartite and quasi-bipartite cases. Note that if $m_i$ is an integer, 
$\pi_{r_i,s_i}^{(\epsilon_i)}(m_i)$ is non-zero only if $s_i$ and $1-\epsilon_i$ have the same parity, a property 
which, for $s_i\geq 0$ and 
$1-\epsilon_i\in\{0,1\}$ implies that $s_i\geq 1-\epsilon_i$. If all the $m_i$'s are integers, the required constraint
$\sum_{i=1}^n(1-\epsilon_i)=\sum_{i=1}^n s_i+2$ in the definition of the set $I_n$ cannot be fulfilled, hence the 
sum vanishes. If exactly two of the $m_i$'s, say for $i=i_1$ and $i=i_2$, are half-integers, this same constraint 
imposes that $\epsilon_{i_1}=\epsilon_{i_2}=0$ and $s_{i_1}=s_{i_2}=0$ while $s_i=1-\epsilon_i$ 
for all $i\neq i_1,i_2$, hence $\left(\sum_{i=1}^n \epsilon_i s_i \right)=\sum_{i\neq i_1,i_2}(1-\epsilon_i)\epsilon_i=0$
and the first sum again vanishes.
To summarize, for an even value of $k$, exactly one of the three lines in the right-hand side of \eqref{eq:themaintheorem} 
is not identically zero. Each line corresponds to one of the mutually exclusive situations
$k\geq 4$ (first line), $k=2$ (second line) and $k=0$ (third line), corresponding to Theorem~\ref{thm:quasimain}, \ref{thm:quasibiptheorem} and \ref{thm:maintheorem} respectively.

\section{Conclusion}
\label{sec:conclusion}

In this paper, we have provided an explicit expression for the planar lattice
count quasi-polynomial $N_{0,n}(2m_1,\ldots,2m_n)$, by extending the
slice decomposition  to
the case of planar tight maps. Note that, in contrast with previous
work such as \cite{irredmaps,Budd2020a}, our
approach is not based on generating functions but rather on direct
bijective enumeration techniques. 

We note that our most general formula, stated in Theorem
\ref{thm:themaintheorem}, still involves a complicated sum over a no
less complicated set $I_n$. As discussed in the Appendix
\ref{app:symmetry}, there are in fact many identities that can be used
to rewrite it, and it is not impossible that it admits a
significantly simpler expression yet to be unveiled. A similar remark
applies to our extension of Tutte's slicings formula given in Theorem
\ref{thm:extendedslicings}. 

We believe that the methodology of Sections~\ref{sec:bijective} and
\ref{sec:extens-non-bipart} is quite robust and may be adapted to
other map enumeration problems.  A first problem one may think of is a
model of maps with \emph{continuous} edge lengths, whose set can be
associated with a natural volume measure. The volumes of such measures
have been considered by Kontsevich \cite{Kontsevich92}, and, as noted
by Norbury, they correspond to the homogeneous top degree part of the
lattice count polynomials,
see~\cite[Theorem~3]{Norbury2010}. Combining with our
Theorem~\ref{thm:maintheorem}, we obtain the following result (see
Norbury's paper for the definition of $V_{0,n}$).
\begin{proposition}
  The genus zero volume polynomial $V_{0,n}$ reads explicitly
  \begin{equation}
    V_{0,n}(b_1,\ldots,b_n) = \frac{(n-3)!}{2^{2n-7}}
    \sum_{k_1,\ldots,k_n\geq 0\atop k_1+\cdots+k_n=n-3}
    \left( \frac{b_1^{k_1} \cdots b_n^{k_n}}{k_1! \cdots k_n!}\right)^2.
  \end{equation}
\end{proposition}
This amounts to the known expression for the genus zero intersection
numbers, see for instance~\cite[Proposition~4.6.10]{Lando2004}
and~\cite[Corollary~1]{Norbury2010}. Still, a direct construction of
maps with continuous edge lengths, using continuous paths and
mimicking the above discrete paths encoding for slices, should also be
possible and interesting, and we plan to investigate this question in
future work. A second extension is motivated by the work of Budd
\cite{Budd2020a} who generalized Norbury's results to the case of
\emph{irreducible maps}, i.e.~maps with a girth constraint. The slice
decomposition of these maps was discussed in \cite{irredmaps} in the
non-tight case, and we expect it to be adaptable to the tighness
constraint. Finally, combining these two ideas, namely, passing to
maps with continuous edge lengths \emph{and} adding an irreducibility
constraint, one obtains \emph{irreducible metric maps} which have been
considered by Budd \cite{budd2020irreducible}, who showed that their
volumes are related to the Weil-Petersson volumes of hyperbolic
surfaces. We plan to investigate the slice decomposition of these
maps, which might shed new light on these questions.

In another direction, the bijective techniques presented in this paper
should pave the road to the study of continuum limits of random planar
tight maps, when the number of faces tends to infinity. We believe in
particular that, as soon as the face degrees are well-behaved in a
certain sense, the Gromov-Hausdorff limit of appropriately renormalized planar tight
maps, seen as discrete metric spaces, should be given by the
\emph{Brownian sphere}. To this purpose, the recent approach by
Marzouk \cite{Marzouk18,marzouk2022}, dealing with limits of planar non necessarily
tight maps with prescribed face degrees, should be particularly
relevant.

Note that our approach is currently restricted to the planar
(i.e.~genus $0$) case. This is a current limitation of the slice
decomposition. We hope however that this limitation will be challenged
by further investigations. A first result in this vein is the
bijective study of pairs of pants (planar maps with three boundaries)
done in \cite{triskell}, and the fact that general surfaces can be decomposed
into pairs of pants gives some support to our hope.

\appendix

\section{Enumeration of one- and two-type labeled plane
  forests}\label{sec:enumeration-two-type}

This appendix lists the forest enumeration results that we need in
this paper. We consider \emph{plane forests} (sequences of rooted
plane trees) which are \emph{labeled} (distinct labels are assigned to
the vertices). We start with the case of forests with one type of
vertices.

\begin{proposition}
  \label{prop:onetypeforest}
  Let $n,k_1,\ldots,k_n$ be non-negative integers such that
  $k_1+\cdots+k_n<n$. Then, there are exactly $(n-1)!$ plane forests
  with $n$ vertices labeled $\{1,\ldots,n\}$, such that vertex $i$ has
  $k_i$ children for all $i=1,\ldots,n$ and such that vertex $1$
  appears in the first tree. (Such forests consist necessarily of
  $k_0:=n-k_1-\cdots-k_n$ trees.)
\end{proposition}

\begin{proof}
  The case of trees ($k_0=1$) is given explicitly in~\cite[Section 5,
  Equation~(18)]{Bernardi2014}. It implies the general case since, for
  $k_0>1$, there is a straightforward bijection between the set of
  forests at hand and the set of labeled plane trees such that vertex
  $1$ has $k_1+k_0-1$ children, the number of children of the other
  vertices being unmodified.
\end{proof}

\begin{remark}
  \label{rem:onetypeforest}
  Proposition~\ref{prop:onetypeforest} can alternatively be proved
  directly by exhibiting a bijection between the set of forests at
  hand and the set of cyclic orders on $\{1,\ldots,n\}$. Such a
  bijection is obtained by simply listing the vertex labels of a
  labeled plane forest in depth-first order, giving a linear, hence
  a cyclic, order on $\{1,\ldots,n\}$.  Conversely, from a cyclic order
  and the data of the $k_i$'s, we construct a conjugacy class of
  Łukasiewicz words---see e.g.\ \cite[Section~5.3]{StanleyEC2}---whose
  letters are labeled and which contains exactly one word coding for a
  plane forest having vertex $1$ in the first tree.
\end{remark}

We now turn our attention to labeled plane forests with two types of
vertices, say $A$ and $B$. Our purpose is to enumerate such
\emph{two-type} forests in which, for every vertex, we prescribe not
only its type but also the sequence formed by the types of all its
children, read in the planarity order. Similar counting problems, for
an arbitrary number of types, have been previously considered in the
literature---see e.g.\ \cite{Bacher2013,Bernardi2014,Chaumont2016} and
references therein---but since the general formulas are quite
complicated we provide a self-contained derivation for two types. Our
approach is closely related with that of Chottin~\cite{Chottin1981}
who treated the case of two-type trees. Handling forests with several
components however involves an extra difficulty, which we circumvent
by specializing the general multitype approach of Bacher and
Schaeffer~\cite{Bacher2013}. Note that the latter two references
consider plane forests which are unlabeled, but adding vertex labels
does not fundamentally change the problem since plane forests have no
symmetries.

To state our result we need some definitions and notation. Let us
consider a two-type plane forest whose vertices are labeled
$\{1,\ldots,n\}$. For every vertex $i=1,\ldots,n$, we denote by $k_i$
its number of children, and we let
$w^{(i)}=(w_0^{(i)},w_1^{(i)},\ldots,w_{k_i}^{(i)}) \in
\{A,B\}^{k_i+1}$ be the sequence such that $w_0^{(i)}$ is the type of
$i$ and such that, for every $j=1,\ldots,k_i$, $w_j^{(i)}$ is the type
of the $j$-th child of $i$ in the planarity order.  We also define a
sequence $w^{(0)}=(w_0^{(0)},w_1^{(0)},\ldots,w_{k_0}^{(0)})$ where
$k_0$ is the number of trees of the forest, $w_0^{(0)}$ is a third
type denoted $O$ and, for $j=1,\ldots,k_0$, $w_j^{(0)}$ is the type of
the $j$-th \emph{root}, i.e. of the root vertex of the $j$-th tree of
the forest. In this sense, $0$ can be seen as the label of a
\emph{super-root} of type $O$, which is the parent of all the roots
(which have the usual types $A$ or $B$).

The collection $\boldsymbol{w}=(w^{(0)},w^{(1)},\ldots,w^{(n)})$ is
called the \emph{type array} of the forest. Denoting by $\iv{\cdot}$
the Iverson bracket ($\iv{P}$ is equal to $1$ if $P$ is true, and to
$0$ otherwise), we have necessarily
\begin{equation}
  \label{eq:Aconsis}
  \sum_{i=1}^n \iv{w_0^{(i)}=A} = \sum_{i=0}^n \sum_{j=1}^{k_i} \iv{w_j^{(i)}=A} := a
\end{equation}
as seen by expressing in two different ways the number $a$ of type $A$
vertices. Similarly, the number $b$ of type $B$ vertices is given by
\begin{equation}
  \label{eq:Bconsis}
  \sum_{i=1}^n \iv{w_0^{(i)}=B} = \sum_{i=0}^n \sum_{j=1}^{k_i} \iv{w_j^{(i)}=B} := b.
\end{equation}
Note that, by adding these two equations, we obtain the relation
$n=k_0+k_1+\cdots+k_n$ already seen in the context of one-type forests
in Proposition~\ref{prop:onetypeforest}. An array
$\boldsymbol{w}=(w_j^{(i)})^{i=0,\ldots,n}_{j=0,\ldots,k_i}$ with
$w_0^{(0)}=O$ and $w_j^{(i)}\in \{A,B\}$ for $i,j$ not both zero is
said \emph{consistent} if it satisfies \eqref{eq:Aconsis} and
\eqref{eq:Bconsis}.

\begin{proposition}
  \label{prop:twotypeforest}
  Let $\boldsymbol{w}$ be a consistent array, and let $a$ and $b$ be
  as defined in~\eqref{eq:Aconsis} and \eqref{eq:Bconsis},
  respectively. Define furthermore the integers
  \begin{align}
    a^O &:= \sum_{j=1}^{k_0} \iv{w_j^{(0)}=A}, & a^B &:= \sum_{i=1}^n \sum_{j=1}^{k_i} \iv{w_0^{(i)}=B} \iv{w_j^{(i)}=A}, \\
    b^O &:= \sum_{j=1}^{k_0} \iv{w_j^{(0)}=B}, & b^A &:= \sum_{i=1}^n \sum_{j=1}^{k_i} \iv{w_0^{(i)}=A} \iv{w_j^{(i)}=B},
  \end{align}
  namely $a^O$ and $b^O$ correspond to the number of type $A$ and type
  $B$ roots, respectively, while $a^B$ is the number of type $A$
  vertices with a type $B$ parent, and vice versa for $b^A$.

  Then, we have the following enumerative formulas.
  \begin{itemize}
  \item (General enumeration) The number of two-type labeled plane
    forests of type array $\boldsymbol{w}$---thus containing $a$ type
    $A$ and $b$ type $B$ vertices---is equal to
    \begin{equation}
      \label{eq:twotypeforestgen}
      \begin{cases}
        (a^O b^O + a^O b^A + a^B b^O) (a-1)! (b-1)!
        & \text{if $a>0$ and $b>0$,}\\
        a^O (a-1)! & \text{if $a>0$ and $b=0$,}\\
        b^O (b-1)! & \text{if $a=0$ and $b>0$,}\\
        1 & \text{if $a=0$ and $b=0$.}
      \end{cases}
    \end{equation}
  \item (Constrained enumeration) Assume that, say, $a>1$ and
    $w_0^{(1)}=A$, i.e.\ vertex $1$ has type $A$. Then, the number of
    two-type labeled plane forests of type array $\boldsymbol{w}$ such
    that vertex $1$ is in the first tree is equal to
    \begin{equation}
      \label{eq:twotypeforestconst}
      \begin{cases}
        (b^O + b^A) (a-1)! (b-1)! & \text{if $b>0$ and $w_1^{(0)}=A$,}\\
        a^B (a-1)! (b-1)! & \text{if $b>0$ and $w_1^{(0)}=B$,}\\
        (a-1)! & \text{if $b=0$.}
      \end{cases}
    \end{equation}
    (The first and second cases correspond to a first root of type $A$
    and $B$, respectively.)
  \end{itemize}
\end{proposition}

\begin{proof}
  Let us first note that the cases where $a$ or $b$ vanish follow
  immediately from Proposition~\ref{prop:onetypeforest} (for the
  general enumeration, we perform a circular permutation of the trees
  to lift the constraint that vertex $1$ is in the first tree, giving
  the extra factor $a^O$ or $b^O$). Hence, we assume from now on that
  $a$ and $b$ are both positive.

  We claim that the general enumeration formula follows from the
  constrained one. Indeed, we may partition the set of forests of type
  array $\boldsymbol{w}$ according to the index $j=1,\ldots,k_0$ of
  the tree containing vertex $1$. Upon doing a circular permutation of
  the trees, we deduce from~\eqref{eq:twotypeforestconst} that, for
  $w_0^{(1)}=A$, the number of forests with a given value of $j$ is
  equal to $(b^O + b^A) (a-1)! (b-1)!$ if $w_j^{(0)}=A$, and to
  $a^B (a-1)! (b-1)!$ if $w_j^{(0)}=B$.  As the first (resp.\ second)
  case occurs for $a^O$ (resp.\ $b^O$) values of $j$, summing over $j$
  gives the first line of~\eqref{eq:twotypeforestgen}. The case
  $w_0^{(1)}=B$ is deduced by exchanging the roles of $A$ and $B$.

  \begin{figure}[t]
    \centering
    \includegraphics[width=\textwidth]{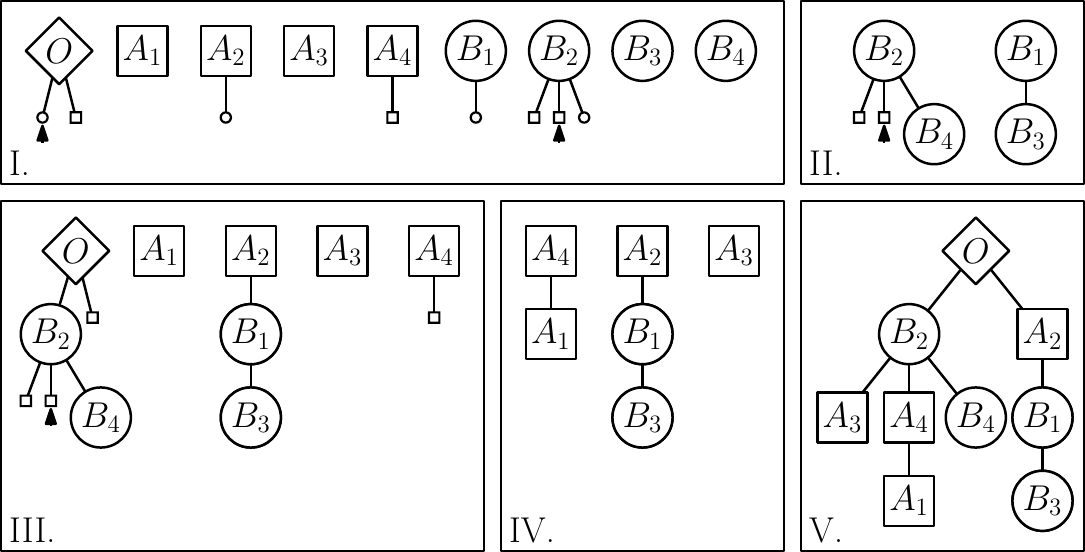}
    \caption{Construction of a two-type forest in the case
      $w_0^{(1)}=B$. See the main text for a precise description of
      the stages I-V. We display the type $A$ (resp.\ $B$) vertices
      and dangling $A$-edges (resp.\ $B$-edges) with squares (resp.\
      circles). The distinguished dangling edges are indicated with
      arrows.}
    \label{fig:twotypeforest}
  \end{figure}
  
  It remains to prove the constrained enumeration formula. We will do
  so by giving an algorithm to construct any forest of type array
  $\boldsymbol{w}$, where it will be manifest that the number of
  possibilities is given by~\eqref{eq:twotypeforestconst}. The
  algorithm is easier to visualize if, instead of working with labels
  $0,1,\ldots,n$, we relabel the vertices as
  $O,A_1,\ldots,A_a,B_1,\ldots,B_b$ to make their type apparent (as
  the type array $\boldsymbol{w}$ is fixed, this may be done by fixing
  a bijection between $\{0,1,\ldots,n\}$ and
  $\{O,A_1,\ldots,A_a,B_1,\ldots,B_b\}$ hence does not change the
  counting problem).  The general idea, illustrated on
  Figure~\ref{fig:twotypeforest}, is to proceed in several stages, by
  first ``assembling'' the type $B$ vertices together, before dealing
  with the types $A$ and $O$ vertices. Let us describe the different
  stages of the algorithm in detail.
  \begin{itemize}
  \item[I.] We start with isolated vertices
    $O,A_1,\ldots,A_a,B_1,\ldots,B_b$, to which we attach sequences of
    \emph{dangling $A$-edges} and \emph{$B$-edges}: these dangling
    edges will be connected later to type $A$ and $B$ vertices,
    respectively. The type array $\boldsymbol{w}$ tells us precisely
    the sequence which we have to attach to each vertex.  We define a
    total order on the dangling edges by listing those incident to
    $O$, then those incident to $A_1$, etc. We distinguish the
    \emph{first} edge incident to $O$: it has type $w_0^{(1)}$. We
    then distinguish another dangling edge of the opposite type:
    \begin{itemize}
    \item if $w_0^{(1)}=A$, then we distinguish a $B$-edge incident
      either to $O$ or to a type $A$ vertex: there are $b^O+b^A$
      possible choices,
    \item if $w_0^{(1)}=B$, then we distinguish an $A$-edge incident
      to a type $B$ vertex: there are $a^B$ possible choices.
    \end{itemize}
  \item[II.] We form a plane forest with the vertices
    $B_1,\ldots,B_b$, by attaching them together via their incident
    $B$-edges. By Proposition~\ref{prop:onetypeforest}, there are
    $(b-1)!$ ways to do so (with $B_1$ in the first tree). A simple
    computation shows that the resulting forest is made of $b^O+b^A$
    trees. If $w_0^{(1)}=B$, we permute the trees circularly so that
    the tree containing the distinguished $A$-edge comes first.
  \item[III.] We attach the roots of the forest constructed at stage
    II to the $B$-edges (in number $b^O+b^A$) dangling from
    $O,A_1,\ldots,A_a$. Precisely, we attach the root of the first
    tree to the distinguished $B$-edge, and we then proceed circularly
    using the order on dangling edges defined at stage I.  Note that
    all the $B$-edges have now been matched to type $B$ vertices. We
    end up with a sequence of \emph{supernodes}, which are trees with
    roots $O,A_1,\ldots,A_a$, subtrees made of type $B$ vertices, and
    dangling $A$-edges. Observe that the distinguished $A$-edge always
    belongs to the supernode with root $O$ (precisely, the first edge
    incident to $O$ is either the distinguished $A$-edge when
    $w_0^{(1)}=A$, or leads to a subtree which contains it when
    $w_0^{(1)}=B$).
  \item[IV.] Viewing the supernodes as compound vertices, we form a
    plane forest with $a^O$ trees by assembling the supernodes with roots
    $A_1,\ldots,A_a$ via their dangling $A$-edges. Note that the
    number of ``children'' of a supernode is prescribed by its number
    of dangling $A$-edges. By Proposition~\ref{prop:onetypeforest},
    there are $(a-1)!$ ways to do so (with $A_1$ in the first tree).
  \item[V.] We then complete the construction by attaching these trees
    to the $A$-edges of the supernode with root $O$, starting with the
    first tree which we attach to the distinguished $A$-edge, and then
    proceeding circularly using the order of stage I. All dangling
    edges have now been matched, and vertex $1$ is in the first tree
    by construction.
  \end{itemize}
  It is plain from stage I that we obtain a forest of type array
  $\boldsymbol{w}$. Furthermore, each such forest is obtained in
  precisely one way, as we may check that it is obtained for a unique
  choice of the second distinguished dangling edge at stage I
  \footnote{Precisely, if $w_1^{(0)}=A$ then the distinguished
    $B$-edge corresponds to the edge closest to $B_1$ attached to a
    parent not of type $B$, and if $w_1^{(0)}=B$ then the
    distinguished $A$-edge corresponds to the parent edge of the
    oldest type $A$ ancestor of $A_1$.} and of the one-type forests at
  stages II and IV.
\end{proof}

\begin{remark}
  The factor $a^O b^O + a^O b^A + a^B b^O$
  in~\eqref{eq:twotypeforestgen} corresponds to a sum over the three
  Cayley trees on the set $\{O,A,B\}$. For $k$ types of vertices there
  would be as many terms as Cayley trees on a set with $k+1$
  elements~\cite{Bacher2013}.
\end{remark}

\section{Symmetrizing the planar lattice count quasi-polynomials}
\label{app:symmetry}
Let us fix an integer $n \geq 3$ and denote by $\Pi(m_1,\ldots,m_n)$
the right-hand side of \eqref{eq:gennonbip}. It is manifest that it is
a quasi-polynomial in $2m_1,\ldots,2m_n$ of degree $n-3$, since
$\pi_{r,s}^{(\epsilon)}(m)$ is a univariate quasi-polynomial in $2m$
of degree $r+s$. The purpose of this appendix is to show that
$\Pi(m_1,\ldots,m_n)$ is in fact \emph{symmetric} in $m_1,\ldots,m_n$,
which we will do by rewriting it in a manifestly symmetric form. Note
that the expression \eqref{eq:gennonbip} displays a symmetry in
$m_4,\ldots,m_n$ only.  Here we assume only that $m_i\in \Z/2$,
$i=1,\ldots,n$ without further restriction.

It is useful to introduce compact notations for the
high-dimensional sums appearing in \eqref{eq:gennonbip}. Let $I$ be the finite subset
of $\{0,1\}^n \times \Z_{\geq 0}^n \times \Z_{\geq 0}^n$ defined by
\begin{equation}\label{eq:Idef}
  I := \left\{ \left(
      \begin{array}{c}
        \epsilon_1,\ldots,\epsilon_n\\r_1,\ldots,r_n\\s_1,\ldots
       ,s_n
     \end{array}
   \right):
    \begin{array}{l}
      \sum\limits_{i=1}^n \epsilon_i=\sum\limits_{i=1}^n r_i +1  \\
      \sum\limits_{i=1}^n (1-\epsilon_i)=\sum\limits_{i=1}^n s_i+2
    \end{array}
\right\}\, .
\end{equation}
Given
a tuple in $I$, we set
\begin{equation}
  r := \sum_{i=1}^n r_i, \qquad s:= \sum_{i=1}^n s_i,
\end{equation}
and note that $r+s=n-3$ by the definition of $I$. Let $I_{r\geq 1}$,
$I_{r=0}$, $I_{s\geq 1}$, and $I_{s=0}$, be the subsets of $I$ consisting
of tuples such that $r\geq 1$, $r=0$, $s\geq 1$, and $s=0$, respectively. 
Note that the set $I_n$ introduced in \eqref{eq:7} corresponds in our
present notations to $I_{s\geq 1}$. Note also that 
$I_{r=0}$ consists of tuples such that $r_i=0$ for all $i$,
exactly one $\epsilon_i$ is equal to $1$, and
$s_1+\cdots+s_n=s=n-3$. Similarly, $I_{s=0}$ consists of tuples such
that $s_i=0$ for all $i$, exactly two $\epsilon_i$ are equal to $0$,
and $r_1+\cdots+r_n=r=n-3$.

For $e_1,e_2 \in \{0,1\}$, we let $I^{(e_1,e_2)}$ be the subset of $I$
consisting of tuples such that $\epsilon_1=e_1$ and $\epsilon_2=e_2$.
We also allow for the value $e_1=-1$, which means that we consider
tuples such that $\epsilon_1=-1$, keeping the sum condition in
\eqref{eq:Idef} unchanged.
The notations $I^{(e_1,e_2)}_{r\geq 1}$, etc, should hopefully be
self-explanatory. Note that $I^{(-1,1)}_{s=0}$ consists of tuples such
that $s_i=0$ for all $i$, $\epsilon_1=-1$ and $\epsilon_i=1$ for all $i \geq 2$, and
$r_1+\cdots+r_n=r=n-3$. We finally use, for $\epsilon\in \{0,1\}$, the shorthand notation
$\bar \epsilon :=1-\epsilon$ (which we shall never use for $\epsilon=-1$).

Armed with all these notations, we can rewrite
$\Pi(m_1,\ldots,m_n)$\footnote{Here we use the shorthand notation $\sum\limits_{J}(\cdot) $ for $\sum\limits_{\left(\substack{\epsilon_1,\ldots,\epsilon_n\\r_1,\ldots,r_n\\s_1,\ldots
       ,s_n}\right)\in J}(\cdot)$.} as
\begin{equation}
  \label{eq:Pirew}
  \begin{split}
  \Pi(m_1,\ldots,m_n) =&
  \sum_{I^{(-1,1)}_{s\geq 1} \cup I^{(0,0)}_{s\geq 1}} r! (s-1)! \left(\epsilon_3 (s_1+s_2) + \sum_{j=3}^n \epsilon_j s_j \right) \prod_{i=1}^n \pi_{r_i,s_i}^{(\epsilon_i)}(m_i) \\
  +&
  \sum_{I^{(0,1)}_{r\geq 1} \cup I^{(1,0)}_{r\geq 1}} (r-1)! s! \left(\bar\epsilon_3 (r_1+r_2) + \sum_{j=3}^n \bar\epsilon_j r_j \right) \prod_{i=1}^n \pi_{r_i,s_i}^{(\epsilon_i)}(m_i) \\
  +& \sum_{I^{(-1,1)}_{s=0} \cup I^{(0,0)}_{s=0} \cup
    I^{(0,1)}_{r=0} \cup I^{(1,0)}_{r=0}} r! s! \prod_{i=1}^n
  \pi_{r_i,s_i}^{(\epsilon_i)}(m_i).
  \end{split}
\end{equation}
Here, the first (resp.~second) sum corresponds to the first (resp.~second) sum in \eqref{eq:gennonbip} when 
at least one of the $s_i$'s (resp.~one of the $r_i$'s) is non zero, hence when $s\geq 1$ (resp.~$r\geq 1$);
the third sum accounts for the conventional values in \eqref{eq:gennonbip}: $(\epsilon_3 (s_1+s_2) + \sum_{i=3}^n \epsilon_i s_i) ( \sum_{i=1}^n s_i -1)!\to 1=s!$ when all the $s_i$'s are zero (or equivalently $s=0$) and 
$(\bar{\epsilon}_3 (r_1+r_2) + \sum_{i=3}^n \bar{\epsilon}_i r_i) ( \sum_{i=1}^n r_i -1)!\to 1=r!$ when
all the $r_i$'s are zero (or equivalently $r=0$).
 
 \medskip
The main result of this appendix is:
\begin{proposition}
  \label{prop:Pifin}
  The quasi-polynomial $\Pi(m_1,\ldots,m_n)$ is symmetric in
  $m_1,\ldots,m_n$ and admits the expression
  \begin{multline}
    \label{eq:Pifin}
    \Pi(m_1,\ldots,m_n) = \\ \sum_{I_{s\geq 1}} r! (s-1)! \left(\sum_{j=1}^n
      \epsilon_j s_j \right)
    \prod_{i=1}^n \pi_{r_i,s_i}^{(\epsilon_i)}(m_i) \quad +
     \sum_{I^{(-1,1)}_{s=0} \cup I_{s=0}} r! s! \prod_{i=1}^n
    \pi_{r_i,s_i}^{(\epsilon_i)}(m_i).
  \end{multline}
  Note that the sum over $I^{(-1,1)}_{s=0}$ is equal to $(n-3)! \sum\limits_{r_1+\cdots+r_n=n-3} \pi^{(-1)}_{r_1,0}(m_1)\prod\limits_{i=2}^n \pi_{r_i,0}^{(1)}(m_i)$ which is equal to the symmetric
  polynomial $(n-3)! p_{n-3}(m_1,\ldots,m_n)$ when all $m_i$ are
  integers, and to zero otherwise, hence it is symmetric like the
  rest.
\end{proposition}

The expression above for $\Pi(m_1,\ldots,m_n)$ is precisely the
right-hand side of \eqref{eq:themaintheorem}.  Indeed, we already
noticed that $I_{s\geq 1}$ is nothing but the set $I_n$ defined in
\eqref{eq:7}, hence the sum over $I_{s\geq 1}$ gives the first term in
the right-hand side of \eqref{eq:themaintheorem}. Furthermore, the set
$I_{s=0}$ consists of tuples such that $s_i=0$ for all $i$, exactly
two $\epsilon_i$ are equal to $0$, and $r_1+\cdots+r_n=r=n-3$, hence
the sum over $I_{s=0}$ corresponds to the second term in
\eqref{eq:themaintheorem}.  Finally, the sum over $I^{(-1,1)}_{s=0}$
corresponds to the third term in \eqref{eq:themaintheorem}.

\medskip
In order to prove Proposition~\ref{prop:Pifin}, we will first record
the following:
\begin{proposition}
  The univariate polynomials $\pi_{r,s}^{(\epsilon)}(m)$ satisfy the
  relations
  \begin{align}
    \label{eq:nspieps}
    s \pi_{r,s}^{(\epsilon)}(m) &= (r+1) \pi_{r+1,s-1}^{(\epsilon+1)}(m), \\
    \label{eq:npi1} \pi_{r,s}^{(1)}(m) &= \pi_{r,s}^{(-1)}(m) + \pi_{r-1,s}^{(-1)}(m), \\
    \label{eq:nspi1} s \pi_{r,s}^{(1)}(m) &= (r+1) \pi_{r+1,s-1}^{(0)}(m) + r \pi_{r,s-1}^{(0)}(m),
  \end{align}
  valid for all $r,s \geq 0$ and all integer $\epsilon$, with the convention that
  $\pi^{(\epsilon)}_{-1,s} =\pi^{(\epsilon)}_{r,-1}=0$.
\end{proposition}

\begin{proof}
  The first relation follows immediately from the mere
  definition~\eqref{eq:pieps} of $\pi_{r,s}^{(\epsilon)}(m)$. For the
  two other ones, we make use of the ``dilaton-like'' equation
  \begin{equation}
    k  p_{k,e+2}(m) = k p_{k,e}(m) + (k-e) p_{k-1,e}(m)
  \end{equation}
  which can be checked from the definition~\eqref{eq:pkedef1} of
  $p_{k,e}(m)$ and is valid for all $k \geq 0$, with the convention
  that $p_{-1,e}(m)=0$. This dilaton-like equation implies
  immediately~\eqref{eq:npi1}, and to get~\eqref{eq:nspi1} we first
  apply~\eqref{eq:nspieps} at $\epsilon=1$ then the dilaton-like
  equation to go back from $\pi^{(2)}$'s to $\pi^{(0)}$'s.
\end{proof}

\begin{lemma}[Transmutation relations]
  \label{lem:trans}
  For any $j=1,\ldots,n$, we have
  \begin{equation}
    \label{eq:trans1}
    \sum_{I_{s\geq 1}} r! (s-1)! (\bar\epsilon_j s_j)   \prod_{i=1}^n
    \pi_{r_i,s_i}^{(\epsilon_i)}(m_i) =
    \sum_{I_{r\geq 1}} (r-1)! s! (\epsilon_j r_j)  \prod_{i=1}^n
    \pi_{r_i,s_i}^{(\epsilon_i)}(m_i).
  \end{equation}
We also have the four identities
\begingroup
\allowdisplaybreaks
\begin{align}
    \label{eq:trans21}
\sum_{I^{(-1,1)}_{s\geq 1}} r! (s-1)!\  \epsilon_3 s_1 \prod_{i=1}^n \pi_{r_i,s_i}^{(\epsilon_i)}(m_i) &= 
\sum_{I^{(0,1)}_{r\geq 1}} (r-1)! s!\  \epsilon_3 r_1 \prod_{i=1}^n \pi_{r_i,s_i}^{(\epsilon_i)}(m_i),\\ \label{eq:trans22}
\sum_{I^{(-1,1)}_{s\geq 1}} r! (s-1)!\  \epsilon_3 s_2 \prod_{i=1}^n \pi_{r_i,s_i}^{(\epsilon_i)}(m_i) &\overset{(*)}{=} \!\!
\sum_{I^{(1,-1)}_{s\geq 1}} \!\! r! (s-1)!\  \epsilon_3 s_2 \prod_{i=1}^n \pi_{r_i,s_i}^{(\epsilon_i)}(m_i) \nonumber \\ &=  
\sum_{I^{(1,0)}_{r\geq 1}} (r-1)! s!\  \epsilon_3 r_2 \prod_{i=1}^n \pi_{r_i,s_i}^{(\epsilon_i)}(m_i),\\ \label{eq:trans23}
\sum_{I^{(1,0)}_{s\geq 1}}  (r-1)! s!\  \bar{\epsilon}_3 r_1 \prod_{i=1}^n \pi_{r_i,s_i}^{(\epsilon_i)}(m_i) &= 
\sum_{I^{(0,0)}_{r\geq 1}} r! (s-1)!\  \bar{\epsilon}_3 s_1 \prod_{i=1}^n \pi_{r_i,s_i}^{(\epsilon_i)}(m_i),\\ \label{eq:trans24}
\sum_{I^{(0,1)}_{s\geq 1}}  (r-1)! s!\  \bar{\epsilon}_3 r_2 \prod_{i=1}^n \pi_{r_i,s_i}^{(\epsilon_i)}(m_i) &= 
\sum_{I^{(0,0)}_{r\geq 1}} r! (s-1)!\  \bar{\epsilon}_3 s_2 \prod_{i=1}^n \pi_{r_i,s_i}^{(\epsilon_i)}(m_i).
  \end{align}
\endgroup
  Finally, for $j \geq 3$ we have
  \begin{equation}
    \label{eq:trans3}
    \sum_{I^{(-1,1)}_{s\geq 1}} r! (s-1)! (\epsilon_j s_j) \prod_{i=1}^n
    \pi_{r_i,s_i}^{(\epsilon_i)}(m_i) =
    \sum_{I^{(1,1)}_{r\geq 1}} (r-1)! s! (\bar\epsilon_j r_j)  \prod_{i=1}^n
    \pi_{r_i,s_i}^{(\epsilon_i)}(m_i).
  \end{equation}
\end{lemma}

\begin{proof}
  In the left-hand side of~\eqref{eq:trans1}, the only contributing
  tuples are those with $\epsilon_j=0$. By applying the
  relation~\eqref{eq:nspieps} to the factor
  $s_j \pi_{r_j,s_j}^{(0)}(m_j)$ appearing in the product, and by
  performing the change of variables $\bar\epsilon_j \to \epsilon_j$,
  $r_j+1 \to r_j$, $s_j-1 \to s_j$ (leaving all other elements of the
  tuples unchanged), we obtain precisely the nonzero terms of the
  right-hand side.

  The proofs of the identities~\eqref{eq:trans21}
  and~\eqref{eq:trans23} are entirely similar, applying
  now~\eqref{eq:nspieps} to the factor
  $s_1 \pi_{r_1,s_1}^{(\epsilon_1)}(m_1)$ (with $\epsilon_1=-1,0$)
  appearing in the left-hand side of~\eqref{eq:trans21} (with
  $\epsilon_1=-1$) and in the right-hand side of~\eqref{eq:trans23}
  (with $\epsilon_1=0$) .  For~\eqref{eq:trans22}
  and~\eqref{eq:trans24}, we proceed in the same way (now with
  $\epsilon_2=-1,0$), after noting that that it is possible to first
  replace $I^{(-1,1)}_{s\geq 1}$ by $I^{(1,-1)}_{s\geq 1}$
  in~\eqref{eq:trans22} (as indicated by the $\overset{(*)}{=}$ sign)
  using~\eqref{eq:npi1} and appropriate changes of variables.

  For~\eqref{eq:trans3}, we apply~\eqref{eq:nspi1} to the factor
  $s_j \pi_{r_j,s_j}^{(1)}(m_j)$ in the left-hand side (as only the
  tuples with $\epsilon_j=1$ contribute), and we apply~\eqref{eq:npi1}
  to the factor $\pi_{r_1,s_1}^{(1)}(m_1)$ in the right-hand
  side. Appropriate changes of variables then show that the difference
  vanishes.
\end{proof}

\begin{proof}[Proof of Proposition~\ref{prop:Pifin}]
Note that the left-hand sides of \eqref{eq:trans21}-\eqref{eq:trans24} all appear in  \eqref{eq:Pirew}, up to $\epsilon_3$ or $\bar{\epsilon}_3$ prefactors.
Changing them into the corresponding right-hand sides, and using $\epsilon_3+\bar{\epsilon}_3=1$, allows to write \eqref{eq:Pirew} as 
  \begin{equation}
    \label{eq:Pirew2}
    \begin{split}
      \Pi(m_1,\ldots,m_n) =&
      \sum_{I^{(-1,1)}_{s\geq 1}} r! (s-1)! \left(\sum_{j=3}^n \epsilon_j s_j \right) \prod_{i=1}^n \pi_{r_i,s_i}^{(\epsilon_i)}(m_i) \\
     +&
     \sum_{I^{(0,0)}_{s\geq 1}} r! (s-1)! \left(s_1+ s_2 + \sum_{j=3}^n \epsilon_j s_j \right) \prod_{i=1}^n \pi_{r_i,s_i}^{(\epsilon_i)}(m_i) \\
      +&
      \sum_{I^{(0,1)}_{r\geq 1} \cup I^{(1,0)}_{r\geq 1}} (r-1)! s! \left( \sum_{j=1}^n \bar\epsilon_j r_j \right) \prod_{i=1}^n \pi_{r_i,s_i}^{(\epsilon_i)}(m_i) \\
      +& \sum_{I^{(-1,1)}_{s=0} \cup I^{(0,0)}_{s=0} \cup
        I^{(0,1)}_{r=0} \cup I^{(1,0)}_{r=0}} r! s! \prod_{i=1}^n
      \pi_{r_i,s_i}^{(\epsilon_i)}(m_i).
    \end{split}
  \end{equation}
  By the transmutation relation~\eqref{eq:trans3}, we may replace the
  sum over $I^{(-1,1)}_{s\geq 1}$ in the first line by a sum over
  $I^{(1,1)}_{r\geq 1}$ in the third line. This gives a sum over
  $I^{(0,1)}_{r\geq 1} \cup I^{(1,0)}_{r\geq 1} \cup I^{(1,1)}_{r\geq 1}$ which we
  can rewrite as a sum over $I_{r\geq 1}$ minus a sum over
  $I^{(0,0)}_{r\geq 1}$. We claim that this latter sum will almost cancel
  the sum over $I^{(0,0)}_{s\geq 1}$ in the second line. Indeed,
  by writing
  $s_1+ s_2 + \sum_{j=3}^n \epsilon_j
  s_i = s - \sum_{j=3}^n \bar\epsilon_j s_j$ in the sum over
  $I^{(0,0)}_{s\geq 1}$, and
  $\sum_{j=1}^n \bar\epsilon_j r_j = r - \sum_{j=3}^n \epsilon_j r_j$
  in the sum over $I^{(0,0)}_{r\geq 1}$, we see using the transmutation
  relation~\eqref{eq:trans1} that their difference evaluates to
  \begin{equation}
  \label{eq:diffid}
    \left( \sum_{I^{(0,0)}_{s\geq 1}} - \sum_{I^{(0,0)}_{r\geq 1}} \right) r! s! \prod_{i=1}^n
    \pi_{r_i,s_i}^{(\epsilon_i)}(m_i) =
    \left( \sum_{I^{(0,0)}_{r=0}} - \sum_{I^{(0,0)}_{s=0}} \right) r! s! \prod_{i=1}^n
    \pi_{r_i,s_i}^{(\epsilon_i)}(m_i).
  \end{equation}
  Observe that the sum over $I^{(0,0)}_{s=0}$ precisely cancels the
  one appearing in the last line of~\eqref{eq:Pirew2}, and the sums
  over $I^{(0,0)}_{r=0}$, $I^{(0,1)}_{r=0}$ and $I^{(1,0)}_{r=0}$
  combine to form a sum over $I_{r=0}$. We arrive at the expression
  \begin{multline}
    \label{eq:Piinter}
      \Pi(m_1,\ldots,m_n) = \\ \sum_{I
      _{r\geq 1}} (r-1)! s! \left( \sum_{j=1}^n \bar\epsilon_j r_j \right) \prod_{i=1}^n \pi_{r_i,s_i}^{(\epsilon_i)}(m_i) \quad + \sum_{I^{(-1,1)}_{s=0} \cup I_{r=0}} r! s! \prod_{i=1}^n
    \pi_{r_i,s_i}^{(\epsilon_i)}(m_i)
  \end{multline}
  which is interesting on its own, since it is already symmetric in
  $m_1,\ldots,m_n$, see again the remark below~\eqref{eq:Pifin}. To
  obtain the wanted final expression, we write
  $\sum_{j=1}^n \bar\epsilon_j r_j = r - \sum_{j=1}^n \epsilon_j r_j$
  in the sum over $I_{r\geq 1}$ and apply again the transmutation
  relation~\eqref{eq:trans1} and the identity \eqref{eq:diffid},
  changing the sum over $I_{r\geq 1}$ and that over $I_{r=0}$ into the
  sum over $I_{s\geq 1}$ and that over $I_{s=0}$ of~\eqref{eq:Pifin},
  leaving the sum over $I^{(-1,1)}_{s=0}$ unchanged.
\end{proof}


\section*{Acknowledgements}

We thank Axel Bacher and Gilles Schaeffer for valuable discussions. JB
acknowledges the hospitality of the Laboratoire de Physique of ENS de
Lyon, where part of this work was completed.


%
%

\bibliographystyle{alphaurl}
\bibliography{polytightmaps}

\end{document}